\documentclass[review,hidelinks,onefignum,onetabnum]{siamart220329}

\usepackage{lipsum}
\usepackage{amsfonts}
\usepackage{graphicx}
\usepackage{epstopdf}
\usepackage{algorithmic}
\usepackage{xr-hyper}
\usepackage{cleveref}
\usepackage{nicefrac}
\usepackage{bm}
\usepackage{bbm}
\usepackage{amssymb}
\usepackage{amsmath}
\usepackage{mathabx}
\usepackage{mathrsfs}
\usepackage{mathtools}
\usepackage{empheq}
\usepackage{color}
\usepackage[normalem]{ulem}
\usepackage{stmaryrd}

\ifpdf
  \DeclareGraphicsExtensions{.eps,.pdf,.png,.jpg}
\else
  \DeclareGraphicsExtensions{.eps}
\fi

\def\softd{{\leavevmode\setbox1=\hbox{d}%
          \hbox to 1.05\wd1{d\kern-0.4ex{\char039}\hss}}}

\newcommand\eref[1]{(\ref{#1})}

\newcommand{\bla}[1]{\textcolor{blue}{#1}}

\newcommand{\jph}{{j+\frac{1}{2}}}
\newcommand{\jmh}{{j-\frac{1}{2}}}
\newcommand{\dx}{\Delta x}

\newcommand{\dt}{\Delta t}

\newcommand*\xbar[1]{%
  \hbox{%
    \vbox{%
      \hrule height 0.5pt 
      \kern0.4ex
      \hbox{%
        \kern-0.05em
        \ensuremath{#1}%
        \kern-0.00em
      }%
    }%
  }%
}

\newsiamremark{rmk}{Remark}
\newsiamthm{thm}{Theorem}
\newsiamremark{hypothesis}{Hypothesis}
\crefname{hypothesis}{Hypothesis}{Hypotheses}
\newsiamthm{claim}{Claim}
\newsiamthm{prop}{Proposition}

\headers{Well-Balanced Schemes Combining Different Forms}{R. Abgrall \& Y. Liu}

\graphicspath{{figures/}{./}}

\title{A New Approach for Designing Well-Balanced Schemes for the Shallow Water Equations: A Combination of Conservative and Primitive Formulations\thanks{Submitted to the editors DATE.
\funding{The work of Y. Liu was supported in part by SNSF grant 200020\_204917 and FZEB-0-166980.}}}

\author{R\'{e}mi Abgrall\thanks{Institute of Mathematics, University of Z\"{u}rich, 8057 Z\"{u}rich, Switzerland
  (\email{remi.abgrall@math.uzh.ch}).}
\and Yongle Liu\thanks{Institute of Mathematics, University of Z\"{u}rich, 8057 Z\"{u}rich, Switzerland
  (\email{yongle.liu@math.uzh.ch; liuyl2017@mail.sustech.edu.cn}).}}

\usepackage{amsopn}

\ifpdf
\hypersetup{
  pdftitle={A New Approach for Designing Well-Balanced Schemes for the Shallow Water Equations: A Combination of Conservative and Primitive Formulations},
  pdfauthor={R. Abgrall and Y. Liu}
}
\fi

\nolinenumbers

\begin{document}

\maketitle

\begin{abstract}
In this paper, we introduce a novel approach for constructing robust, fully well-balanced numerical methods for the one-dimensional Saint-Venant system, both with and without the Manning friction term. Inspired by the method presented in [R. Abgrall, Commun. Appl. Math. Comput. 5 (2023), pp. 370-402], we first combine the conservative and primitive formulations of the studied hyperbolic system in a natural way. The solution is globally continuous and described by a combination of point and average values. The point and average values will then be evolved by two different forms of the studied Saint-Venant system and we will show how to deal with the yield conservative and primitive forms in a well-balanced manner. The developed schemes are capable of exactly preserving both the still-water and moving-water equilibria. We demonstrate the behavior of the proposed new scheme on several challenging examples.
\end{abstract}



\section{Introduction}
The paper focus on the one-dimensional (1-D) Saint-Venant system of shallow water equations with non-flat bottom topography and Manning friction term, which is widely used to model the water flows in rivers, canals, lakes, reservoirs, and coastal areas. The studied system reads as
\begin{equation}\label{1.1}
  \left\{\begin{aligned}
  &h_t+q_x=0,\\
  &q_t+\left(hu^2+\frac{g}{2}h^2\right)_x=-ghZ_x-\frac{g n^2}{h^{\frac{7}{3}}}|q|q,
  \end{aligned}\right.
\end{equation}
where $g$ is the constant acceleration due to gravity, $x$ is the spatial variable, and $t$ is the time. $h(x,t)$ is the water depth, $u(x,t)$ is the velocity, and $q(x,t):=h(x,t)u(x,t)$ is the discharge. $Z(x)$ is the bottom elevation and $n$ is the Manning coefficient.

The system \eref{1.1} is a hyperbolic system of balance laws. Developing robust and accurate numerical methods for such a system is a challenging task as this nonlinear hyperbolic system admits both smooth and nonsmooth solutions. In addition, a good numerical scheme should be able to accurately respect a delicate balance between the flux and source term in \eref{1.1}. In particular, one is interested in developing well-balanced (WB) numerical schemes with the ability to preserve certain steady states at the discrete level. If a stable but non-well-balanced scheme is used, then the numerical errors may trigger the appearance of artificial waves of magnitude larger than the size of the waves to be captured. In order to ameliorate this situation, one can do mesh refinement so that the numerical errors will be controlled by the truncation errors. However, such mesh refinement may be too computationally expensive or even unaffordable in practical applications. Therefore, one needs to develop WB schemes.

Steady states are of great practical importance since they minimize, if they are stable, the energy, and many physically relevant solutions of \eref{1.1} are, in fact, small perturbations of these steady states. The smooth steady-state solutions of \eref{1.1} satisfy the following time-independent system:
\begin{equation}\label{1.2}
  \left\{\begin{aligned}
  & q_x=0,\\
  & \big(hu^2+\frac{g}{2}h^2\big)_x=-ghZ_x-\frac{gn^2}{h^{\frac{7}{3}}}|q|q,
  \end{aligned}\right.
\end{equation}
which can be integrated to obtain
\begin{equation}\label{1.3}
  q\equiv\widehat{q}={\rm Const}, \quad {E}:=\frac{u^2}{2}+g(h+Z+Q)\equiv\widehat{E}={\rm Const},
\end{equation}
where
\begin{equation}\label{1.4}
  Q=\int\limits_{\widehat{x}}^{x} \frac{n^2}{h^{\frac{10}{3}}}|q|q\, {\rm d}\xi
\end{equation}
is a global quantity with $\widehat{x}$ being an arbitrary number. In the case when $u\equiv0$, the steady states \eref{1.3} reduce to a particular class of steady states known as still-water equilibria (``lake at rest''):
\begin{equation}\label{1.5}
  u\equiv0,\quad w:=h+Z\equiv{\rm Const},
\end{equation}
where $w$ is the water surface. In fact, still-water steady state \eref{1.5} is a special case of the general moving-water steady state \eref{1.3}. If the bottom friction is neglected, that is, if $n=0$, the general smooth moving-water equilibria are given by
\begin{equation}\label{1.6}
  q\equiv\widehat{q}={\rm Const}, \quad E:=\frac{u^2}{2}+g(h+Z)\equiv\widehat{E}={\rm Const}.
\end{equation}

In the past two decades, many WB methods for the system \eref{1.1} have been developed. Some of them are capable of preserving still-water steady states \eref{1.5} only (see, e.g.,
\cite{ABBKP,BCKN,BNL,FMT11,JW,KPshw,Ric15,RB09,LeV98,GPC,TMK,LNK_2007,EPD_08,CSDT,XZS,KL,BerthonF,ToumaK,Audusse15,Puppo16,Dumbser19,Guo2023,Zhao2023,Yan2023,BB}), while others can preserve moving-water equilibria as well (see, e.g., \cite{BC_FWB, Victor2016,Victor2017,CPMP07,CLPwb13,CCHKW_18,CCKW,CKLLW,NXS,XS14,CK16,LLTX,CKLX,CTR,Xing14,MBerthon24,CMMOT,XS_AWENO,XS_DGGL,MRA_WB,Bueno,Zhang2023,BerthonM_23}). One may notice that it is significantly more difficult to obtain WB schemes for moving-water equilibria, since a special way to compute the numerical fluxes and a WB quadrature rule of the source term is essential to design such schemes. A key step in constructing a WB evaluation of the numerical fluxes, as demonstrated in the pioneering work \cite{NXS},  involves transforming the conservative variables $(h, q)^\top$ into the equilibrium variables $(q, E)^\top$, and vice versa. We focus on illustrating the key idea presented in certain second-order finite volume methods, high-order methods will be more involved. Given the cell averages of $(h, q)^\top$, one can compute the point values of $(q, E)^\top$ according to the equilibrium variables' definition in\eref{1.3}. When the solution $(h, q)^\top$ is at the steady state, the equilibrium variables $(q, E)^\top$ remain constants. Then, computing numerical fluxes necessitates calculating the one-sided point values of $(h, q)^\top$ at the cell interfaces. To make the scheme fully WB, the (generalized) hydrostatic reconstruction was employed. Rather than directly reconstructing the conservative variables, reconstruction is performed on the equilibrium variables to obtain the one-sided cell boundary point values of $(q, E)^\top$. Next, the one-sided cell boundary point values of the conservative variables can be recovered from the point values of the equilibrium variables. At this stage, one usually has to solve nonlinear equations, as it can be clearly seen from \eref{1.3} or \eref{1.6} that the transformation is nonlinear. In order to recover the point values of conservative variables from the nonlinear transformation and reconstructed point values of equilibrium variables, a nontrivial root-finding mechanism should be carefully developed. In addition, the solution of nonlinear equation is usually nonunique, one has to design a proper principle to single out the unique physical-relevant solution. Finally, WB evaluations of the numerical flux and approximation of the source term are computed using the obtained reconstructed point values of conservative variables.

Very recently, the author in \cite{Abgrall_camc} proposed a new class of schemes that can combine several writings of the same hyperbolic problem and showed how to deal with both the conservative and non-conservative forms in a natural way. In these schemes, the solution is globally continuous and described by a combination of point values and average values. The author demonstrated that one can recover the conservation property from how the average is updated and has more flexibility to deal directly with the primitive form of the system. In \cite{Abgrall_camc}, the author proved that this new class of schemes satisfies a Lax-Wendroff-like theorem and also provided numerical evidence through several challenging numerical experiments on the Euler equations to show that this new scheme possesses the capability to converge towards the correct weak solution.

In this paper, inspired by the method introduced in \cite{Abgrall_camc}, we develop a novel way for designing fully WB schemes, which exactly preserve all the steady
state solutions, both at rest and moving steady states, of the one-dimensional shallow water model. In addition to the conservative system \eref{1.1}, we simultaneously study the following primitive formulation:
\begin{equation}\label{1.7}
\left\{\begin{aligned}
&h_t+hu_x+uh_x=0,\\
&u_t+uu_x+g(h+Z)_x=-\frac{gn^2}{h^{\frac{4}{3}}}|u|u,\\
  \end{aligned}\right.
\end{equation}
which can be rewritten as
\begin{equation}\label{1.7a}
\left\{\begin{aligned}
&h_t+q_x=0,\\
&u_t+E_x=0,\\
  \end{aligned}\right.
\end{equation}
by employing the definition of equilibrium variables $q$ and $E$ in \eref{1.3}.
Since we have more flexibility to work on the primitive form, unlike most of the existing works mentioned above, this new class of WB schemes does not require a nontrivial root-finding mechanism. More precisely, we design a WB finite-volume method for system \eref{1.1} to evolve the cell averages of conservative variables $(h, q)^\top$ and simultaneously propose a WB residual distribution (RD) method for system \eref{1.7a} to evolve the point values of primitive variables $(h,u)^\top$. Given the average and point values of conservative variables (the point values of conservative variables can be obtained from the point values of the primitive variables), we can apply a third-order parabolic interpolant on the conservative variables and then compute the point values of conservative variables at the cell center. Having these point values of conservative variables allows us to compute the point values of equilibrium variables at the cell center and interfaces. Next, we do the spatial discretizations on the equilibrium variables, which helps to ensure the WB property. It follows from \eref{1.3} that the approximations of spatial derivatives appearing in \eref{1.7a} will vanish when the solution is at the steady sate. Unlike many existing methods, at this stage, we do not rely on any special root-finding mechanism. Instead, we compute the fluctuation of the system \eref{1.7a} directly based on the equilibrium variables, and split the residual onto the nodes through an upwind consideration. For numerically solving \eref{1.1}, we simply compute the flux using the point values of primitive variables at the cell interfaces and a WB approximation of source term is given using the high-order extrapolation technique and a moving-water preserving correction term. In addition to the WB property, the developed schemes are positivity-preserving. This attributes to the fact that we have applied a MOOD paradigm from \cite{CDL,Vilar} equipped with a first-order scheme that has a Local Lax-Friedrichs' flavour to the solvers of \eref{1.1} and \eref{1.7a}.

The rest of the paper is organized as follows. In \cref{sec2}, we introduce the new WB positive-preserving schemes for the conservative formulation \eref{1.1} and the primitive formulation \eref{1.7a}. In \cref{sec3}, we present the numerical results, which demonstrate the WB property, positivity-preserving property, high resolution, and robustness of the proposed schemes. Finally, in \cref{sec4}, we give some concluding remarks and comments on the future development and applications of the proposed WB schemes.

\section{Schemes}\label{sec2}
In this section, we describe the third-order semi-discrete numerical schemes for the 1-D Saint-Venant system \eref{1.1} and its primitve version \eref{1.7a}. To this end, we first rewrite \eref{1.1} in a vector form as
\begin{equation}\label{2.1}
  \bm U_t+\bm F(\bm U)_x=\bm S(\bm U), \quad \bm U=(h, q)^\top,
\end{equation}
where $\bm U\in \mathbb{R}^d$ and
\begin{equation}\label{2.2}
   \bm F(\bm U)=\begin{pmatrix}q\\hu^2+\frac{g}{2}h^2\end{pmatrix},\quad \bm S(\bm U)=\begin{pmatrix}0\\-ghZ_x-ghQ_x\end{pmatrix}.
\end{equation}
Simultaneously, we also rewrite \eref{1.7a} in the following vector form:
\begin{equation}\label{2.3}
  \bm V_t+\frac{\partial \bm E}{\partial x}=\bm 0, \quad \bm V=(h, u)^\top, \quad \bm E=(q, E)^\top,
\end{equation}
where $\bm V=\Psi(\bm U)$ with
\begin{equation*}
  \Psi: \mathcal{D}_{\bm U}\rightarrow\mathcal{D}_{\bm V}, ~ \mathcal{D}_{\bm U}=\{\bm U=(h,q)^\top\in\mathbb{R}^2: h\geq0\},~ \mathcal{D}_{\bm V}=\{\bm V=(h,u)^\top\in\mathbb{R}^2: h\geq0\}.
\end{equation*}

We divide the computational domain into a set of uniform (for simplicity of presentation) cells $K_\jph:=[x_j, x_{j+1}]$ of size $\dx=x_{j+1}-x_j$, which are centered at $x_\jph=(x_j+x_{j+1})/2$. We assume that at a certain time level $t$, the numerical solution of \eref{2.1}--\eref{2.2}, realized in terms of its cell average  $\,\xbar{\bm U}_\jph(t)\approx\frac{1}{\dx}\int_{K_\jph}\bm U(x,t)\,{\rm d}x$, is available. We denote the numerical solutions of \eref{2.3} by $\bm V_j(t)\approx\bm V(x_j,t)$, which are also known at time $t$. Following the method introduced in \cite{Abgrall_camc}, we evolve the cell averages $\xbar{\bm U}_\jph$ in time by solving the following system of time-dependent ODEs:
\begin{equation}\label{2.5}
  \frac{\rm d}{{\rm d}t}\,\xbar{\bm U}_\jph=-\frac{\bm{\mathcal{F}}_{j+1}-\bm{\mathcal{F}}_j}{\dx}+\xbar{\bm S}_\jph,
\end{equation}
where
\begin{equation}\label{2.6}
   \bm{\mathcal F}_j=\bm F(\bm U_j)\stackrel{\eref{2.2}}{=}\begin{pmatrix}q_j\\h_ju_j^2+\frac{g}{2}h_j^2\end{pmatrix}
\end{equation}
and
\begin{equation}\label{2.7}
  \xbar{\bm S}_\jph:\approx\frac{1}{\dx}\int\limits_{K_\jph}{\bm S}(\bm U(x,t)){\rm d}x
\end{equation}
is an approximation of the cell average of source term. It is noteworthy that the evolution of the cell average in \eref{2.5} is obtained as usual by integrating the conservative formulation given by \eref{2.1} over the cell and over a time step, and by applying Gauss' theorem. The flux required at the cell interface can be immediately evaluated using the available point values. Moreover, the source term average \eref{2.7} can be approximated using some quadrature rules. At the same time, the point values $\bm V_j$ are also evolved by solving a system of time-dependent ODEs given by
\begin{equation}\label{2.8}
  \frac{\rm d}{{\rm d}t}\,{\bm V}_j=-(\overleftarrow{\Phi}_\jph^{\bm V}+\overrightarrow{\Phi}_\jmh^{\bm V}),
\end{equation}
where $\overleftarrow{\Phi}_\jph^{\bm V}+\overrightarrow{\Phi}_\jmh^{\bm V}$ is a consistent approximation of $\frac{\partial \bm E}{\partial x}(x_j)$. The update of the point values allows for more flexibility as is not subject to the constraint of conservation. The primary restriction in computing $\overleftarrow{\Phi}_\jph^{\bm V}$ and $\overrightarrow{\Phi}_\jmh^{\bm V}$ is no more than ensuring the stability of the resulting method, thus some form of upwinding is incorporated. Notice that all of the indexed quantities in \eref{2.5}--\eref{2.7}, and \eref{2.8} are time-dependent, but from here on, we omit this dependence for the sake of brevity.

In general, the residuals $\overleftarrow{\Phi}_\jph^{\bm V}$ and $\overrightarrow{\Phi}_\jmh^{\bm V}$ need to depend on some ${\bm E}_\ell$ and ${\bm E}_{\ell+\frac{1}{2}}$, and thus on some ${\bm V}_\ell$ and ${\bm V}_{\ell+\frac{1}{2}}\approx{\bm V}(x_{\ell+\frac{1}{2}})$. Hence, it is necessary to first recover the missing information at the half points. Given the point values $\bm U_j=\Psi(\bm V_j)$, $\bm U_{j+1}=\Psi(\bm V_{j+1})$, and the cell average $\xbar{\bm U}_\jph$, we begin with constructing a continuous function represented by $\bm U_j$, $\bm U_{j+1}$, and $\xbar{\bm U}_\jph$. This function, found in \eref{2.12} from \cref{appb} and denoted by $\widetilde{\bm U}(x)$, is a third-order parabolic interpolant of $\bm U$ defined on $K_\jph$. Subsequently, we compute the point value $\bm U_\jph$ at the cell center, given by
\begin{equation}\label{2.15x}
  \bm U_\jph:=\widetilde{\bm U}(x_\jph)=\frac{3}{2}\,\xbar{\bm U}_\jph-\frac{1}{4}(\bm U_j+\bm U_{j+1}).
\end{equation} 
Note that this relation is simply $\xbar{\bm U}_\jph=\frac{1}{6}(\bm U_j+\bm U_{j+1}+4\bm U_\jph)$, which is essentially Simpson's formula, ensuring third-order accuracy.

We proceed with computing the point values of $Q$. To this end, we rewrite the global quantity $Q$ defined in \eref{1.4} in a recursive way and approximate the yield integration part by the Simpson's rule. Namely, we take  
\begin{equation}\label{2.16Q1}
\begin{aligned}
Q_\jph&=Q_j+\int\limits_{x_j}^{x_{\jph}} \frac{n^2}{h^{\frac{10}{3}}}|q|q\, {\rm d}x\\
&\approx Q_j+\frac{ n^2\dx}{12}\left(\frac{|q_j|q_j}{(h_j)^{10/3}}+\frac{|q_\jph|q_\jph}{(h_\jph)^{10/3}}+4\frac{|q_{j+\frac{1}{4}}|q_{j+\frac{1}{4}}}{(h_{j+\frac{1}{4}})^{10/3}}\right)
  \end{aligned}
\end{equation}
and
\begin{equation}\label{2.16Q2}
\begin{aligned}
  Q_{j+1}&=Q_j+\int\limits_{x_j}^{x_{j+1}} \frac{n^2}{h^{\frac{10}{3}}}|q|q\, {\rm d}x\\
  &\approx Q_j+\frac{ n^2\dx}{6}\left(\frac{|q_j|q_j}{(h_j)^{10/3}}+\frac{|q_{j+1}|q_{j+1}}{(h_{j+1})^{10/3}}+4\frac{|q_\jph|q_\jph}{(h_\jph)^{10/3}}\right),
  \end{aligned}
\end{equation}
where $\bm U_\jph$ is given in \eref{2.15x}, $\bm U_{j+\frac{1}{4}}$ can be computed from \eref{2.12}:
\begin{equation}\label{upv}
  \bm U_{j+\frac{1}{4}}=\frac{3}{16}\bm U_j+\frac{18}{16}\xbar{\bm U}_\jph-\frac{5}{16}\bm U_{j+1},
\end{equation}
and $Q_1=0$ as we simply take $\widehat{x}=x_1$ in \eref{1.4}.

We then compute the point values $\bm E_j=(q_j,E_j)^\top$ and ${\bm E}_\jph=(q_\jph, E_\jph)^\top$ for all $j$ as follows:
\begin{equation}\label{2.16a}
\begin{aligned}
 &q_j:=h_ju_j,\quad  E_j=\frac{u_j^2}{2}+g(h_j+Z_j+Q_j),\\
 &E_\jph=\frac{u_\jph^2}{2}+g(h_\jph+Z_\jph+Q_\jph).
 \end{aligned}
\end{equation}
In \eref{2.16a}, $Z_j=Z(x_j)$, $Z_\jph=Z(x_\jph)$, and $u_\jph:=\mathfrak{D}(q_\jph,h_\jph)$ with $\mathfrak{D}(a,b)$ being the desingularization function defined by \eref{Desing} in \cref{appc}. Equipped with these obtained point values, we use the discretization provided in \cite[\S 3]{Abgrall_camc} to obtain:
\begin{equation}\label{2.17}
    \delta_j^+ \bm E=\frac{1}{\dx}\big(\bm E_{j-1}-4\bm E_\jmh+3\bm E_j\big),\quad \delta_j^- \bm E=\frac{1}{\dx}\big(-3\bm E_j+4\bm E_\jph-\bm E_{j+1}\big),
\end{equation}
which will be used to compute the residuals in \eref{2.8}. Instead of employing the Jacobian splitting method as described in \cite{Abgrall_camc}, we utilize flux vector splitting to compute 
\begin{equation}\label{rnew}
\overleftarrow{\Phi}_\jph^{\bm V}=(\widetilde J(\bm V_j))^-\delta_j^-{\bm E},\quad \overrightarrow{\Phi}_\jmh^{\bm V}=(\widetilde J(\bm V_j))^+\delta_j^+{\bm E},
\end{equation}
where $(\widetilde J(\bm V_j))^\pm$ are defined by
\begin{equation}\label{2.10a}
\widetilde J(\bm V_j)^\pm=\Upsilon_j\begin{pmatrix}\frac{\lambda_1^\pm}{\lambda_1}&0\\0&\frac{\lambda_2^\pm}{\lambda_2}\\\end{pmatrix}\Upsilon_j^{-1}
  :=\Upsilon_j\widetilde{\Lambda}^\pm\Upsilon_j^{-1}.
\end{equation}
In \eref{2.10a}, $\lambda_{1,2}=u_j\mp\sqrt{gh_j}$ are the corresponding eigenvalues of the Jacobian,  $\Upsilon_j=\Upsilon(\bm V_j)$ is the corresponding eigenvector matrix given by
\begin{equation*}
  \Upsilon_j=\begin{pmatrix}-\sqrt{\frac{h_j}{g}}&\sqrt{\frac{h_j}{g}}\\1&1\end{pmatrix},
\end{equation*} 
and when $\widetilde{\Lambda}^\pm$ is computed we have used Haten-Yee entropy fix to prevent calculation on small eigenvalues. We would like to emphasize that the adoption of flux vector splitting in \eref{rnew} not only meets the necessity to upwind but also facilitates achieving well-balanced, as $\delta_j^\pm{\bm E}$ vanishes at the steady state.

So far, we have completed the computation of all quantities required for the evolution of \eref{2.8}; see \eref{2.17}--\eref{2.10a}. We finally describe a WB discretization of $\xbar{\bm S}_j$ so that the computation of the necessary quantities for the evolution of \eref{2.5} is also closed. In fact, since $S^{(1)}\equiv0$, we only need to design a WB approximation for the second component:
\begin{equation}\label{2.20}
  \xbar{S}_\jph^{(2)}\approx\frac{1}{\dx}\int\limits_{K_\jph}S^{(2)}{\rm d}x=\frac{1}{\dx}\int\limits_{K_\jph}\left(-ghZ_x-ghQ_x\right){\rm d}x.
\end{equation}
This can be done by applying a WB fourth-order extrapolation technique derived in \cite{NXS} within the context of WENO schemes. Specially, we compute 
\begin{equation}\label{2.23}
  \xbar{S}_\jph^{(2)}:=\frac{4S_2-S_1}{3\dx}
\end{equation}
with $S_1$ and $S_2$ are defined as
\begin{equation}\label{S12}
\begin{aligned}
  S_1&=-g\overline{\overline h}_{j,j+1}[Z+Q]_{j,j+1}+\widehat{s}_{j,j+1}^{\rm int},\\
  S_2&=-g\overline{\overline h}_{j,\jph}[Z+Q]_{j,\jph}+\widehat{s}_{j,\jph}^{\rm int}-g\overline{\overline h}_{\jph,j+1}[Z+Q]_{\jph,j+1}+\widehat{s}_{\jph,j+1}^{\rm int}.
  \end{aligned}
\end{equation}
Here, $\overline{\overline h}_{i_1,i_2}=\frac{1}{2}(h_{i_1}+h_{i_2})$, $[Z+Q]_{i_1,i_2}=Z_{i_2}+Q_{i_2}-(Z_{i_1}+Q_{i_1})$, and $\widehat{s}^{\rm int}_{i_1,i_2}$ is the moving-water correction term defined as
\begin{equation}\label{sint}
  \widehat{s}^{\rm int}_{i_1,i_2}=\frac{\Gamma_{i_1,i_2}}{4}[h]_{i_1,i_2}\big([u]_{i_1,i_2}\big)^2,
\end{equation} 
where $[h]_{i_1,i_2}=h_{i_2}-h_{i_1}$, $[u]_{i_1,i_2}=u_{i_2}-u_{i_1}$, and $\Gamma_{i_1,i_2}\in[0,1]$ is a parameter defined by 
\begin{equation}\label{Gamma}
  \Gamma_{i_1,i_2}=\mathcal{H}(\eta_{i_1,i_2}),\quad \eta_{i_1,i_2}=\Vert \bm E_{i_2}-\bm E_{i_1} \Vert,
\end{equation}
with $\mathcal{H}(\eta)$ being a smooth cutoff function given by
\begin{equation}\label{2.43ah}
  \mathcal{H}(\eta):=\frac{1}{1+(C\eta)^m}.
\end{equation}
Note that when the computed solution is locally (almost) at/near a steady state ($\eta\approx 0$), $\Gamma_{i_1,i_2}$ is $1$. Therefore, the correction term $\frac{1}{4}[h]_{i_1,i_2}[u]^2_{i_1,i_2}$, crucial for preserving general moving-water equilibria \cite{NXS}, becomes active to guarantee the WB property. Conversely, when the solution deviates from the steady state, $\Gamma_{i_1,i_2}$ rapidly decreases to $0$ and thus the correction term is excluded. This mechanism aids in recovering consistency, particularly in scenarios involving flat bottom topography and the absence of Manning friction. In such cases, where the actual source term $-ghZ_x-ghQ_x$ vanishes, the additional correction term is necessary to maintain consistency with zero. Alternative techniques for ensuring the consistency of the correction term with zero are also available and can be found, e.g., \cite{Victor2016, Victor2017, DM_FWBC, BMLS, NXS, BC_FWB}. In all of the numerical experiments reported in \cref{sec3}, we have used the constants $C=2$ and $m=10$.

We now prove that the numerical schemes given by \eref{2.5}, \eref{2.6}, \eref{2.23}--\eref{sint} and \eref{2.8}, \eref{rnew} are WB in the sense that the following theorem holds.
\begin{thm}[Well-balanced for general steady state]\label{thm2}
  If the solution $\bm U$ is in equilibrium \eref{1.3} at all Gauss-Lobatto points $\{x_j,x_\jph,x_{j+1}\}$, $\forall K_\jph$, namely, the discrete data $\{\bm U_j, \bm U_\jph, \bm U_{j+1}\}$ satisfy the following relations:
  \begin{equation}\label{2.33}
    q_j=q_\jph=q_{j+1}\equiv\widehat{q}={\rm Const},\quad E_j=E_\jph=E_{j+1}\equiv\widehat{E}={\rm Const},\quad \forall K_\jph,
  \end{equation}
then, the proposed schemes given by \eref{2.5}, \eref{2.6}, \eref{2.23}--\eref{sint} and \eref{2.8}, \eref{rnew} are WB, provided that the right-hand sides of \eref{2.5} and \eref{2.8} all vanish, i.e.,
\begin{equation}\label{2.25}
  \frac{\rm d}{{\rm d}t}\,\xbar{\bm U}_\jph=\bm 0 \quad \mbox{and} \quad \frac{\rm d}{{\rm d}t}\,{\bm V}_j=\bm 0.
\end{equation}
\end{thm}
{{\bf Proof.}} On one hand, when the data satisfies \eref{2.33}, the equilibrium variable $\bm E=(\widehat q, \widehat{E})^\top$ and thus $\delta_j^\pm \bm E\equiv0$ follows from \eref{2.17}. Plugging these values into \eref{rnew}, we easily get $\overleftarrow{\Phi}_\jph^{\bm V}=\overrightarrow{\Phi}_\jmh^{\bm V}=\bm 0$, which implies that $\frac{\rm d}{{\rm d}t}\,{\bm V}_j=\bm 0$. 

On the other hand, at the steady state \eref{2.33}, following from \eref{Gamma} and \eref{2.43ah}, we have $\Gamma_{j,j+1}=1$, $\Gamma_{j,\jph}=1$, and $\Gamma_{\jph,j+1}=1$. Therefore, \eref{sint} becomes
\begin{equation}\label{sint2}
  \widehat{s}^{\rm int}_{i_1,i_2}=\frac{1}{4}[h]_{i_1,i_2}\big([u]_{i_1,i_2}\big)^2.
\end{equation} 
Moreover, using the definition of $E$ in \eref{1.3} together with the relations in \eref{2.33}, we can obtain
\begin{equation*}
\begin{aligned}
  &g[Z+Q]_{j,j+1}=-\frac{1}{2}[u^2]_{j,j+1}-g[h]_{j,j+1},\\
  &g[Z+Q]_{j,\jph}=-\frac{1}{2}[u^2]_{j,\jph}-g[h]_{j,\jph},\\
  &g[Z+Q]_{\jph,j+1}=-\frac{1}{2}[u^2]_{\jph,j+1}-g[h]_{\jph,j+1}
  \end{aligned}
\end{equation*}
Substituting these results along with \eref{sint2} into \eref{S12} and after some simple algebraic calculations, we have
\begin{equation*}
S_1=[hu^2]_{j,j+1}+\frac{g}{2}[h^2]_{j,j+1},\quad
S_2=[hu^2]_{j,j+1}+\frac{g}{2}[h^2]_{j,j+1},
\end{equation*}
which are obtained by also using the fact that $h_ju_j=h_\jph u_\jph=h_{j+1}u_{j+1}$ at the general steady state. Consequently, the second component source term approximation \eref{2.23} reduces to
\begin{equation}\label{2.33a}
  \xbar{S}_\jph^{(2)}=\frac{4S_2-S_1}{3\dx}=\frac{1}{\dx}\Big[h_{j+1}u_{j+1}^2-h_ju_j^2+\frac{g}{2}\big(h_{j+1}^2-h_j^2\big)\Big]
\end{equation}
at the general steady state \eref{2.33}. Finally, we compute the right-hand sides of \eref{2.5} in a component-wise manner:
\begin{equation*}
\begin{aligned}
  \frac{\rm d}{{\rm d}t}\,\xbar{U}^{(1)}_\jph&=-\frac{\mathcal{F}^{(1)}_{j+1}-\mathcal{F}^{(1)}_j}{\dx}=-\frac{q_{j+1}-q_j}{\dx}\stackrel{\eref{2.33}}{=}0,\\
  \frac{\rm d}{{\rm d}t}\,\xbar{U}^{(2)}_\jph&=-\frac{\mathcal{F}^{(2)}_{j+1}-\mathcal{F}^{(2)}_j}{\dx}+\xbar{S}^{(2)}_\jph\\
  &=-\frac{h_{j+1}u_{j+1}^2-h_ju_j^2+\frac{g}{2}(h_{j+1}^2-h_j^2)}{\dx}+\xbar{S}^{(2)}_\jph
 \stackrel{\eref{2.33a}}{=}0,
  \end{aligned}
\end{equation*}
which implies that $\frac{\rm d}{{\rm d}t}\,\xbar{\bm U}_\jph=\bm 0$.
 $\hfill\blacksquare$

\subsection{Positivity-Preserving}\label{sec22}
In addition to achieving the WB property, it is also crucial to preserve the positivity of water depth. This is not only necessary for physical significance, but also essential for theoretical analysis and numerical stability. For the proposed schemes, both the cell averages $\{\xbar{h}_\jph\}$ and the point values $\{h_j\}$ should keep positivity.

As pointed out in \cite{Abgrall_camc}, the proposed schemes is at most linearly stable and not positivity-preserving, with a CFL condition based on the fine grid. In order to ensure the positivity of the point values of water depth, we apply the MOOD paradigm which was used in \cite{AT,Abgrall_camc,CDL,Vilar} to switch the numerical schemes to the parachute ones. The idea is to work with several schemes ordered from the most accurate/less stable one to the low order/more reliable one, which is a parachute scheme and able to provide solutions fulfilling some criteria. These detection criteria are set up according to the computer, physical, and numerical admissibility. In all of the numerical examples reported in \cref{sec3}, we have used the following detection criteria to detect the troubled cells:
\begin{enumerate}
  \item We check if all of the primitive variables are numbers (i.e., not NaN and not Inf). If this holds, the criterion is set to .TRUE. and we look for the next criterion, else it is set to .FALSE. and we jump to the next element;
  \item We check if the water depth is positive. If not, we set the criterion to .FALSE. and jump to the next element. Otherwise, we set the criterion to .TRUE. and check next criterion.
  \item We check if the solution is not locally constant in the element $K_\jph$. We take $\varepsilon=\dx^3$ and check if
  \begin{equation*}
    |\max\limits_{\ell\in S} h_\ell-\min\limits_{\ell\in S} h_\ell|\leq \varepsilon,\quad S=\{\ell|\ell=j+\frac{1}{2}(i-3), 1\leq i\leq7, i\in \mathbb{Z}\}.
  \end{equation*}
  If we observe that the solution is locally constant, the .TRUE. criterion is not modified.
  \item We check if a new extrema in $h$ is created or not, by comparing with the solution at the previous time step (for the sake of convenience, we assume it corresponds $t=t^n$). Namely, we first compute $\min_\jph h$ (resp. $\max_\jph h$) the minimum (resp. maximum) of the water depth $h$ at $t=t^n$ on $K_\jmh$, $K_\jph$, and $K_{j+\frac{3}{2}}$. Then:
      \begin{enumerate}
             \item If $h_\ell^{n+1}\in[\min_\jph h-\delta,\max_\jph h+\delta],~~\forall \ell\in\{j,\jph,j+1\}$, we set the criterion to .TRUE. and jump to the next element. Here, $\delta$ is a relaxation parameter and we take $\delta=\min(10^{-4},10^{-3}|\max_\jph h-\min_\jph h|)$.

             \item Else, we denote the Lagrange interpolation polynomial that interpolates $h_\jph(x)$ by $\mathfrak{H}_\jph$ and 
        \begin{itemize}
          \item we compute $\zeta'=\mathfrak{H}_\jph'(x_\jph)$, $\zeta'_L=\mathfrak{H}_\jph'(x_j)$, and the minimum and maximum values of the derivative around $x_j$:
              \begin{equation*}
                \zeta'^{,j}_{\min/\max}=\min/\max(\mathfrak{H}'_\jph(x_\jph),\mathfrak{H}'_\jmh(x_\jmh)).
              \end{equation*}
           All of the obtained quantities are then used to define the left detection factor $\beta_L$:
           \begin{equation*}
             \beta_L=\left\{\begin{aligned}
             &\min\Big(1,\frac{\zeta'^{,j}_{\max}-\zeta'}{\zeta'_L-\zeta'}\Big),\quad&& \mbox{if}~\zeta'_L>\zeta',\\
             &1,\quad && \mbox{if}~\zeta'_L=\zeta',\\
             &\min\Big(1,\frac{\zeta'^{,j}_{\min}-\zeta'}{\zeta'_L-\zeta'}\Big),\quad&& \mbox{if}~\zeta'_L<\zeta'.
             \end{aligned}\right.
           \end{equation*}
          \item Analogously, we compute $\zeta'=\mathfrak{H}_\jph'(x_\jph)$, $\zeta'_R=\mathfrak{H}_\jph'(x_{j+1})$, the minimum and maximum values of the derivative around $x_{j+1}$:
           \begin{equation*}
               \zeta'^{,j+1}_{\min/\max}=\min/\max(\mathfrak{H}'_{j+\frac{3}{2}}(x_{j+\frac{3}{2}}),\mathfrak{H}'_\jph(x_\jph)),
              \end{equation*}
              and define the right detection factor $\beta_R$ as
               \begin{equation*}
             \beta_R=\left\{\begin{aligned}
             &\min\Big(1,\frac{\zeta'^{,j+1}_{\max}-\zeta'}{\zeta'_R-\zeta'}\Big),\quad&& \mbox{if}~\zeta'_R>\zeta',\\
             &1,\quad && \mbox{if}~\zeta'_R=\zeta',\\
             &\min\Big(1,\frac{\zeta'^{,j+1}_{\min}-\zeta'}{\zeta'_R-\zeta'}\Big),\quad&& \mbox{if}~\zeta'_R<\zeta'.
             \end{aligned}\right.
           \end{equation*}
          \item Next, we set $\beta=\min(\beta_L, \beta_R)$.
          \item Finally, if $\beta=1$, we will have a true extrema and thus we set the criterion to .TRUE. and jump to the next element. Otherwise, we set the criterion to .FALSE. and jump to the next element.
        \end{itemize}
        \end{enumerate}
\end{enumerate}
The idea behind the step 4 (b) is as described in \cite[Section 3.1]{Vilar}: we want to preserve scheme accuracy in the presence of smooth extrema and try to check if the gradient of the interpolation $\mathfrak{H}_\jph$ lies in the interval
$[\min(\zeta'^{,j}_{\min},\zeta'^{,j+1}_{\min}),\max(\zeta'^{,j}_{\max},\zeta'^{,j+1}_{\max})]$.

\begin{figure}[ht!]
\centerline{\includegraphics[trim=0.01cm 0.01cm 0.01cm 0.01cm,clip,width=12.5cm]{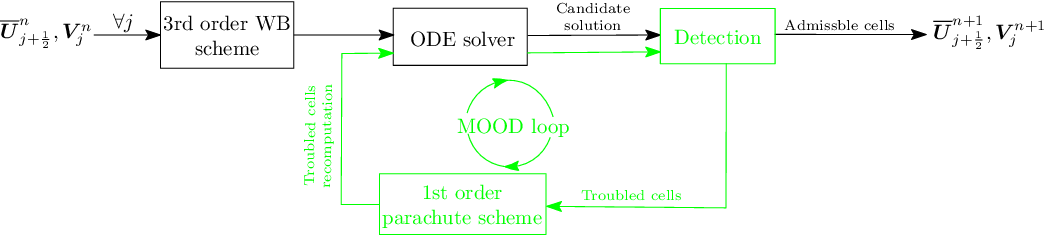}}
\caption{Sketch of the proposed WB active flux like scheme stabilized by a MOOD loop.\label{mood}}
\end{figure}

We expect to use as often as possible the third order scheme, and to use the lower order one to correct potential problems. Therefore, we first perform a step of the higher accurate scheme, and then check the criteria on each cell to detect the areas where the criteria are not met. Next, on these troubled cells, we switch to the parachute scheme, which is a slight modification of the first-order parachute scheme in \cite{Abgrall_camc} and has a Local Lax-Friedrichs' flavor. The scheme reads as
\begin{equation}\label{2.42b}
  \frac{{\rm d}}{{\rm d}t}\,{\bm V}_j=-(\overleftarrow{{\Phi}}_\jph+\overrightarrow{\Phi}_\jmh),
\end{equation}
where
\begin{equation}\label{2.42}
   \overleftarrow{\Phi}_\jph=\frac{1}{\dx}\begin{pmatrix}
                    \xbar{q}_\jph-q_j-\alpha_\jph(\xbar{h}_\jph-h_j) \\
                    \xbar{E}_\jph-E_j-\alpha_\jph(u_\jph-u_j)
                  \end{pmatrix}
\end{equation}
and
\begin{equation}\label{2.43}
 \overrightarrow{\Phi}_\jph=\frac{1}{\dx}\begin{pmatrix}
                    q_{j+1}-\xbar{q}_\jph-\alpha_\jph(\xbar{h}_\jph-h_{j+1}) \\
                    E_{j+1}-\xbar{E}_\jph-\alpha_\jph(u_\jph-u_{j+1})
                  \end{pmatrix}.
\end{equation}
In \eref{2.42} and \eref{2.43},  $u_\jph=\mathfrak{D}(\xbar{q}_\jph,\xbar{h}_\jph)$ and $\alpha_\jph$ is an upper-bound of the spectral radius of $J(V_j)$, $J(V_\jph)$, and $J(V_{j+1})$. $\xbar{E}_\jph$ is obtained as follows:
\begin{equation*}
  \xbar{E}_\jph\stackrel{\eref{1.3}}{=}\frac{(u_\jph)^2}{2}+g(\xbar{h}_\jph+Z_\jph)+\xbar{Q}_\jph, \quad \xbar{Q}_\jph=\frac{Q_j+Q_{j+1}}{2}.
\end{equation*}
As for preserving the positivity of cell averages $\{\xbar h_\jph\}$ for the conservative formulation \eref{2.5}, we will replace the evaluation $\bm {\mathcal{F}}_j$ with the Local Lax-Friedrichs' numerical flux in the troubled cells. The MOOD procedure is illustrated in Figure \ref{mood}. It is quite well-known that the Local Lax-Friedrichs' finite volume scheme for the evolution of the average value of water depth is positivity-preserving under the standard Courant--Friedrichs--Lewy (CFL) condition with a CFL number of 1. In the following theorem, we show that the first-order scheme \eref{2.42b} with \eref{2.42}--\eref{2.43} is also positivity-preserving for the evolution of the point value of water depth. 
\begin{theorem}
  Let $\xbar{h}_\jph$ and $h_j$ be the average and point values at the $n$-th time level, respectively. For the first-order scheme \eref{2.42b} with \eref{2.42}--\eref{2.43}, if $\xbar{h}_\jph\geq0$ and $h_j\geq0$, then $h_j^{n+1}\geq0$, for all $j$, under the CFL condition
  \begin{equation}\label{CFL}
    \dt\leq\frac{\dx}{2a_{\max}},\quad a_{\max}=\max_j\{\alpha_\jph\}.
  \end{equation}
\end{theorem}
Proof: For the sake of brevity, we consider only the forward Euler method in the discretization of the ODE system \eref{2.42b}. The resulting schemes will remain positivity-preserving when the forward Euler method is substituted with a three-stage third-order strong stability preserving Runge-Kutta solver in \cref{sec3}, as such a solver can be written as a convex combination of several forward Euler steps.

Let $\nu = \frac{\dt}{\dx}$ and we aim to prove $h_j^{n+1} \geq 0$. Substituting \eref{2.42}  and \eref{2.43} into \eref{2.42b}, we obtain:
\begin{equation*}
\begin{aligned}
  h_j^{n+1}&=h_j-\nu\Big[\xbar{q}_\jph-\xbar{q}_\jmh-\alpha_\jph\xbar{h}_\jph-\alpha_\jmh\xbar{h}_\jmh+(\alpha_\jph+\alpha_\jmh) h_j\Big]\\
  &=\Big[1-(\alpha_\jph+\alpha_\jmh)\nu\Big]h_j+\nu\Big[\xbar{h}_\jph(\alpha_\jph-{u}_\jph)+\xbar{h}_\jmh(\alpha_\jmh+{u}_\jmh)\Big]\\
  &\geq(1-2a_{\max}\nu)h_j\stackrel{\eref{CFL}}{\geq}0.
  \end{aligned}
\end{equation*} 
The penultimate  inequality follows from the fact that $\alpha_\jph\leq a_{\max}$ and 
\begin{equation*}
  \alpha_\jph=\max\Big(|u_\jph|+\sqrt{gh_\jph},|u_j|+\sqrt{gh_j},|u_{j+1}|+\sqrt{gh_{j+1}}\Big)\geq\pm u_\jph.
\end{equation*}
$\hfill\blacksquare$

\section{Numerical Examples}\label{sec3}
In this section, we demonstrate the performance of the proposed new WB schemes on several numerical examples. In some of the following examples, we compare the numerical results obtained by the developed new WB scheme (referred to as ``WB AF'' scheme) with those computed by the WB second-order central-upwind scheme (referred to as ``WB PCCU'' scheme) from \cite{CKLX}. In all numerical experiments, unless specified otherwise, we assume that the acceleration due to gravity $g=9.812$ and use zero-order extrapolation boundary conditions. In Examples 1--6 and 9, we take the manning coefficient $n=0$, while in Examples 7 and 8 we take $n=0.05$. The figures presented below plot the cell average values of the variables of interest.

In the numerical implementation, we integrate the ODE systems \eref{2.5} and either \eref{2.8} or \eref{2.42b} using the three-stage third-order strong stability preserving (SSP) Runge-Kutta method (see, e.g., \cite{GKS,GST}) with an adaptive time step computed at every time level using the CFL number $0.4$, namely, by taking
\begin{equation*}
\dt=0.4\frac{\dx}{a_{\max}},
\end{equation*}
where $a_{\max}$ is an upper-bound of the spectral radius of $J(\bm V_j)$, $J(\bm V_\jph)$, and $J(\bm V_{j+1})$.

\begin{rmk}
In this section, all numerical examples presented here exclusively involve continuous bottom topography. This is due to the fact that when the bottom topography is discontinuous, the presence of nonconservative product term $-ghZ_x$ makes it impossible to understand weak solutions in the sense of distributions (instead weak solutions can be understood as the Borel measures; see \cite{DLM,LPG2002,LPG2004}). In order to address this issue, one needs to design special technique to accurately treat the nonconservative product term. This remains unclear for the proposed schemes in \cref{sec2} and also goes beyond the scope of the current study. Nevertheless, in \cref{appe}, we provide some numerical tests involving discontinuous bottom topographies. The results obtained from these tests demonstrate that, even without employing special treatment techniques, our proposed method performs comparably well to a path-conservative finite-volume scheme.
\end{rmk}

\subsubsection*{Example 1---Accuracy Test}
The goal of the first example is to numerically check the experimental rate of convergence. To this end, we consider the following initial data and bottom topography taken from \cite{BB}:
\begin{equation*}
  h(x,0)=0.3\Big[1+e^{-\frac{(x-0.5)^2}{0.05^2}}\Big]-0.2\cos(6\pi x),\quad q(x,0)=0,\quad Z(x)=0.2(1+\cos(6\pi x)),
\end{equation*}
which is prescribed in a computational domain $[0,1]$ and subject to the periodic boundary conditions.

We first compute the numerical solution until the final time $t =0.03$ using the WB AF scheme on a sequence of uniform meshes with $\dx=1/64$, $1/128$, $1/256$, $1/512$, $1/1024$, $1/2048$, and
$1/4096$. Then, we measure the discrete $L^1$-errors for the point values $\{h_j\}$ and $\{u_j\}$ at the cell interfaces using the following Runge formulae, which are based on the solutions computed on the three consecutive uniform grids with the size $\dx$, $2\dx$, and $4\dx$ and denoted by $(\cdot)^{\dx}$, $(\cdot)^{2\dx}$, and $(\cdot)^{4\dx}$, respectively:
\begin{equation*}
  \rm{Error}(\dx)\approx\frac{\varsigma_{12}^2}{|\varsigma_{12}-\varsigma_{24}|},
\end{equation*}
where $\varsigma_{12}:=\|(\cdot)^{\dx}-(\cdot)^{2\dx}\|_{L^1}$ and $\varsigma_{24}:=\|(\cdot)^{2\dx}-(\cdot)^{4\dx}\|_{L^1}$. The obtained $L^1$-errors and corresponding convergence rates are reported in Table \ref{tab1}. At the same time, we also compute the discrete $L^1$-errors as well as corresponding convergence rate for the average values $\{\xbar h_\jph\}$ and $\{\xbar q_\jph\}$ and report the results in Table \ref{tab1b}. From these obtained results, we can clearly see that the expected experimental third order of accuracy is achieved.

\begin{table}[ht!]
\caption{\sf Example 1: $L^1$-errors and experimental convergence rates of the point values $\{h_j\}$ and $\{u_j\}$.\label{tab1}}
\begin{center}
\begin{tabular}{|c|c|c|c|c|}
\hline
\multicolumn{1}{|c|}{$\dx$}&\multicolumn{1}{c|}{$L^1$-error in $h_j$}&{Rate}&\multicolumn{1}{c|}{$L^1$-error in $u_j$}&{Rate}\\
\hline
$1/256$&$2.87\times 10^{-4}$&-&$1.85\times 10^{-3}$&-\\
\hline
$1/512$&$7.32\times 10^{-5}$&1.97&$4.71\times 10^{-4}$&1.97\\
\hline
$1/1024$&$1.04\times 10^{-5}$&2.81&$6.56\times 10^{-5}$&2.85\\
\hline
$1/2048$&$1.37\times 10^{-6}$&2.92&$8.88\times 10^{-6}$&2.89\\
\hline
$1/4096$&$1.75\times 10^{-7}$&2.97&$1.23\times 10^{-6}$&2.86\\
\hline
\end{tabular}
\end{center}
\end{table}

\begin{table}[ht!]
\caption{\sf Example 1: $L^1$-errors and experimental convergence rates of the average values $\{\xbar h_\jph\}$ and $\{\xbar q_\jph\}$.\label{tab1b}}
\begin{center}
\begin{tabular}{|c|c|c|c|c|}
\hline
\multicolumn{1}{|c|}{$\dx$}&\multicolumn{1}{c|}{$L^1$-error in $\xbar h_\jph$}&{Rate}&\multicolumn{1}{c|}{$L^1$-error in $\xbar q_\jph$}&Rate\\
\hline
$1/256$&$3.70\times 10^{-4}$&-&$7.41\times 10^{-4}$&-\\
\hline
$1/512$&$6.54\times 10^{-5}$&2.50&$1.37\times 10^{-4}$&2.44\\
\hline
$1/1024$&$9.26\times 10^{-6}$&2.82&$1.95\times 10^{-5}$&2.81\\
\hline
$1/2048$&$1.23\times 10^{-6}$&2.91&$2.58\times 10^{-6}$&2.92\\
\hline
$1/4096$&$1.58\times 10^{-7}$&2.97&$3.29\times 10^{-7}$&2.97\\
\hline
\end{tabular}
\end{center}
\end{table}

\subsubsection*{Example 2---Still-Water Equilibrium and Its Small Perturbation}
In the second example, we validate that the WB AF scheme is capable of exactly preserving the still-water steady state. To this end, we take the bottom topography
\begin{equation*}
  Z(x)=\left\{\begin{aligned}
  &2(\cos(10\pi(x+0.3))+1),\quad && x\in[-0.4,-0.2],\\
  &0.5(\cos(10\pi(x-0.3))+1),\quad && x\in[0.2,0.4],\\
  &0,\quad &&\mbox{otherwise},
  \end{aligned}\right.
\end{equation*}
and consider the following initial data
\begin{equation*}
  h(x,0):=h_{\rm eq}(x)=4.000001-Z(x), \quad q(x,0):=q_{\rm eq}(x)=0,
\end{equation*}
which satisfies \eref{1.5} in the computational domain $[-1,1]$. We run the simulations until a long final time $t=30$ using both the WB AF and WB PCCU schemes with 100 uniform cells. Then we calculate the $L^1$- and $L^\infty$-errors in the average values of $h$ and $q$ and report the obtained errors in Table \ref{tab2}. As one can clearly see, both of the two studied schemes can preserve the still-water steady state at the level of round-off error. 

\begin{table}[ht!]
\caption{\sf Example 2: $L^1$- and $L^\infty$-errors in $h$ and $q$ computed by the WB AF and WB PCCU schemes.\label{tab2}}
\begin{center}
\begin{tabular}{|c|c|c|c|c|}
\hline
Schemes&$L^1$-error in $h$&$L^\infty$-error in $h$&$L^1$-error in $q$&$L^\infty$-error in $q$\\
\hline
WB AF&6.42e-13&3.42e-13&5.52e-13&4.25e-13\\
\hline
WB PCCU&3.32e-15&2.93e-15&7.94e-14&1.06e-13\\
\hline
\end{tabular}
\end{center}
\end{table}

\begin{figure}[ht!]
\centerline{\includegraphics[trim=0.8cm 0.25cm 1.1cm 0.2cm,clip,width=4.0cm]{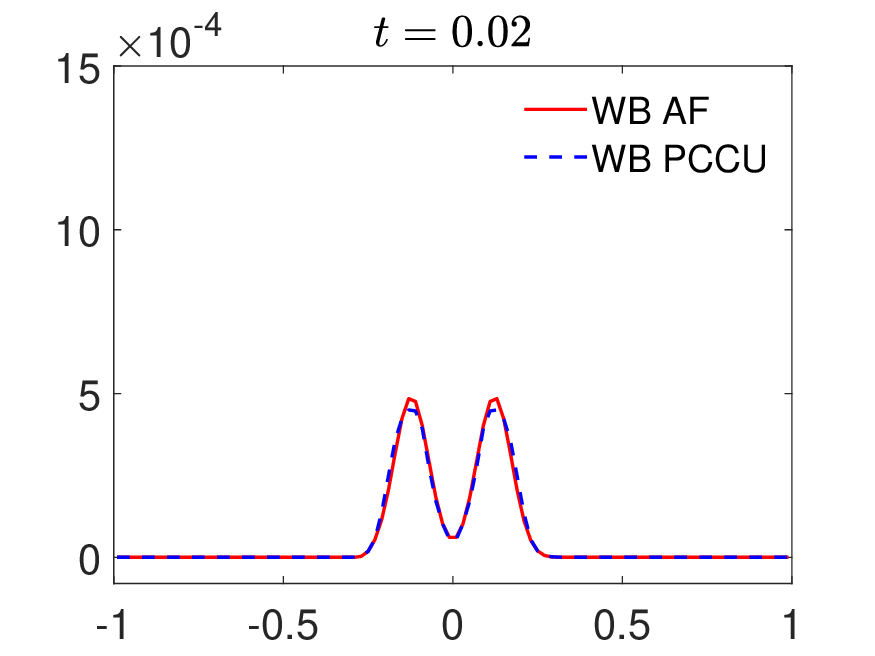}\hspace*{0.15cm}
\includegraphics[trim=0.8cm 0.25cm 1.1cm 0.2cm,clip,width=4.0cm]{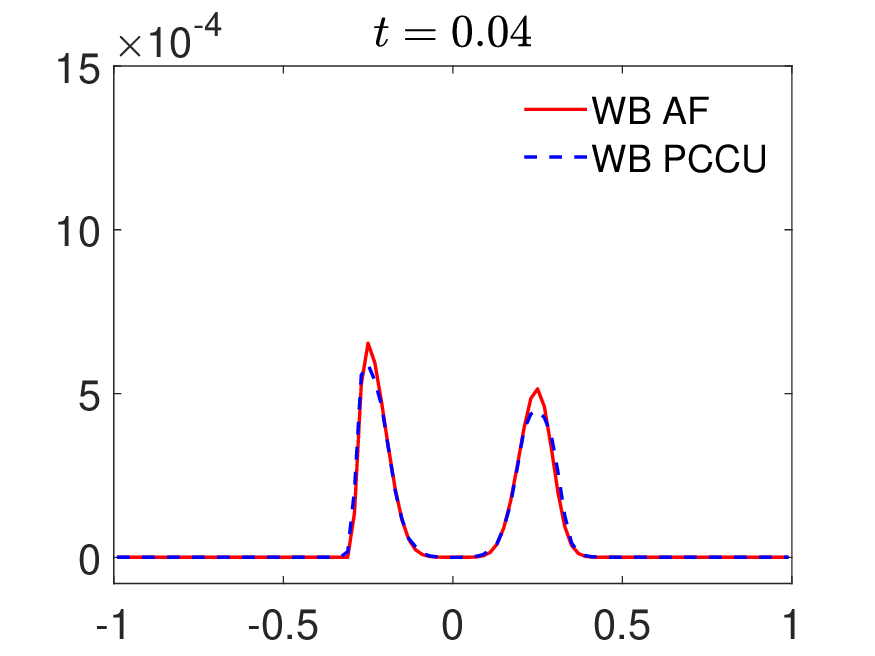}\hspace*{0.15cm}
\includegraphics[trim=0.8cm 0.25cm 1.1cm 0.2cm,clip,width=4.0cm]{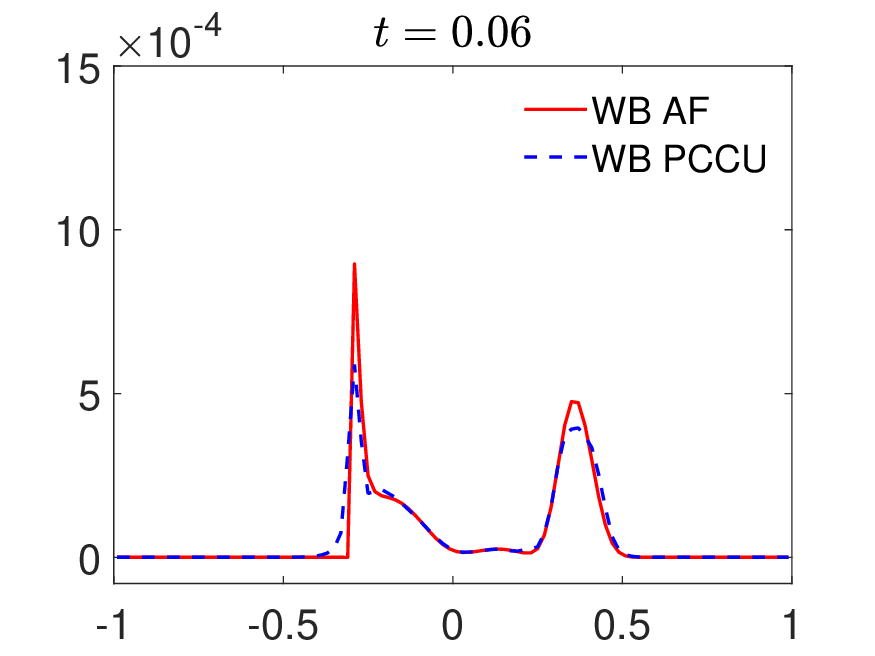}}
\vskip5pt
\centerline{\includegraphics[trim=0.8cm 0.25cm 1.1cm 0.2cm,clip,width=4.0cm]{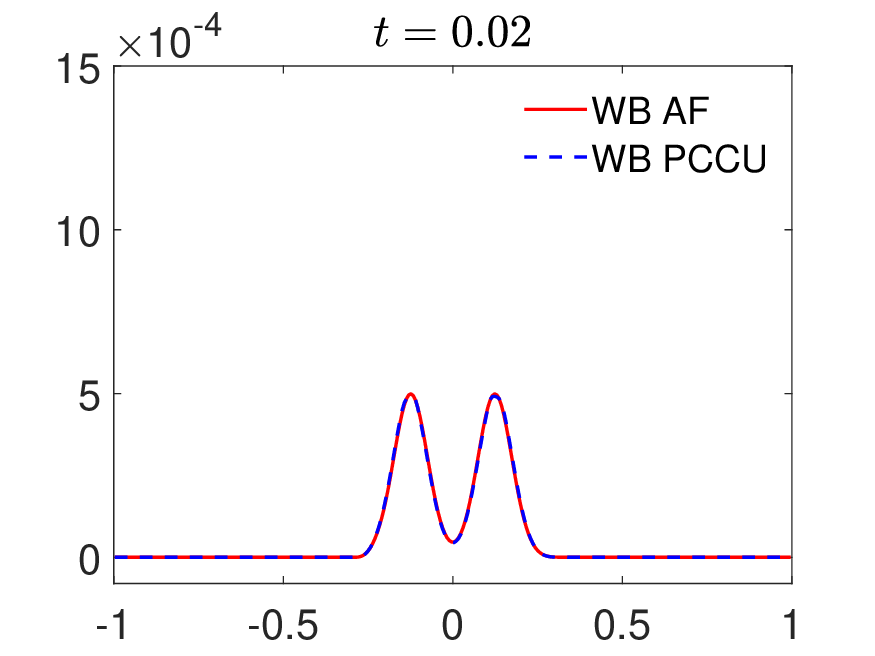}\hspace*{0.15cm}
\includegraphics[trim=0.8cm 0.25cm 1.1cm 0.2cm,clip,width=4.0cm]{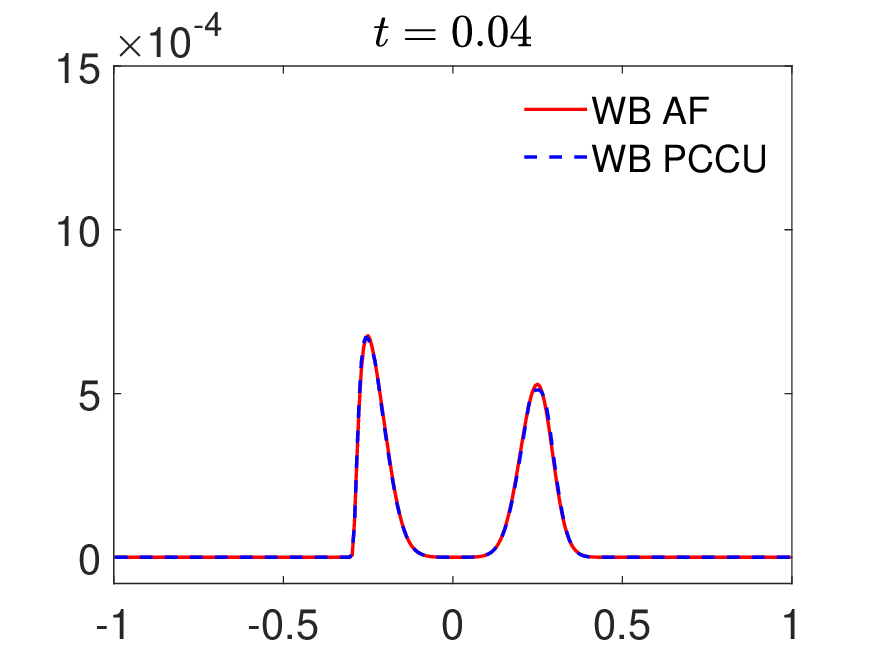}\hspace*{0.15cm}
\includegraphics[trim=0.8cm 0.25cm 1.1cm 0.2cm,clip,width=4.0cm]{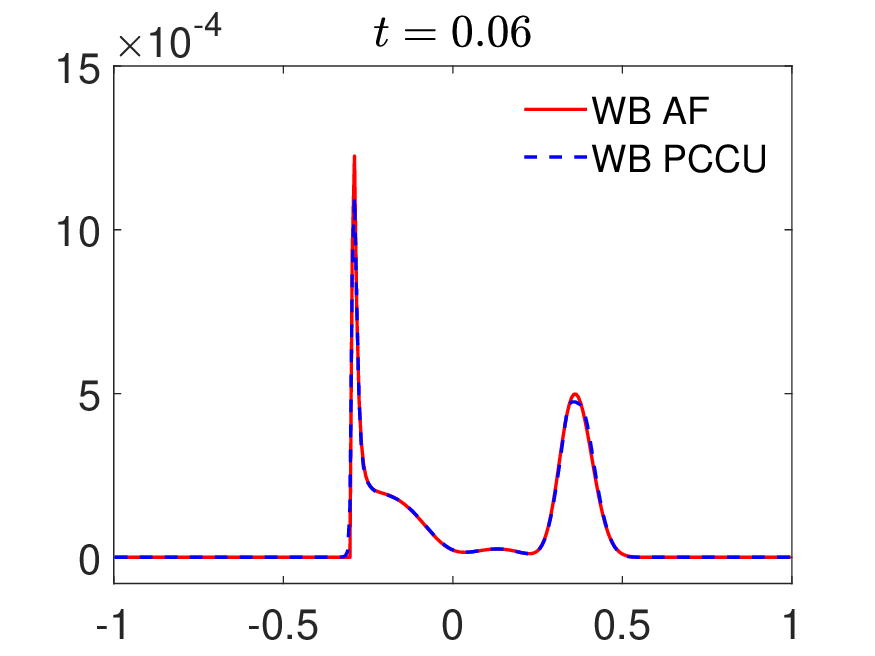}}
\caption{\sf Example 2 (small perturbation): Time snapshots of the difference $h-h_{\rm eq}$ computed by the WB AF and WB PCCU schemes using a coarse
mesh with $100$ cells (top row) and a finer mesh with $300$ cells (bottom row).\label{Ex2_fig1}}
\end{figure}

We now add a small Gaussian-shape perturbation to the stationary water depth and test the ability of the proposed schemes to accurately capture this small perturbation.  That is, we add a small perturbation, which takes the form of $10^{-3}e^{-200x^2}$, to the initial cell average of water depth and compute the numerical solutions at three different times $t=0.02$, $0.04$, and $0.06$ using either $100$ or $300$ uniform cells by both the WB AF and WB PCCU schemes. The differences between the computed and background cell averages of water depth are plotted in Figure \ref{Ex2_fig1}. As one can observe, the initial perturbation splits into two humps, which are then propagating into opposite directions, respectively. When the coarse mesh is used, the WB AF scheme generates much more accurate results. In particular, at the final time $t=0.06$, the result at $x\approx 0.3$ computed by the WB AF scheme is much shaper than its counterpart. When the mesh is refined, both of the two studied schemes converge to (almost) the same solution. 

Finally, we measure the difference in CPU times between the WB AF and WB PCCU schemes across the two different test cases. In the first test case of simulating the lake-at-rest equilibrium, the CPU time comparison indicates that the WB AF scheme, with MOOD criteria checks enabled, exhibits approximately 7.6\% higher computational cost than the WB PCCU scheme. However, excluding all MOOD criteria checks, deemed unnecessary for this test, yields a reduction of roughly 15.8\% in computational demand for the WB AF scheme compared to the WB PCCU scheme. Moving to the second test case, involving a small perturbation over a coarse mesh, the WB AF scheme without MOOD is approximately 7.8\% lower computational demand than the WB PCCU scheme. Conversely, with MOOD criteria checks enabled, the WB AF scheme is 6.8\% more computationally expensive than the WB PCCU scheme. Similar resluts are observed in the CPU time comparisons over the fine mesh. These findings prompt us to explore alternative more efficient strategies to replace the MOOD paradigm, which falls within the scope of the future work.

\subsubsection*{Example 3---Moving-Water Equilibria and Their Small Perturbations}
In the third example, we verify the WB property of the WB AF scheme in the sense that it is capable of exactly preserving moving-water equilibria \eref{1.6}. To this end, we consider a continuous bottom topography,
\begin{equation}\label{zz}
Z(x)=\left\{\begin{aligned}
&0.2-0.05(x-10)^2&&\mbox{if}~8\le x\le12,\\
&0&&\mbox{otherwise},
\end{aligned}\right.
\end{equation}
prescribed in a computational domain is $[0,25]$, and the following three sets of initial data that correspond to three different flow regimes (supercritical, subcritical, and transcritical without shock):
\begin{equation}\label{3.2}
\begin{aligned}
&\mbox{Case (a):}&&q_{\rm eq}(x,0)\equiv24,&&E_{\rm eq}(x,0)\equiv\frac{24^2}{2\times2^2}+9.812\times2;\\
&\mbox{Case (b):}&&q_{\rm eq}(x,0)\equiv4.42,&&E_{\rm eq}(x,0)\equiv\frac{4.42^2}{2\times2^2}+9.812\times2;\\
&\mbox{Case (c):}&&q_{\rm eq}(x,0)\equiv1.53,&&E_{\rm eq}(x,0)\equiv\frac{3}{2}(9.812\times1.53)^{\frac{2}{3}}+9.812\times0.2.
\end{aligned}
\end{equation}
The relevant steady-state water depths $h_{\rm eq}(x,0)$ are computed in the same manner as introduced in \cite[Remark 4.1]{CK16} and given by
\begin{equation}\label{3.3}
\begin{aligned}
&\mbox{Case (a):}~~ h_{\rm eq}(x,0)=-\frac{a_0}{3}\Big[2\cos\big(\frac{\theta+4\pi}{3}\big)+1\Big],\\
&\mbox{Case (b):}~~ h_{\rm eq}(x,0)=-\frac{a_0}{3}\Big[2\cos\frac{\theta}{3}+1\Big],\\
&\mbox{Case (c):}~~ h_{\rm eq}(x,0)=\left\{\begin{aligned}
&-\frac{a_0}{3}\Big[2\cos\big(\frac{\theta+4\pi}{3}\big)+1\Big]&&\quad \mbox{if}~x>10,\\
&-\frac{a_0}{3}\Big[2\cos\frac{\theta}{3}+1\Big]&&\quad \mbox{if}~x<10,\\
&-\frac{2a_0}{3}&&\quad \mbox{else},
\end{aligned}\right.
\end{aligned}
\end{equation}
where $a_0=\frac{gZ(x)-E_{\rm eq}(x,0)}{g}$, $a_2=\frac{(q_{\rm eq}(x,0))^2}{2g}$, and $\theta=\arccos(1+\frac{27a_2}{2a_0^3})$.

We first use the WB AF scheme with 100 uniform cells to compute the numerical solutions until a final time $t=20$ for all the three cases mentioned above. The discrete $L^1$- and $L^\infty$-errors in the average values of $h$ and $q$ are computed and reported in Tables \ref{tab3}. It is clear that the WB AF scheme can exactly preserve all three moving-water equilibria within machine accuracy.

\begin{table}[ht!]
\caption{\sf Example 3: $L^1$- and $L^\infty$-errors in average values of $h$ and $q$ over continuous bathymetry \eref{zz} .\label{tab3}}
\begin{center}
\begin{tabular}{|c|c|c|c|c|}
\hline
&$L^1$-error in $h$&$L^\infty$-error in $h$&$L^1$-error in $q$&$L^\infty$-error in $q$\\
\hline
Case (a)&3.23e-13&3.64e-14&3.58e-12&3.77e-13\\
\hline
Case (b)&2.90e-14&6.22e-15&6.09e-13&6.22e-14\\
\hline
Case (c)&2.65e-14&4.11e-15&7.03e-14&1.22e-14\\
\hline
\end{tabular}
\end{center}
\end{table}

\begin{figure}[ht!]
\centerline{\includegraphics[trim=0.8cm 0.25cm 0.9cm 0.2cm,clip,width=4.0cm]{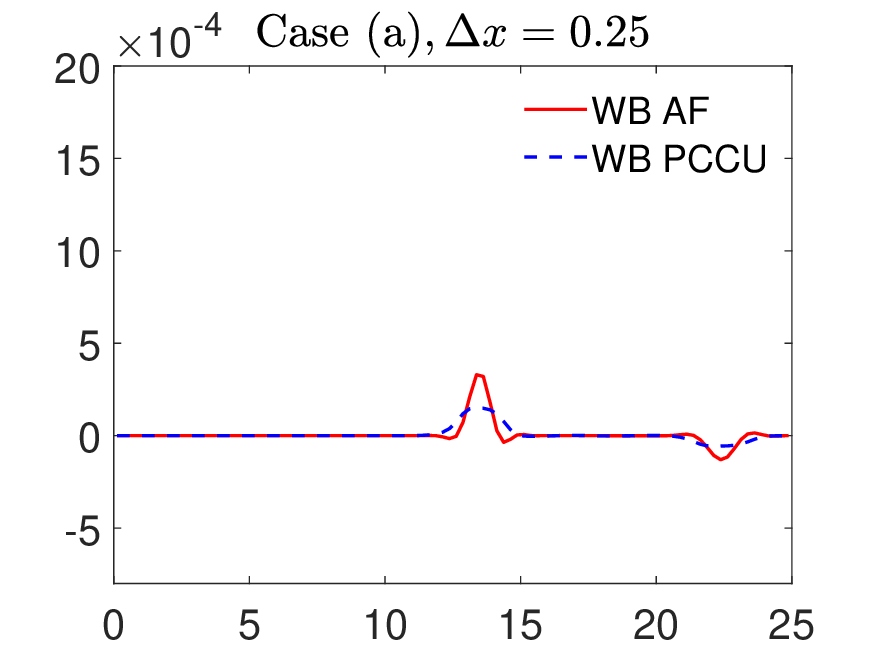}\hspace*{0.15cm}
\includegraphics[trim=0.8cm 0.25cm 0.9cm 0.2cm,clip,width=4.0cm]{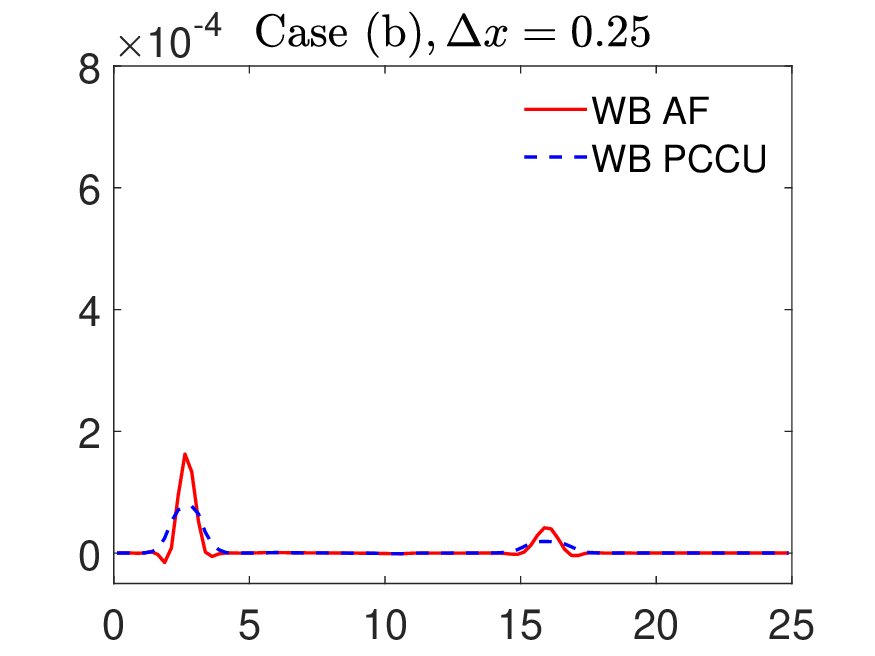}\hspace*{0.15cm}
\includegraphics[trim=0.8cm 0.25cm 0.9cm 0.2cm,clip,width=4.0cm]{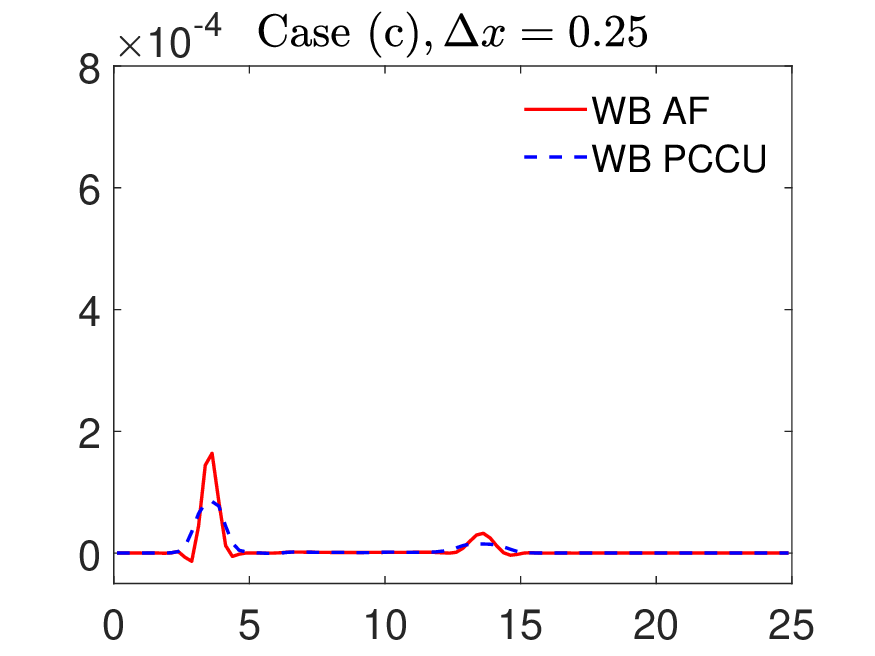}}
\vskip5pt
\centerline{\includegraphics[trim=0.8cm 0.25cm 0.9cm 0.2cm,clip,width=4.0cm]{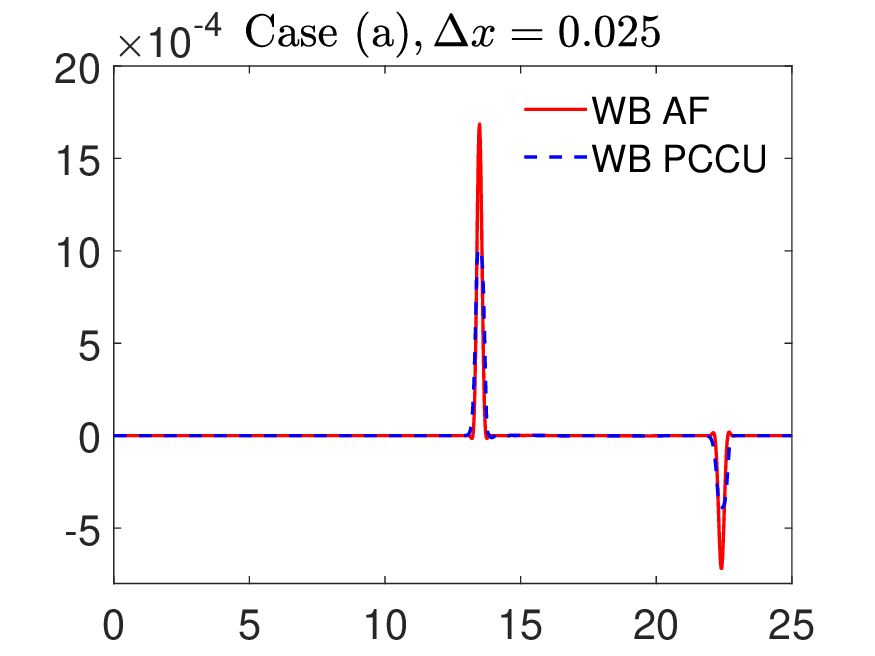}\hspace*{0.15cm}
\includegraphics[trim=0.8cm 0.25cm 0.9cm 0.2cm,clip,width=4.0cm]{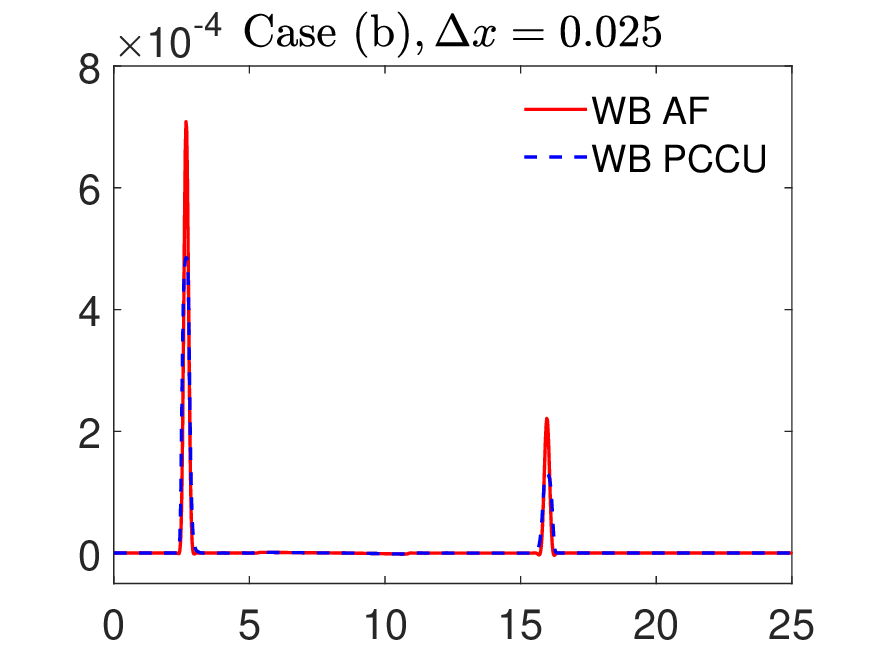}\hspace*{0.15cm}
\includegraphics[trim=0.8cm 0.25cm 0.9cm 0.2cm,clip,width=4.0cm]{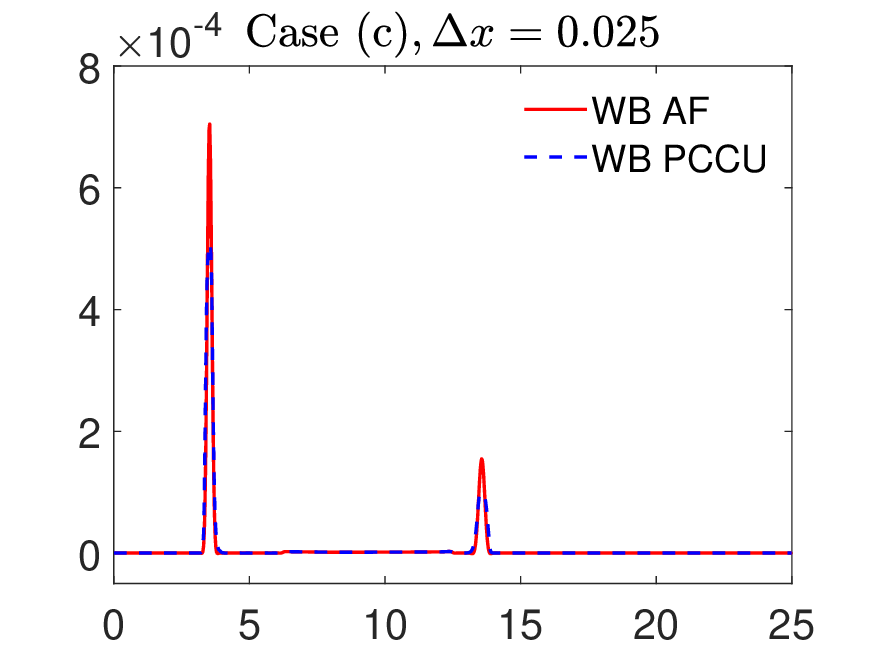}}
\caption{\sf Example 3: The difference between the obtained and the steady-state cell averages of water depth computed by both the WB AF and WB PPCU schemes using uniform cells with mesh size $\dx=0.25$ (top row) and $\dx=0.025$ (bottom row) for Cases (a, left column), (b, middle column) and (c, right column).\label{fig2}}
\end{figure}

We then add a small perturbation, which takes the form of $10^{-3}e^{-80(x-6)^2}$, to the steady-state cell average of water depth $h_{\rm eq}(x,0)$ in \eref{3.3}. We compute the numerical solutions using the WB AF and WB PPCU schemes until $t=1$ in the supercritical case and until $t=1.5$ in the subcritical and transcritical cases. The difference between the obtained and the background moving steady-state cell averages of water depth computed with $100$ and $1000$ uniform cells are plotted in Figure \ref{fig2}. As one can observe, the propagating perturbations are well captured by both studied schemes over both coarse and fine meshes. Comparing the results obtained by the alternative WB PCCU scheme, we can see that the WB AF scheme produces much higher resolution solutions.

\subsection*{Example 4---Convergence to Moving-Water Equilibria}
In the fourth example, we study the convergence in time of the numerical solution computed by the proposed WB AF scheme towards steady flow over the continuous bathemetry given by \eref{zz}. Depending on the initial and boundary conditions, the obtained convergent flow may be supercritical, subcritical, or transcritical without or with a steady shock. The four sets of initial and boundary conditions in terms of water depth and velocity are summarized as follows:

\begin{equation*}
   \begin{aligned}
&\mbox{Case (a):}&&\left\{\begin{array}{l}h(x,0)=2-Z(x),\quad u(x,0)=0,\\h(0,t)=2,\qquad\qquad~~ u(0,t)=12;\end{array}\right.\\
&\mbox{Case (b):}&&\left\{\begin{array}{l}h(x,0)=2-Z(x),\quad u(x,0)=0,\\h(0,t)=2, \qquad\qquad~~ u(0,t)=2.21,\quad h(25,t)=2;\end{array}\right.\\
&\mbox{Case (c):}&&\left\{\begin{array}{l}h(x,0)=0.66-Z(x),\quad~ u(x,0)=0,\\h(0,t)=1.01439,\qquad\quad~ u(0,t)=1.53/h(0,t), \quad h(25,t)=0.66;\end{array}\right.\\
&\mbox{Case (d):}&&\left\{\begin{array}{l}h(x,0)=0.33-Z(x),\quad u(x,0)=0,\\h(0,t)=0.41372,\qquad\quad u(0,t)=0.18/h(0,t), \quad h(25,t)=0.33. \end{array}\right.
\end{aligned}
\end{equation*}

\begin{rmk}
In Case (c), the downstream boundary condition $(h(25,t)=0.66)$ is imposed only when the flow is subcritical.
\end{rmk}

\begin{figure}[ht!]
\centerline{\includegraphics[trim=0cm 0.25cm 0.9cm 0.12cm,clip,width=4.2cm]{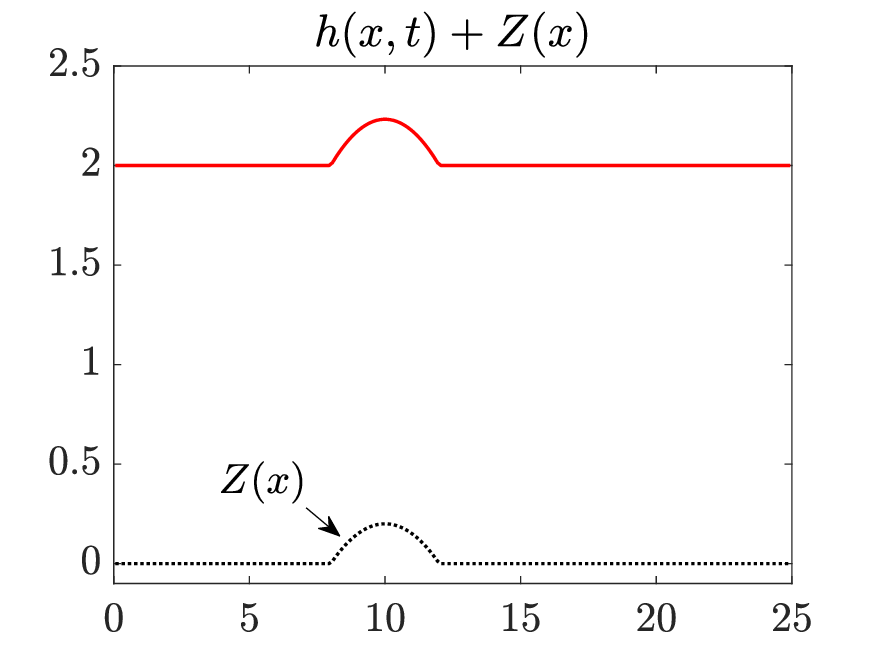}\hspace*{0.15cm}
\includegraphics[trim=0cm 0.25cm 0.9cm 0.12cm,clip,width=4.2cm]{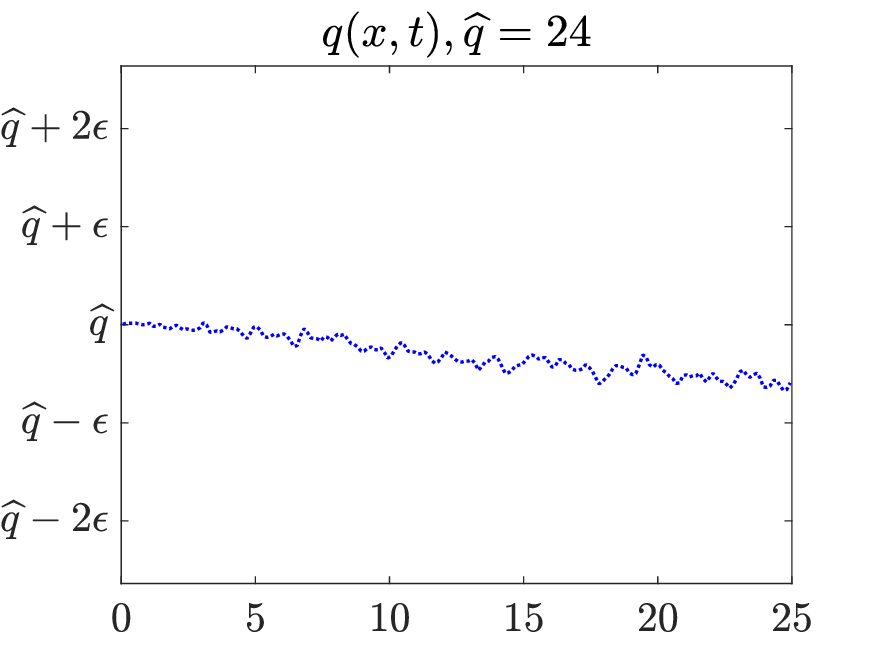}\hspace*{0.15cm}
\includegraphics[trim=0cm 0.25cm 0.9cm 0.12cm,clip,width=4.2cm]{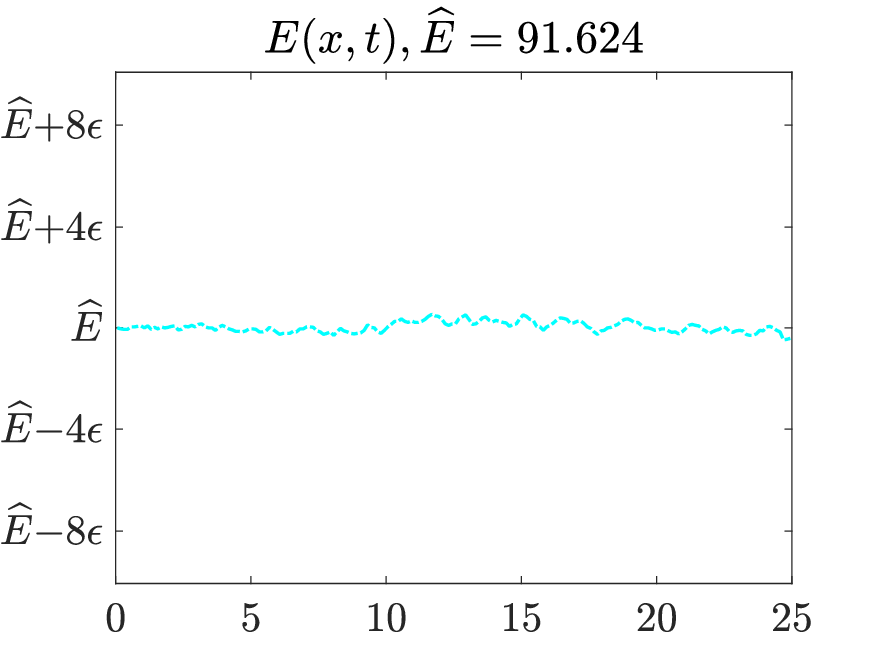}}
\vskip5pt
\centerline{\includegraphics[trim=0cm 0.25cm 0.9cm 0.12cm,clip,width=4.2cm]{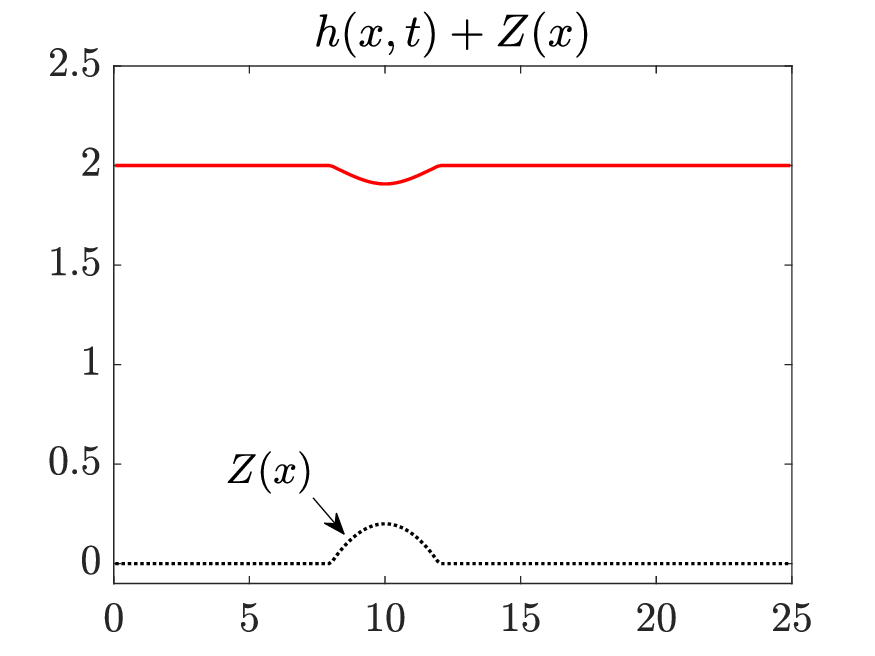}\hspace*{0.15cm}
\includegraphics[trim=0cm 0.25cm 0.9cm 0.12cm,clip,width=4.2cm]{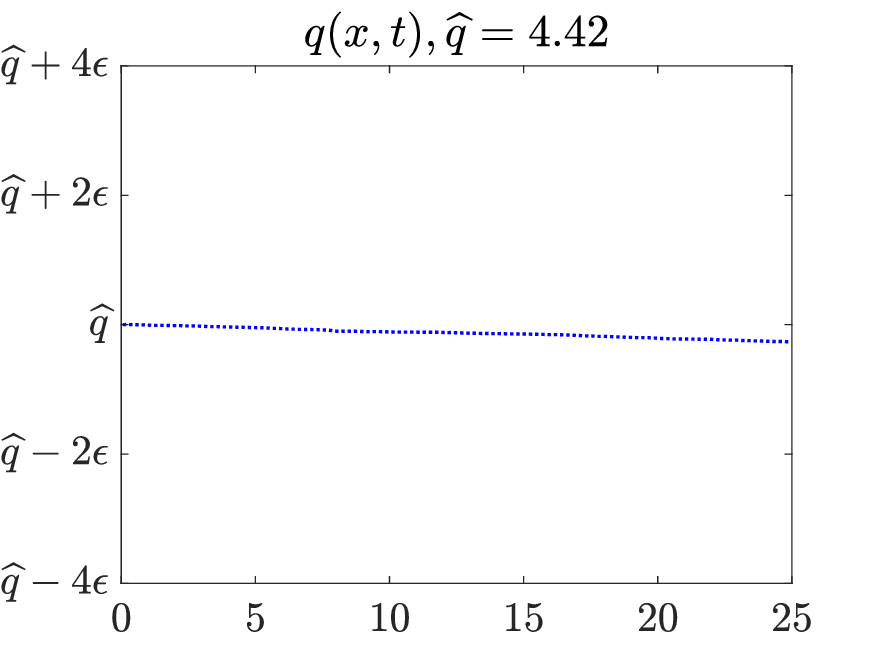}\hspace*{0.15cm}
\includegraphics[trim=0cm 0.25cm 0.9cm 0.12cm,clip,width=4.2cm]{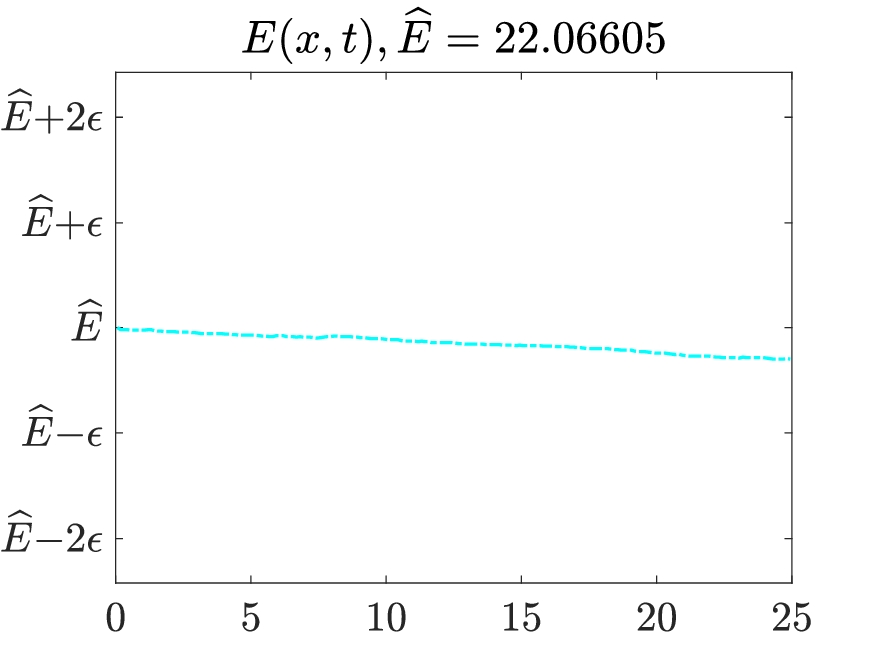}}
\vskip5pt
\centerline{\includegraphics[trim=0cm 0.25cm 0.9cm 0.12cm,clip,width=4.2cm]{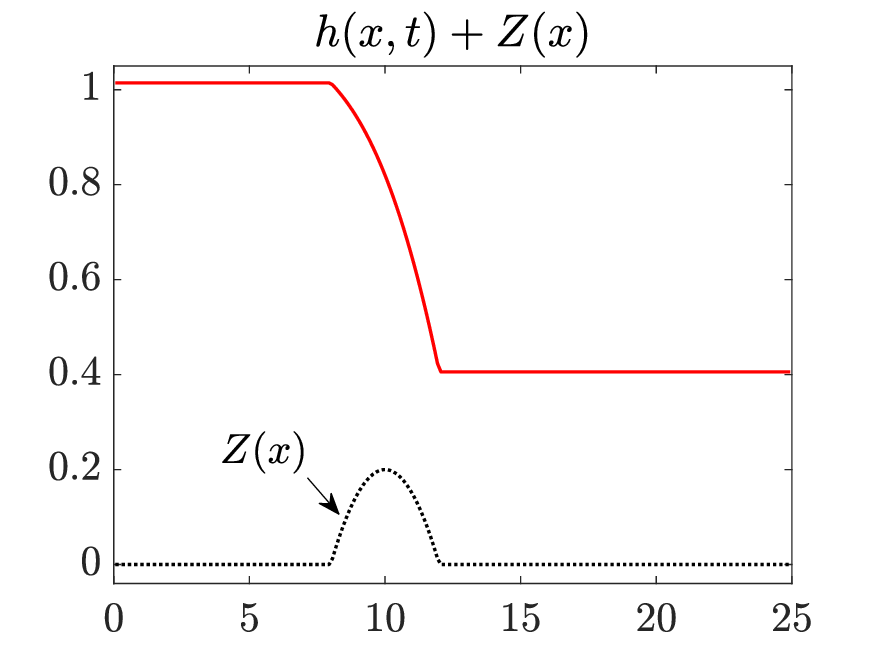}\hspace*{0.15cm}
\includegraphics[trim=0cm 0.25cm 0.9cm 0.12cm,clip,width=4.2cm]{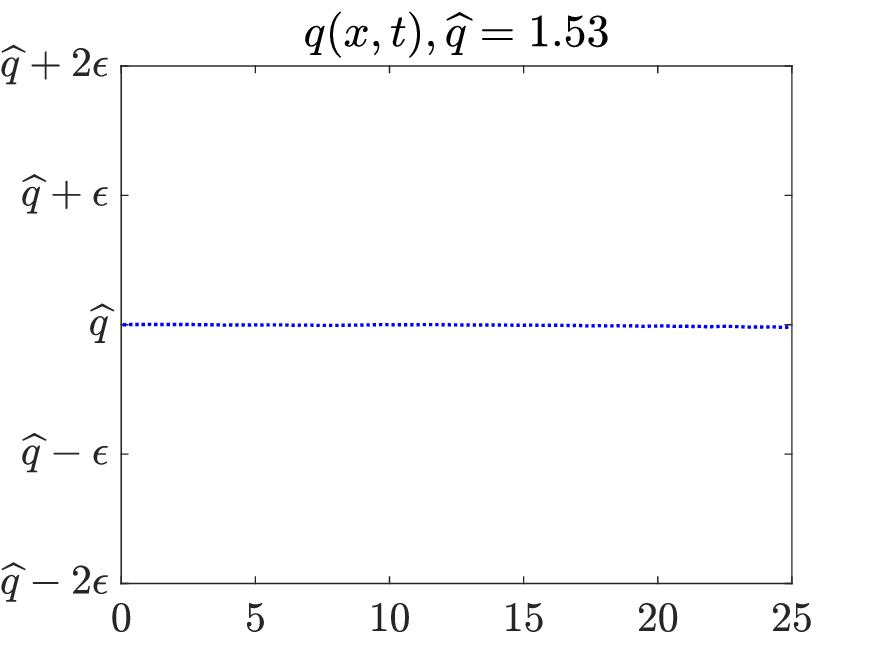}\hspace*{0.15cm}
\includegraphics[trim=0cm 0.25cm 0.9cm 0.12cm,clip,width=4.2cm]{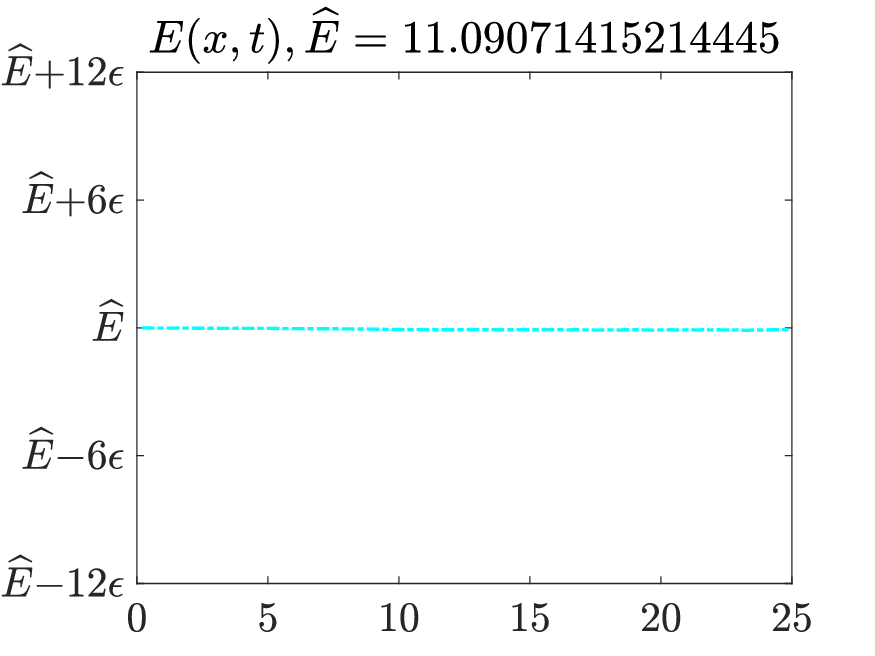}}
\vskip5pt
\centerline{\includegraphics[trim=0cm 0.25cm 0.9cm 0.12cm,clip,width=4.2cm]{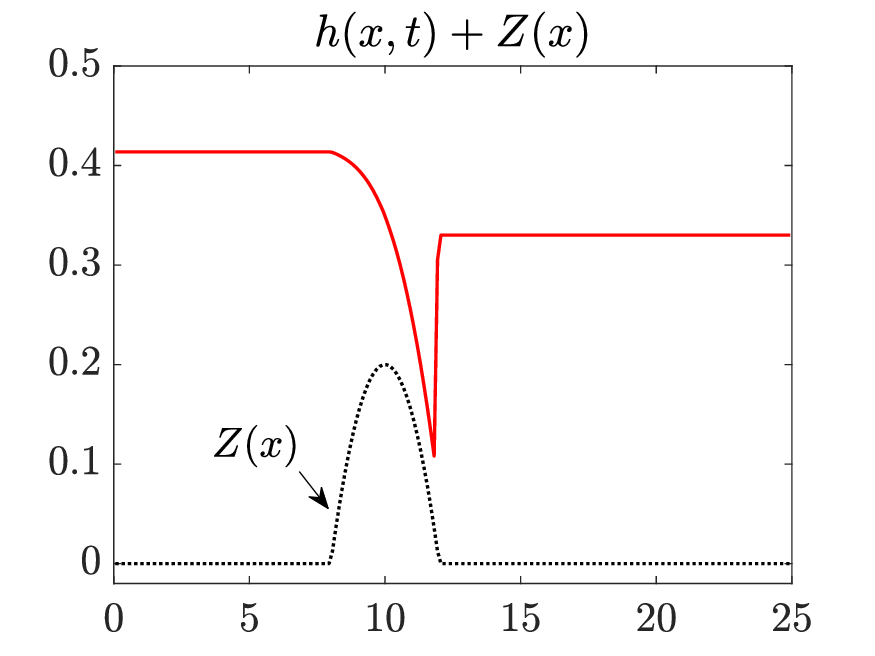}\hspace*{0.15cm}
\includegraphics[trim=0cm 0.25cm 0.9cm 0.12cm,clip,width=4.2cm]{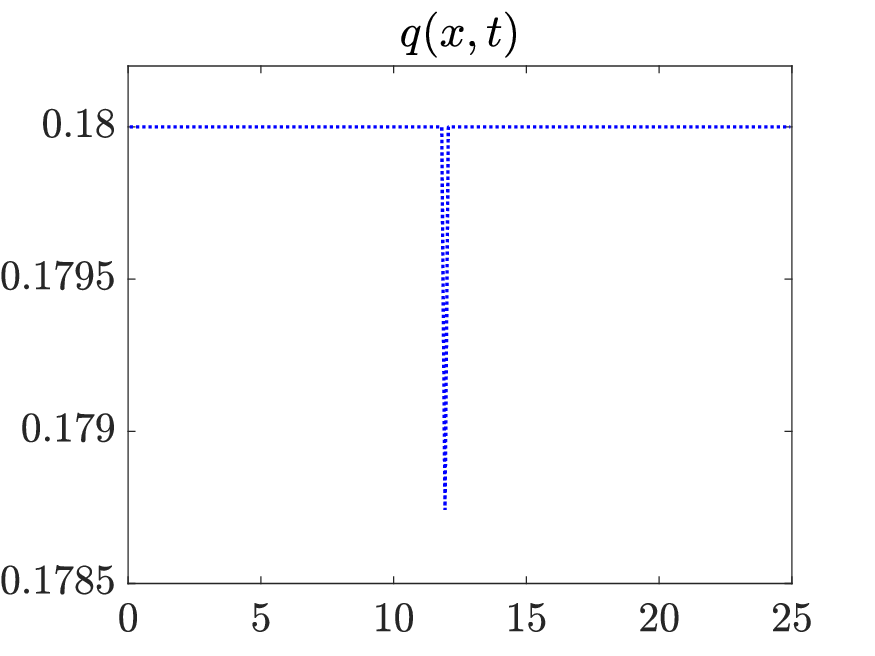}\hspace*{0.15cm}
\includegraphics[trim=0cm 0.25cm 0.9cm 0.12cm,clip,width=4.2cm]{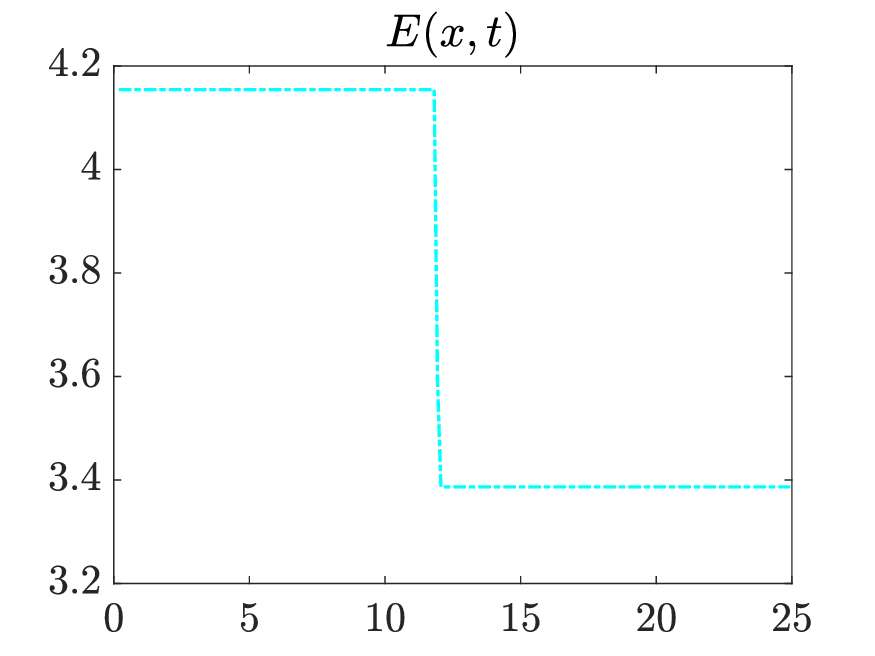}}
\caption{\sf Example 4: Convergent solutions of $h+Z$ together with $Z$ (left column), $q$ (middle column), and $E$ (right column) for Cases (a, first row), (b, second row), (c, third row), and (d, fourth row). The small deviation parameter $\epsilon=10^{-12}$. \label{fig4}}
\end{figure}

The numerical solutions ($h+Z$, $q$, and $E$) at a final time $t=500$ on the computational domain $[0,25]$ computed by the WB AF scheme with 200 uniform cells are depicted in Figure \ref{fig4}. As one can see, the numerical solutions of the first three cases have almost converged to the desired discrete steady states. However, for Case (d), like most of the existing schemes, the proposed scheme also fails to accurately resolve the equilibrium as the errors at the shock are $\mathcal{O}(1)$. It is easy to find that the numerical solutions obtained here are comparable with those reported in \cite{CKLX,LCJKY,CCHKW_18,CK16}.

\subsection*{Example 5---Dam-break Problem}
\begin{figure}[ht!]
\centerline{\includegraphics[trim=0.9cm 0.3cm 0.6cm 0.2cm,clip,width=5.2cm]{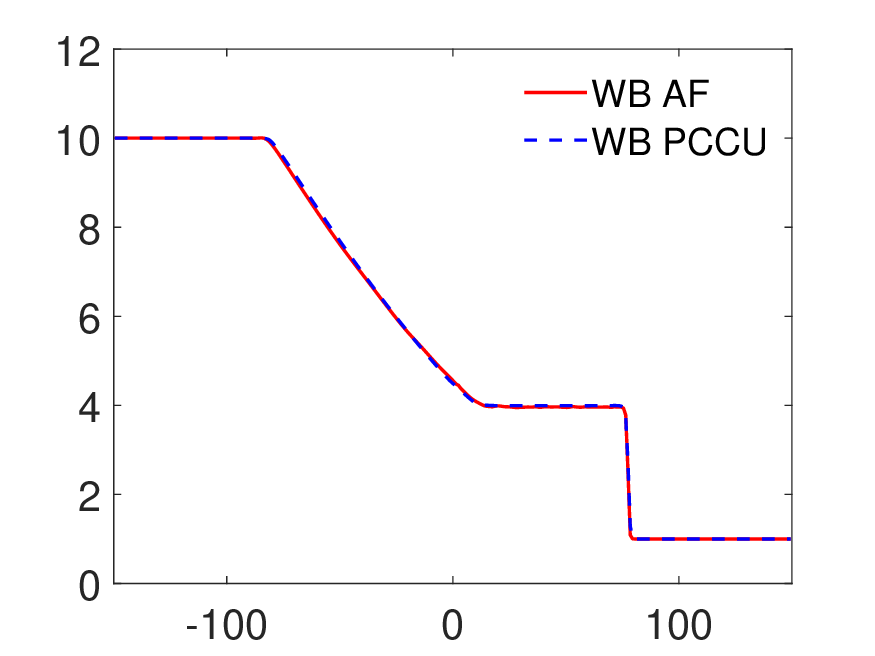}\hspace*{0.3cm}
\includegraphics[trim=0.9cm 0.3cm 0.6cm 0.2cm,clip,width=5.2cm]{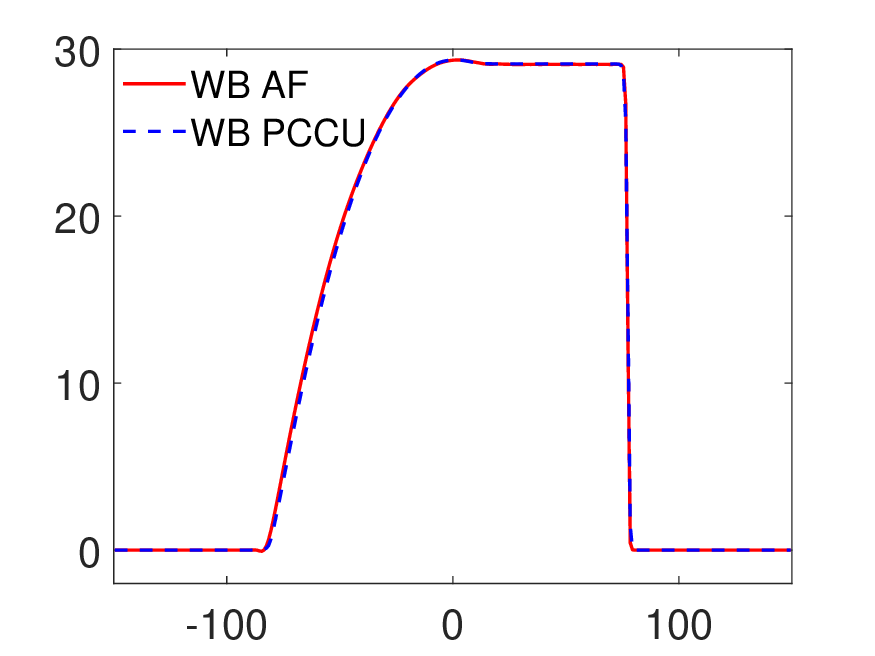}}
\caption{\sf Example 5: Water depth $h$ (left) and discharge $q$ (right) computed by the \bla{WB AF} and \bla{WB PCCU} schemes.\label{Ex5fig}}
\end{figure}

In the fifth example, we study a dam-break problem with a mixed Riemann solution, i.e., the solution consists of one rarefaction and one shock wave. The initial data are
\begin{equation*}
  (h(x,0), u(x,0), Z(x))=\left\{\begin{aligned}
  &(10, 0, 0),\quad &&x<0,\\
  &(1, 0, 0),\quad &&x\geq0.\\
  \end{aligned}\right.
\end{equation*}
The exact solution of this problem is given by a left-propagating rarefaction and a right-propagating shock wave. We use both the \bla{WB AF} and \bla{WB PCCU} schemes to compute the numerical solutions at the final time $t=8$ in the computational domain $[-150,150]$ covered with $300$ uniform cells. Results are shown in Figure \ref{Ex5fig}. As one can see, the obtained profiles of $h$ and $q$ provide expected results with a left-directed rarefaction and a right-directed shock wave. The numerical solutions computed by the two studied schemes are in a satisfactory agreement.

\subsection*{Example 6---Riemann Problem}
\begin{figure}[ht!]
\centerline{\includegraphics[trim=0.5cm 0.25cm 0.9cm 0.2cm,clip,width=4.2cm]{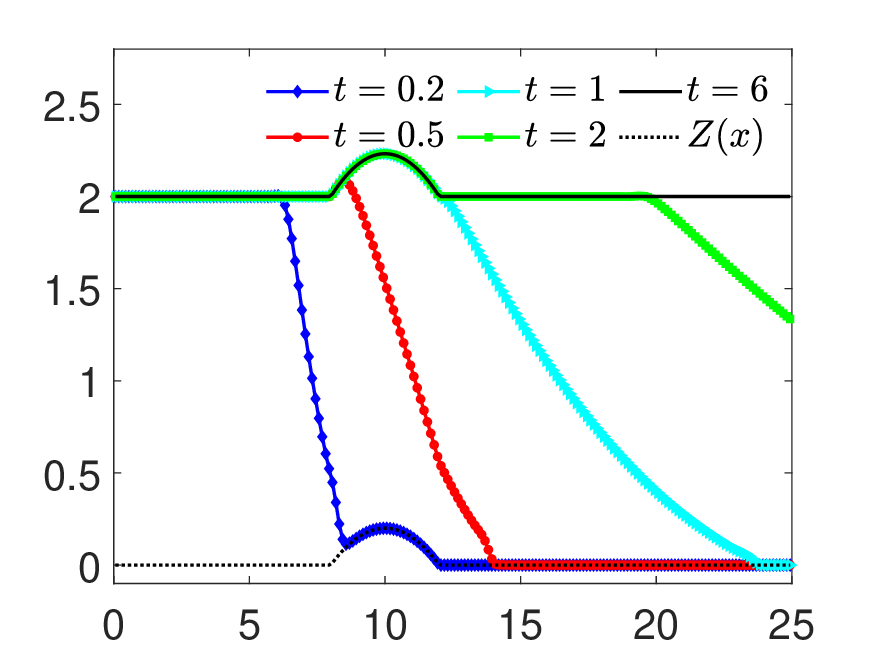}\hspace*{0.15cm}
\includegraphics[trim=0.5cm 0.25cm 0.9cm 0.2cm,clip,width=4.2cm]{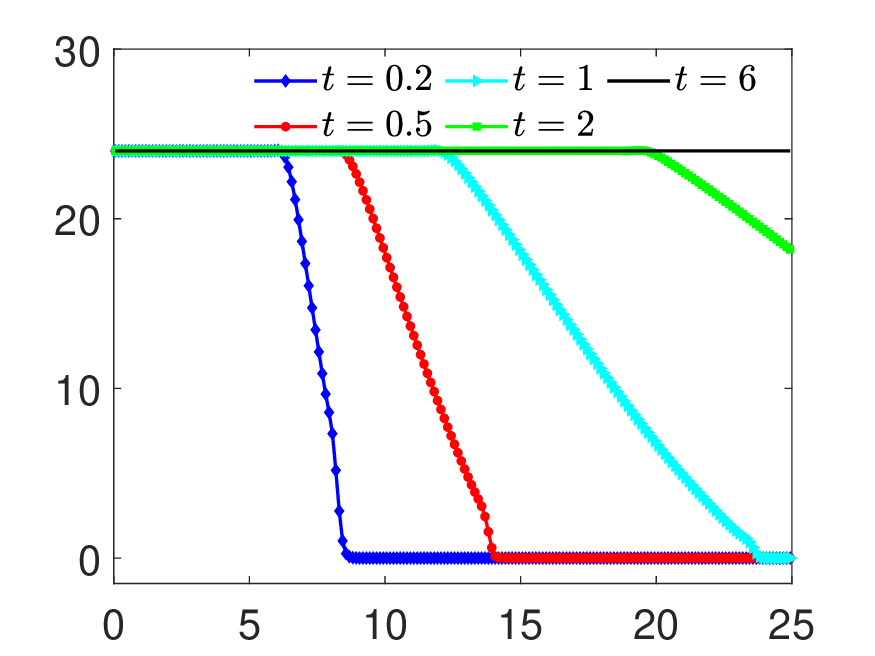}\hspace*{0.15cm}
\includegraphics[trim=0.5cm 0.25cm 0.9cm 0.2cm,clip,width=4.2cm]{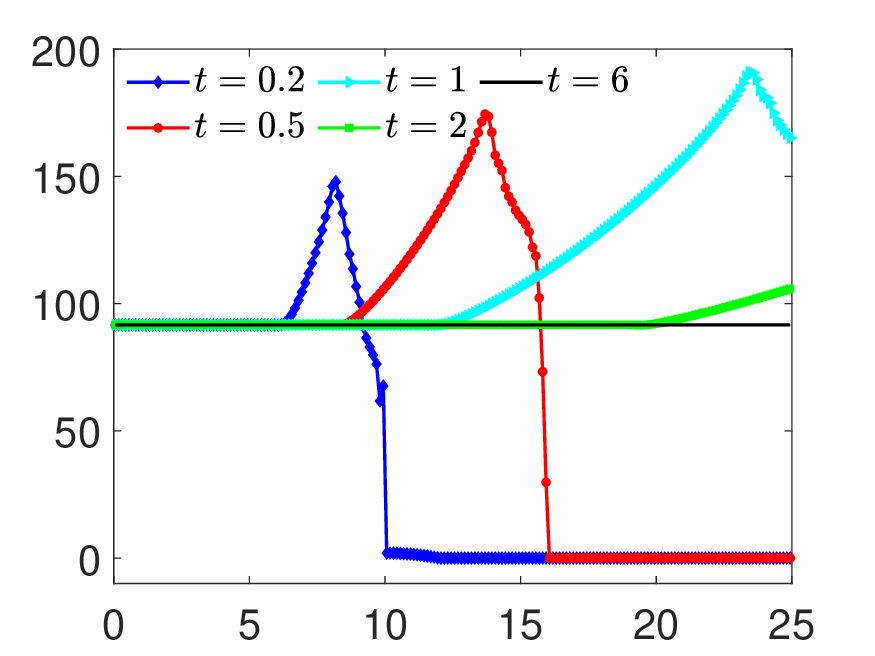}}
\caption{\sf Example 6: Computed solutions of $h+Z$ together with $Z$ (left), $q$ (middle), and $E$ (right).\label{fig7}}
\end{figure}

In the sixth example, we consider the following initial value problem in $[0,25]$ with the Riemann initial data to verify the positivity-preserving property of the \bla{WB AF} scheme:
\begin{equation*}
  h(x,0)=\left\{\begin{aligned}
  &2 \quad &&\mbox{if}~x<5,\\
  &0\quad &&\mbox{otherwise},
  \end{aligned}\right. \quad
  u(x,0)=\left\{\begin{aligned}
  &12 \quad &&\mbox{if}~x<5,\\
  &0\quad &&\mbox{otherwise},
  \end{aligned}\right. 
\end{equation*}
subject to the Dirichlet boundary condition:
\begin{equation*}
  h(0,t)=2,\quad u(0,t)=12.
\end{equation*}
The bottom topography is given by \eref{zz}. We compute the solutions using $200$ uniform cells at five different times: $t=0.2$, $0.5$, $1$, $2$, and $6$. The numerical solutions ($h+Z$, $q$, and $E$) of these time snapshots are plotted in Figure \ref{fig7}. As one can clearly observe that, \bla{with time's progression, the water gradually flows through the bottom hump. At $t=6$, it reaches the same steady state as in case (a) of Example 3, where both the discharge $q$ and the energy $E$ tend to be constants}. The initial data contains dry area $(h=0)$ and the obtained results are stable and thus show that the proposed \bla{WB AF} scheme possesses the positivity-preserving property.

\subsection*{Example 7---Convergence to Moving-Water Equilibria for Frictional case}
In the seventh example, which is a modification of Example 4, we now consider a nonzero Manning friction term and study the convergence in time of the numerical solution towards discrete steady states over a continuous hump given by \eref{zz}. Three cases corresponding to the supercritical [Case (a)], subcritical [Case (b)], and transcritical flow with a steady shock [Case (d)] are taken into account.  We take ``lake-at-rest'' initial conditions and the Dirichlet boundary conditions:
\begin{equation*}
   \begin{aligned}
&\mbox{Case (a):}&&\left\{\begin{array}{l}h(x,0)=2-Z(x),\quad u(x,0)=0,\\h(0,t)=2,\qquad\qquad~~ u(0,t)=12;\end{array}\right.\\
&\mbox{Case (b):}&&\left\{\begin{array}{l}h(x,0)=2-Z(x),\quad u(x,0)=0,\\h(0,t)=2.14618, \quad~~ u(0,t)=4.42/h(0,t),\quad h(25,t)=2;\end{array}\right.\\
&\mbox{Case (d):}&&\left\{\begin{array}{l}h(x,0)=0.33-Z(x),~~ u(x,0)=0,\\h(0,t)=0.44835,\quad\quad~~ u(0,t)=0.18/h(0,t), \quad h(25,t)=0.33. \end{array}\right.
\end{aligned}
\end{equation*}

\begin{figure}[ht!]
\centerline{\includegraphics[trim=0cm 0.25cm 0.9cm 0.15cm,clip,width=4.2cm]{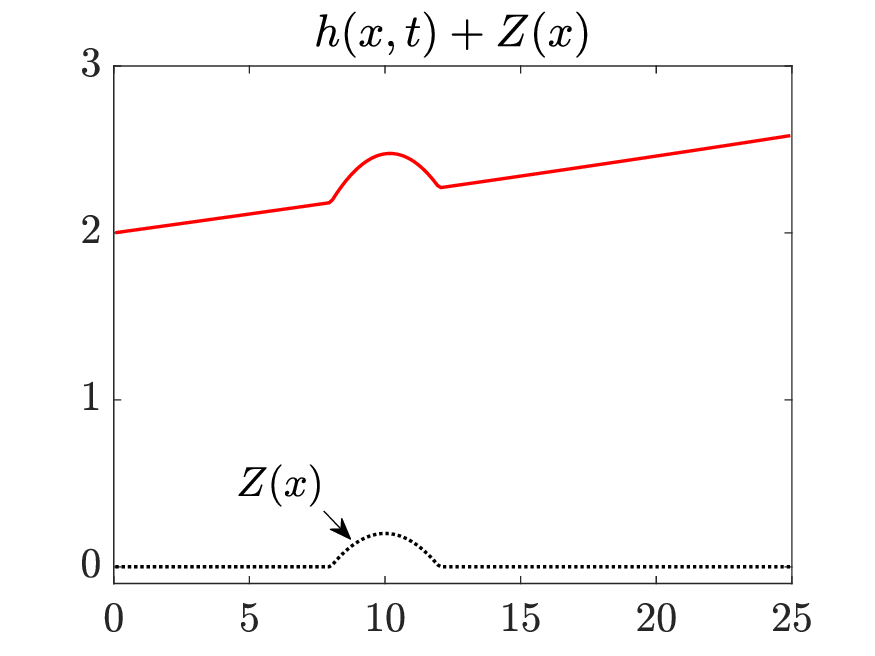}\hspace*{0.15cm}
\includegraphics[trim=0cm 0.25cm 0.9cm 0.15cm,clip,width=4.2cm]{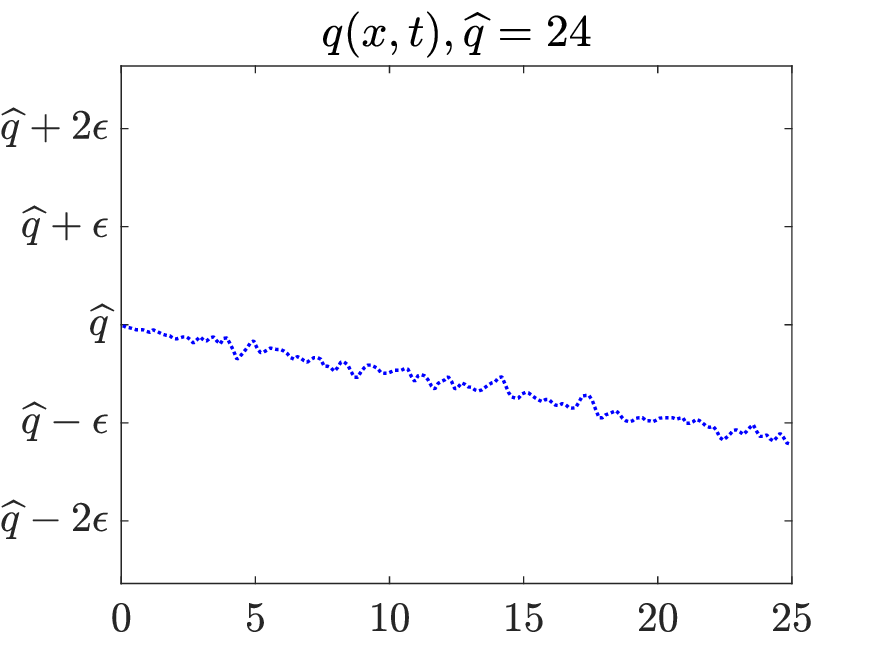}\hspace*{0.15cm}
\includegraphics[trim=0cm 0.25cm 0.9cm 0.15cm,clip,width=4.2cm]{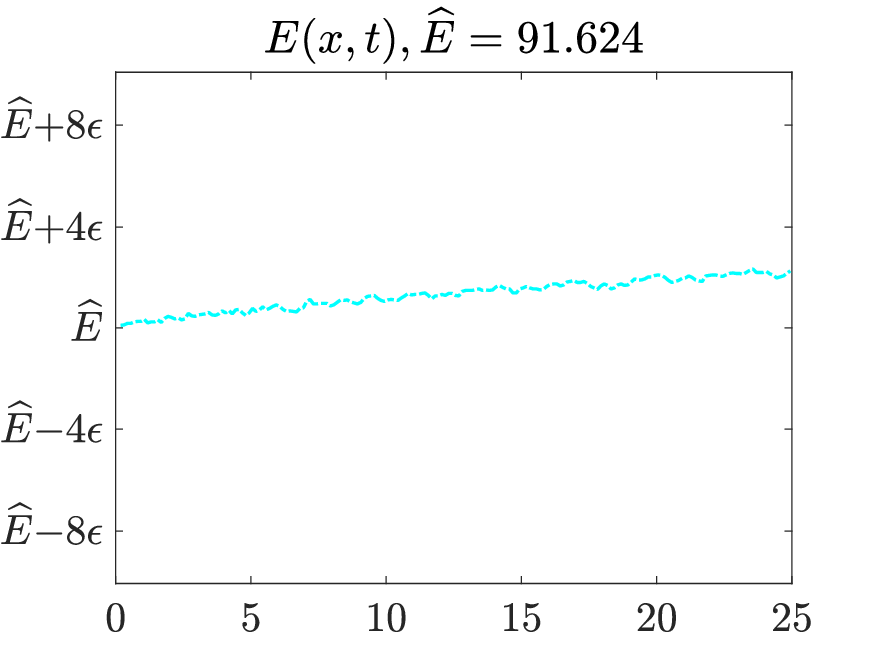}}
\vskip5pt
\centerline{\includegraphics[trim=0cm 0.25cm 0.9cm 0.15cm,clip,width=4.2cm]{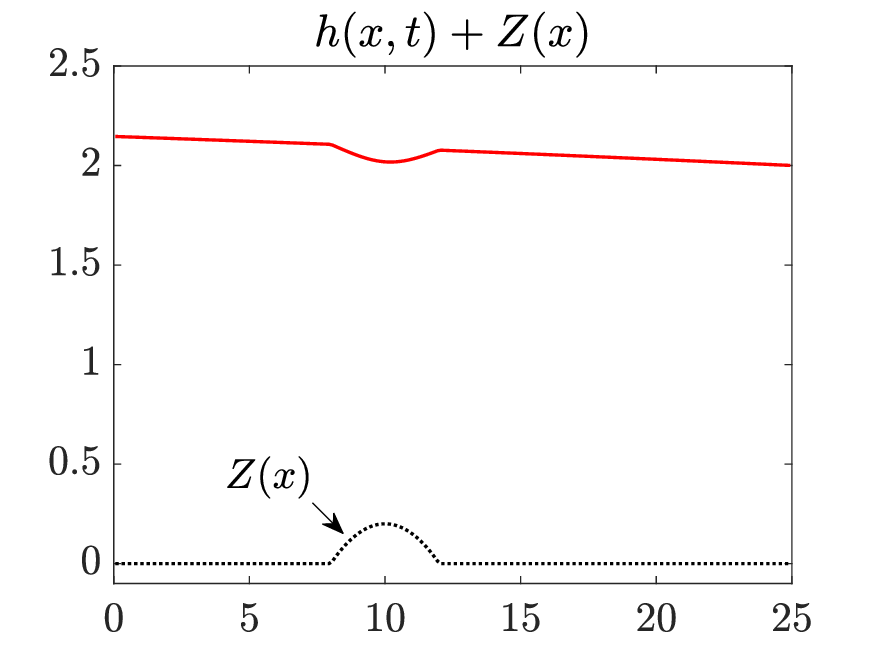}\hspace*{0.15cm}
\includegraphics[trim=0cm 0.25cm 0.9cm 0.15cm,clip,width=4.2cm]{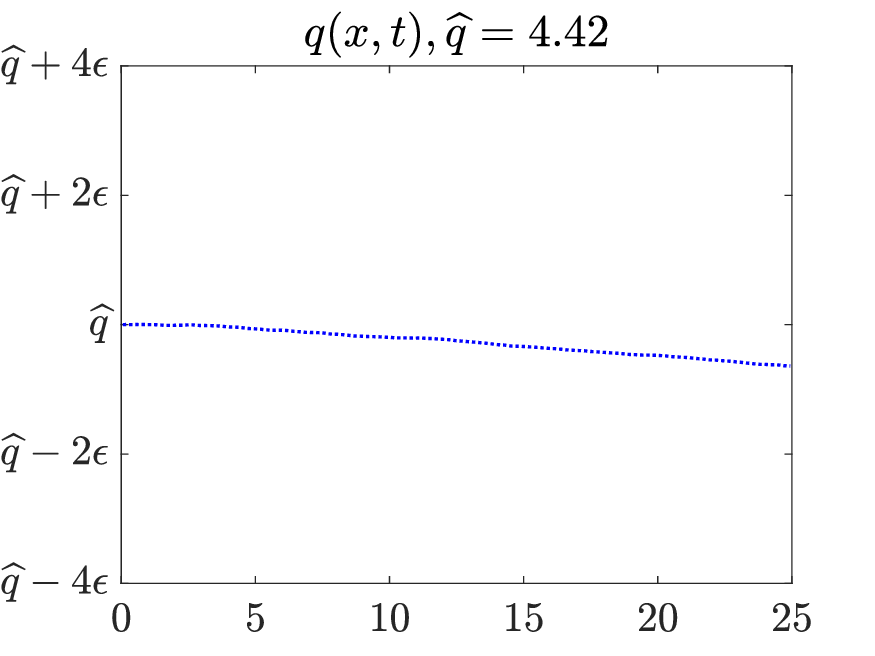}\hspace*{0.15cm}
\includegraphics[trim=0cm 0.25cm 0.9cm 0.15cm,clip,width=4.2cm]{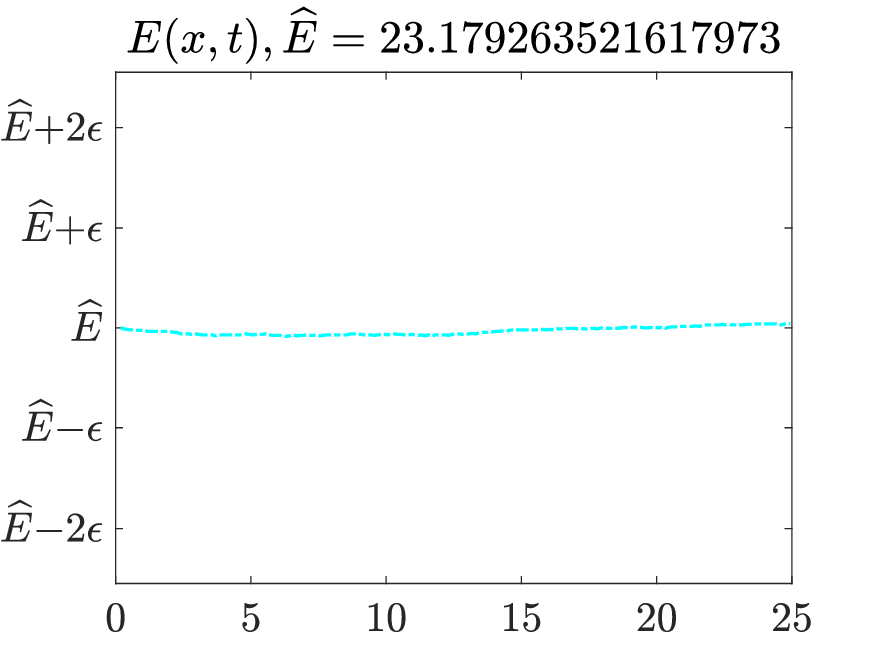}}
\vskip5pt
\centerline{\includegraphics[trim=0cm 0.25cm 0.9cm 0.15cm,clip,width=4.2cm]{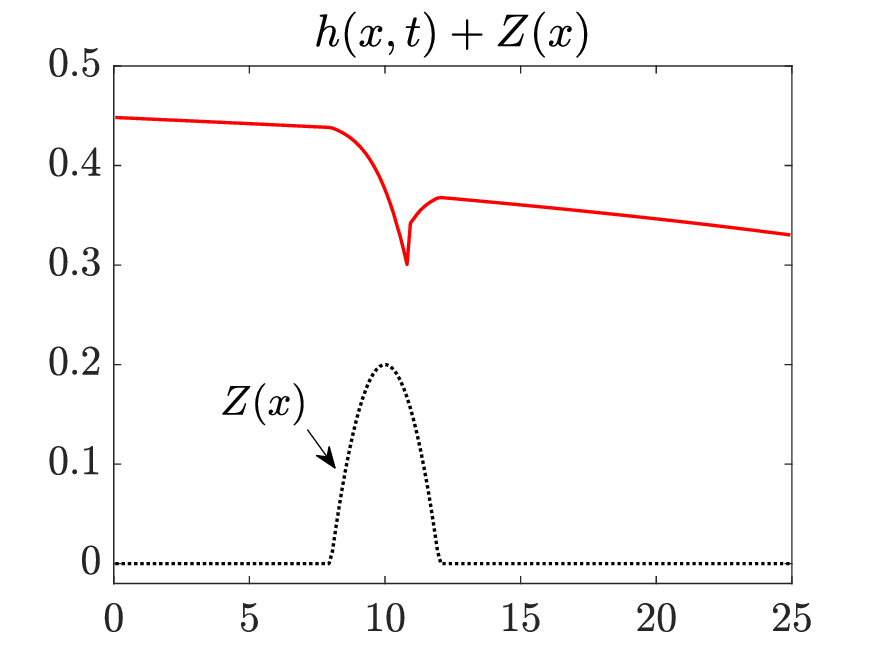}\hspace*{0.15cm}
\includegraphics[trim=0cm 0.25cm 0.9cm 0.15cm,clip,width=4.2cm]{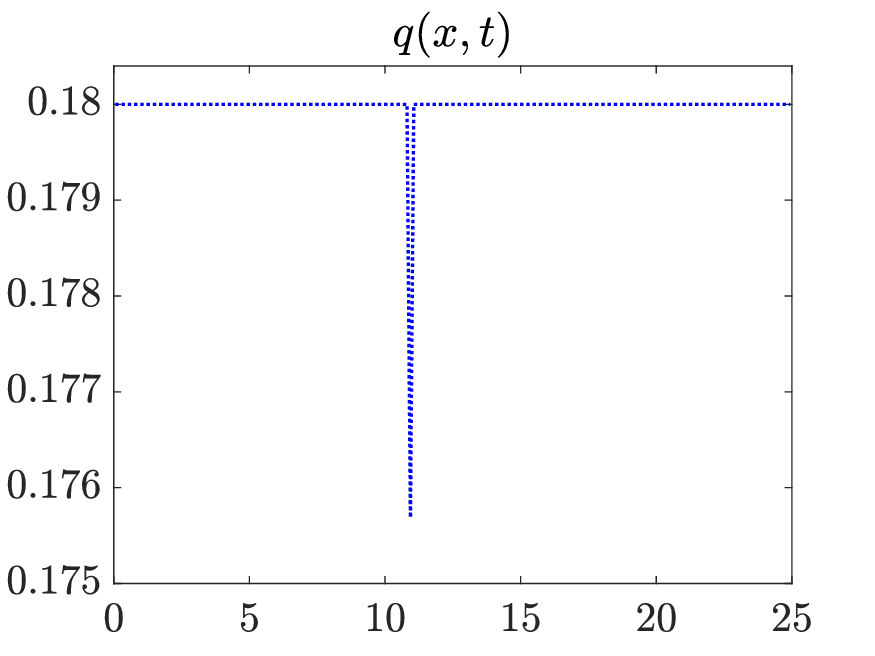}\hspace*{0.15cm}
\includegraphics[trim=0cm 0.25cm 0.9cm 0.15cm,clip,width=4.2cm]{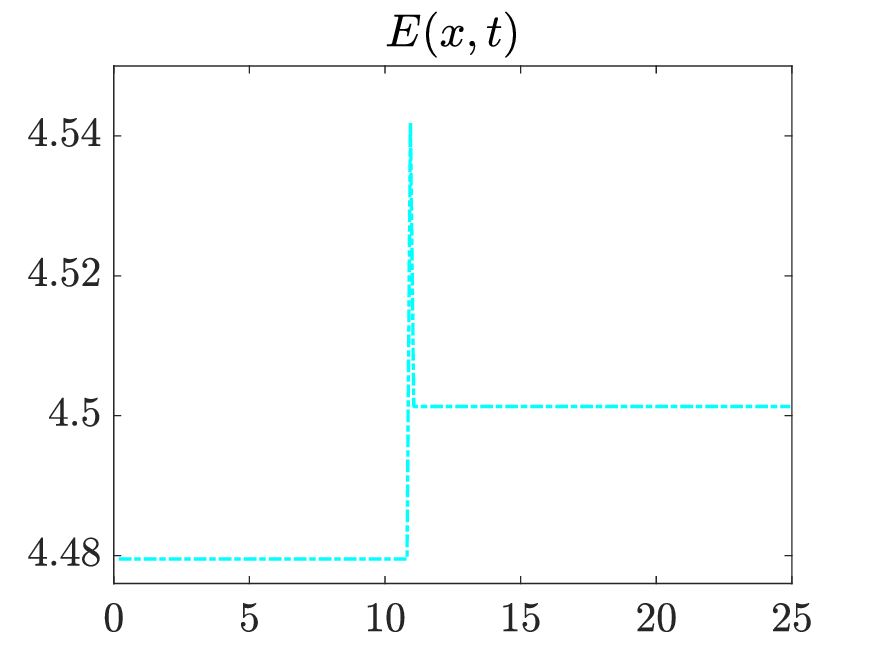}}
\caption{\sf Example 7: Convergent solutions of $h+Z$ together with $Z$ (left column), $q$ (middle column), and $E$ (right column) for Cases (a, first row), (b, middle row), and (d, bottom row) with $n=0.05$ and the deviation number $\epsilon=10^{-12}$.\label{fig6}}
\end{figure}

As in Example 4, we compute the numerical solutions using the proposed \bla{WB AF} scheme at time $t=500$ in the computational domain $[0,25]$ covered with $200$ uniform cells. We plot the obtained numerical solutions ($h+Z$, $q$, and $E$) in Figure \ref{fig6}. As one can see, the obtained numerical results are comparable with those reported in \cite{LCJKY, CCHKW_18}. We can therefore conclude that the proposed \bla{WB AF} scheme is also able to capture the steady states of different flow regimes in the presence of nonzero friction terms. 

\subsection*{Example 8---Small Perturbation of Moving-Water Equilibria for Frictional case}
In the \bla{eighth} example, we test the ability of the proposed \bla{WB AF} scheme to capture small perturbations of the obtained discrete steady states of Cases (a) and (b) as shown in Example 7. Similar to Example 3, we add a small perturbation $10^{-3}e^{-80(x-6)^2}$ to the discrete steady-state cell average of water depth obtained from the following two sets of the discrete moving-water equilibria:
\begin{equation}\label{3.4}
   \begin{aligned}
&\mbox{Case (a):}&&q(x,0)\equiv 24,\quad &&E(x,0)\equiv 91.624;\\
&\mbox{Case (b):}&&q(x,0)\equiv 4.42,\quad &&E(x,0)\equiv 23.179263521617973.
\end{aligned}
\end{equation}

We note that the initial data in \eref{3.4} is given in terms of the equilibrium variable $\bm E=(q, E)^\top$ rather than in $h$ and $u$. However, in order to start the computation at time $t=0$, one has to obtain the average values of $\xbar h_\jph$ and the point values of $h_j$. Unlike Example 3, here we use the Newton's method to numerically solve the following nonlinear equations:
\begin{equation*}
  E_j=\frac{q_j^2}{2h_j^2}+g(h_j+Z_j+Q_j),\quad E_\jph=\frac{q^2_\jph}{2h^2_\jph}+g(h_\jph+Z_\jph+Q_\jph),\quad \forall j,
\end{equation*}
where $Q_\jph$ and $\{Q_j\}_{j>1}$ are computed from the recursive formulae \eref{2.16Q1} and \eref{2.16Q2}. After obtaining the point values of water depth, it is easy to compute the initial average values of $\xbar h_\jph$ using Simpson's rule.

\begin{figure}[ht!]
\centerline{\includegraphics[trim=0.9cm 0.3cm 0.6cm 0.2cm,clip,width=5.5cm]{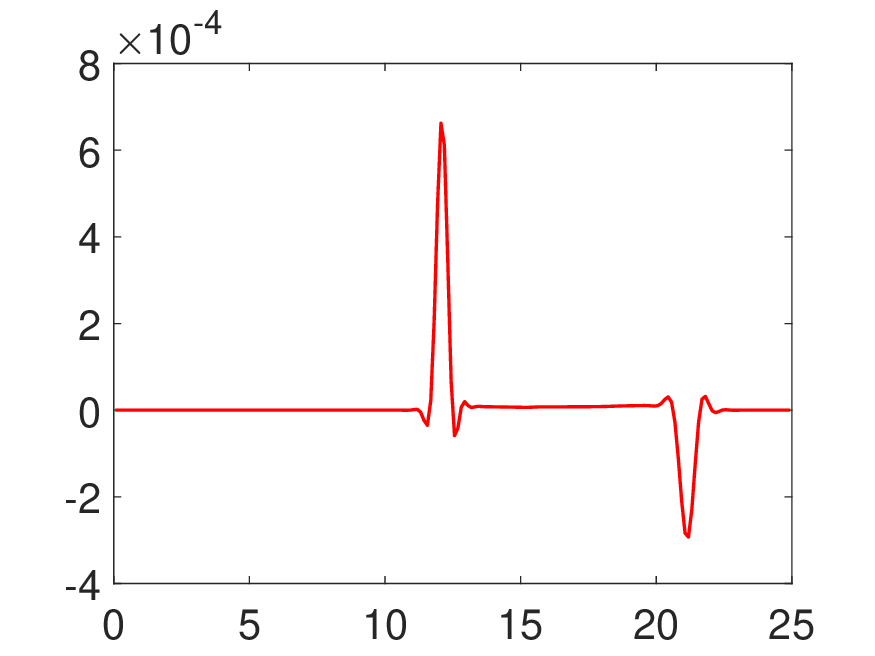}\hspace*{0.3cm}
\includegraphics[trim=0.9cm 0.3cm 0.6cm 0.2cm,clip,width=5.5cm]{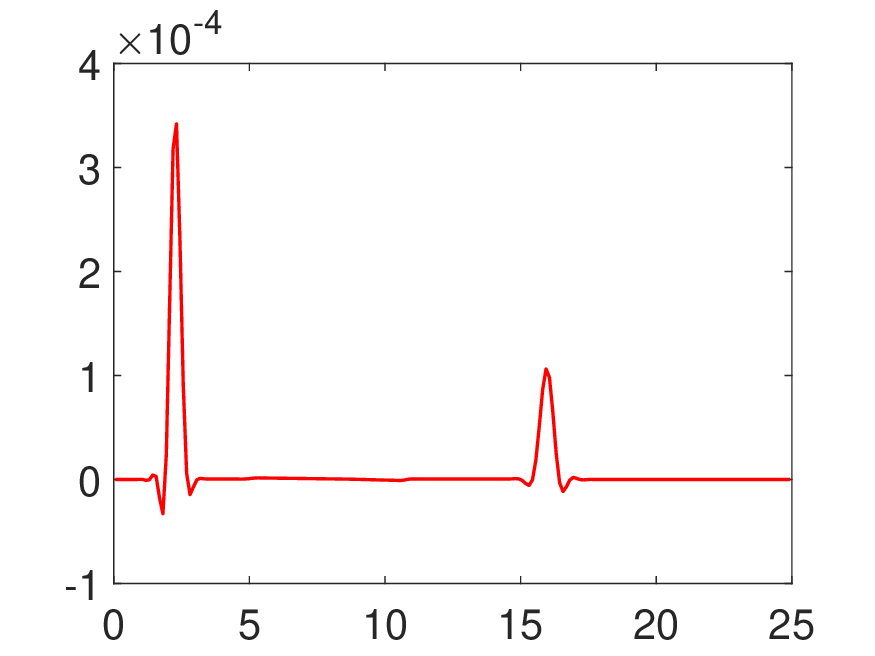}}
\caption{\sf Example 8: The difference between the computed and the background discrete steady-state cell averages of water depth computed using 200 uniform cells for Cases (a, left) and (b, right).\label{fig62}}
\end{figure}

Equipped with the initial values in terms of water depth and velocity, we now compute the numerical solutions using the \bla{WB AF} scheme until $t=1$ and $t=1.5$ for Cases (a) and (b), respectively. The difference between the computed and the background discrete moving steady-state cell averages of water depth are plotted in Figure \ref{fig62}, where one can see that the proposed \bla{WB AF} scheme is able to accurately capture the small perturbation in a robust manner.

\subsection*{Example 9---Conservation test}
\begin{figure}[ht!]
\centerline{\includegraphics[trim=0.5cm 0.1cm 0.6cm 0.1cm,clip,width=6.0cm]{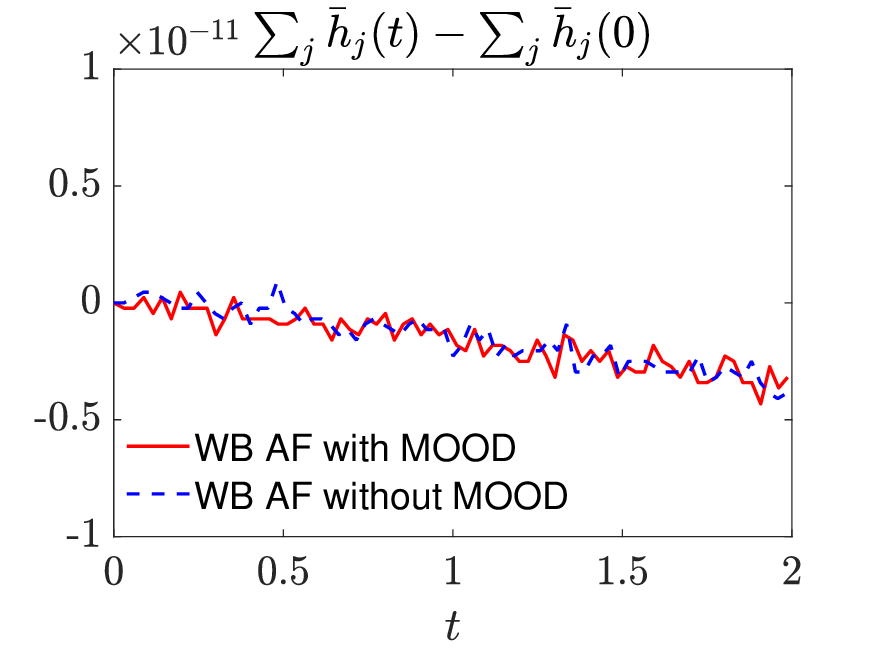}\hspace*{0.3cm}
\includegraphics[trim=0.5cm 0.1cm 0.6cm 0.1cm,clip,width=6.0cm]{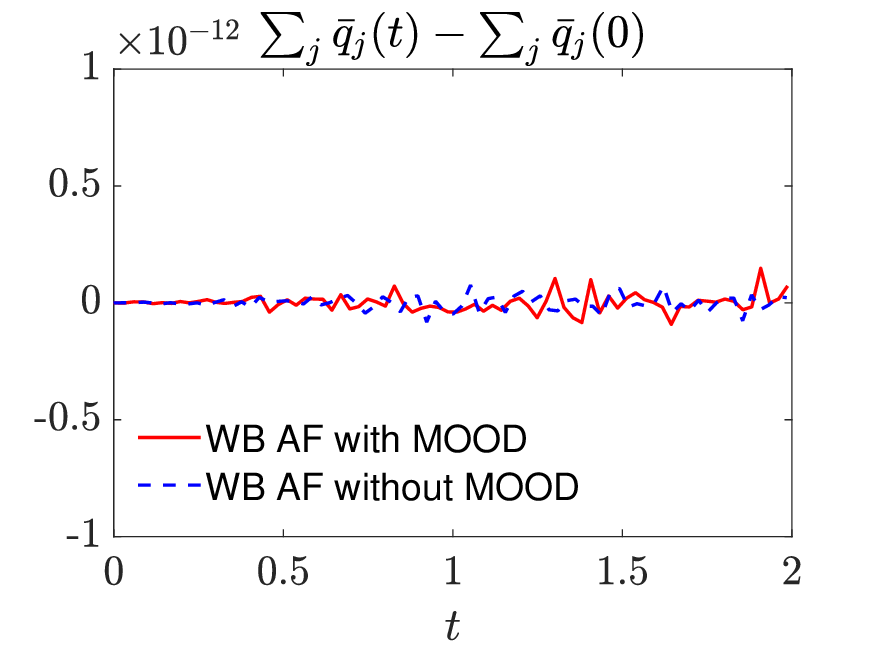}}
\caption{\sf Example 9: Plot of the deviations in the discrete mass and discharge over time up to $t=2$.\label{fig9}}
\end{figure}

In the final example, we consider a flat bottom topography $(Z(x)\equiv0)$ and zero Manning friction ($n=0$). It is known that the scheme should conserve mass and momentum. To test the conservation property, we consider the following initial data
\begin{equation*}
  h(x,0)=\left\{\begin{aligned}
  &5\quad &&\mbox{if}~x\in[-5,5],\\
  &10\quad &&\mbox{otherwise},
  \end{aligned}\right.\quad u(x,0)=0,
\end{equation*}
prescribed in a computational domain $[-100,100]$. We set the periodic boundary condition and use the WB AF scheme with/without MOOD activation to compute the numerical solution until time $t=2$ on a uniform mesh with $200$ cells. The deviations in the discrete mass and discharge comparing with the initial values at every small time step are shown in Figure \ref{fig9}. We can see that the deviations are on the order of machine accuracy and hence the conservation property is verified.

\section{Conclusion}\label{sec4}
In this paper, we have developed and tested new positivity-preserving well-balanced schemes for the one-dimensional (1-D) Saint-Venant system of shallow water equations with the Manning friction term. In order to construct an accurate, robust, and fully well-balanced scheme for the 1-D shallow water equations, it is essential to define a transformation from the equilibrium variables to the conservative ones for many existing numerical methods; see, e.g., \cite{CKLX,CKLLW,CTR,KKLZ,CCHKW_18,CK16,Xing14, NXS}. In the process of transformation, one needs to solve some nonlinear equations and this can be cumbersome. In this work, we consider both the conservative and primitive forms of the studied hyperbolic system and are able to achieve the fully well-balanced property without solving any nonlinear equations. More precisely, we use the conservative form to evolve the cell averages while use the primitive one to evolve the point values. It is noticed that using the introduced equilibrium variables, the primitive formulation will reduce to a very simple form with flux function only represented by the equilibrium variables. We then apply a stable finite difference method to compute the residuals, which will vanish when the solution is at the steady state. In order to update the average values in a WB manner, we use a well-balanced fourth-order approximation from \cite{NXS} to evaluate the average value of the source term. We also propose a way to guarantee the positivity property of water depth.

This study is only concerned with the one-dimensional shallow water equations but provides a new framework that does not require the transformation from equilibrium variables to conservative variables for designing fully WB methods. In future work, we plan to extend the proposed method to several space dimensions and other hyperbolic models including the thermal rotating shallow water equations, shallow water flows in channels, and the blood flow models.

\appendix
\section{A third-order parabolic interpolant}\label{appb} 
In this appendix, we describe a third-order parabolic interpolant which is used in \eref{2.15x}. Given the point value $\bm U_j$ and the cell average $\xbar{\bm U}_\jph$, we can construct a third-order parabolic interpolant:
\begin{equation}\label{2.12}
  \widetilde{\bm U}(x)=\bm U_j\ell_0(\xi)+\xbar{\bm U}_\jph\ell_{\frac{1}{2}}(\xi)+\bm U_{j+1}\ell_1(\xi),\quad \xi=\frac{x-x_j}{\dx}\in[0,1],
\end{equation}
where 
\begin{equation*}
  \ell_0(\xi)=(1-\xi)(1-3\xi),\quad \ell_1(\xi)=\xi(3\xi-2),\quad \ell_{\frac{1}{2}}(\xi)=6\xi(1-\xi).
\end{equation*}
\section{Desingularization}\label{appc}
In this appendix, we describe the desingularization function, which we use whenever a division by zero or by a very small positive number needs to be avoided. Assume that we need to evaluate the quotient $a/b$ and that $b\sim0$. We then use the simplest desingularization from \cite{BCKN} to compute the quotient as
\begin{equation}\label{Desing}
  \mathfrak{D}(a,b)=\left\{\begin{aligned}
  &\frac{a}{b}\quad &&\mbox{if}~b\geq\varepsilon,\\
  &0\quad &&\mbox{otherwise}.
  \end{aligned}\right.
\end{equation}
where $\varepsilon$ is small desingularization parameter which is set to $\varepsilon=10^{-9}$ in all of the reported numerical examples. Alternative desingularization functions can also be employed, the proposed schemes are not tight to these selections.

\section{Numerical tests over discontinuous bathymetries}\label{appe}
\begin{figure}[ht!]
\centerline{\includegraphics[trim=0.8cm 0.25cm 0.9cm 0.2cm,clip,width=4.0cm]{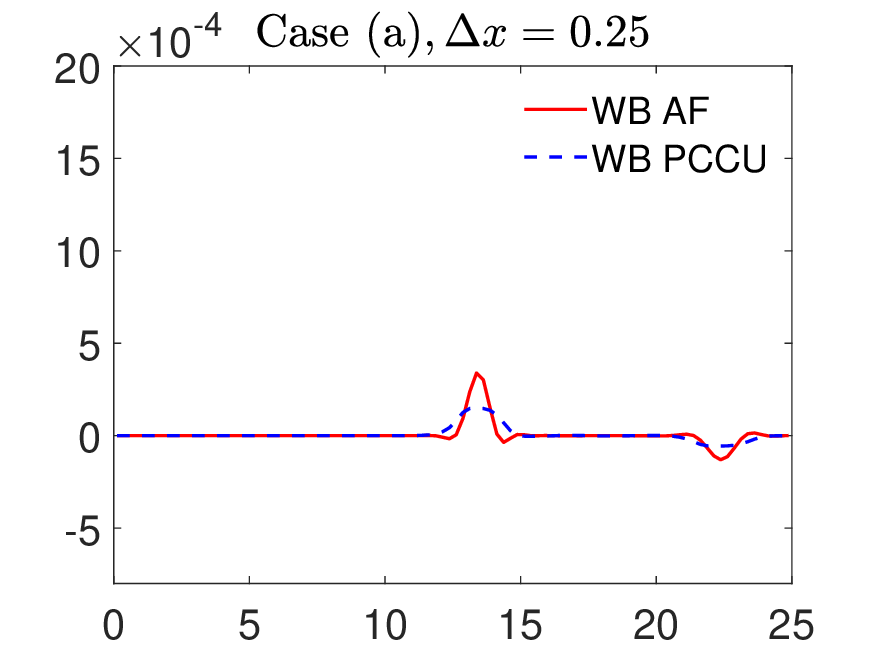}\hspace*{0.15cm}
\includegraphics[trim=0.8cm 0.25cm 0.9cm 0.2cm,clip,width=4.0cm]{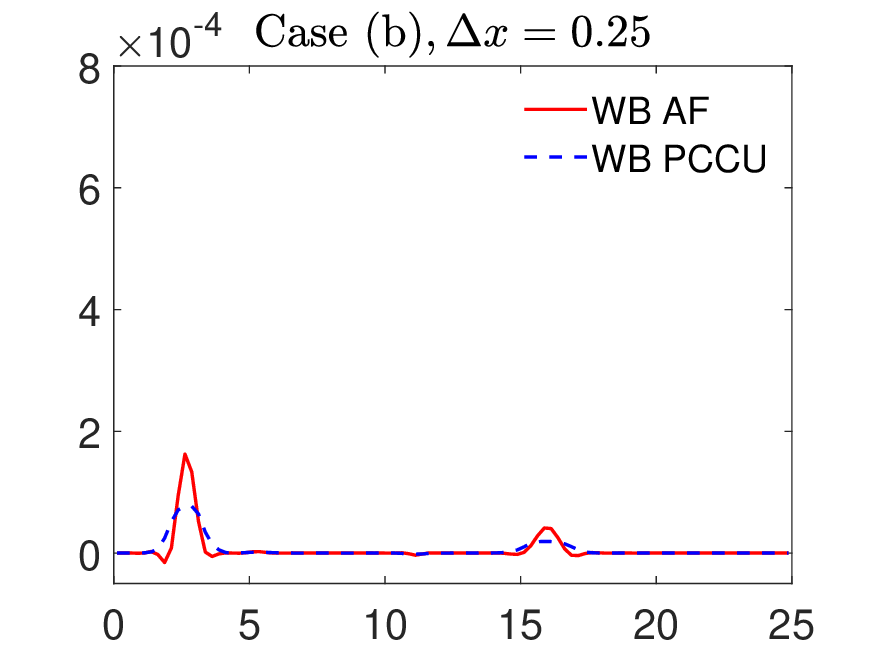}\hspace*{0.15cm}
\includegraphics[trim=0.8cm 0.25cm 0.9cm 0.2cm,clip,width=4.0cm]{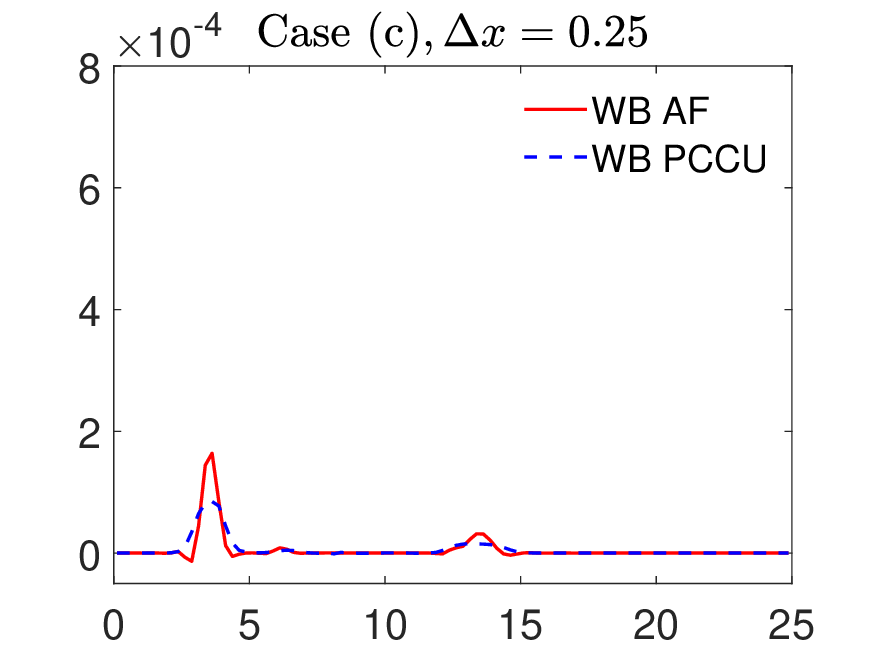}}
\vskip5pt
\centerline{\includegraphics[trim=0.8cm 0.25cm 0.9cm 0.2cm,clip,width=4.0cm]{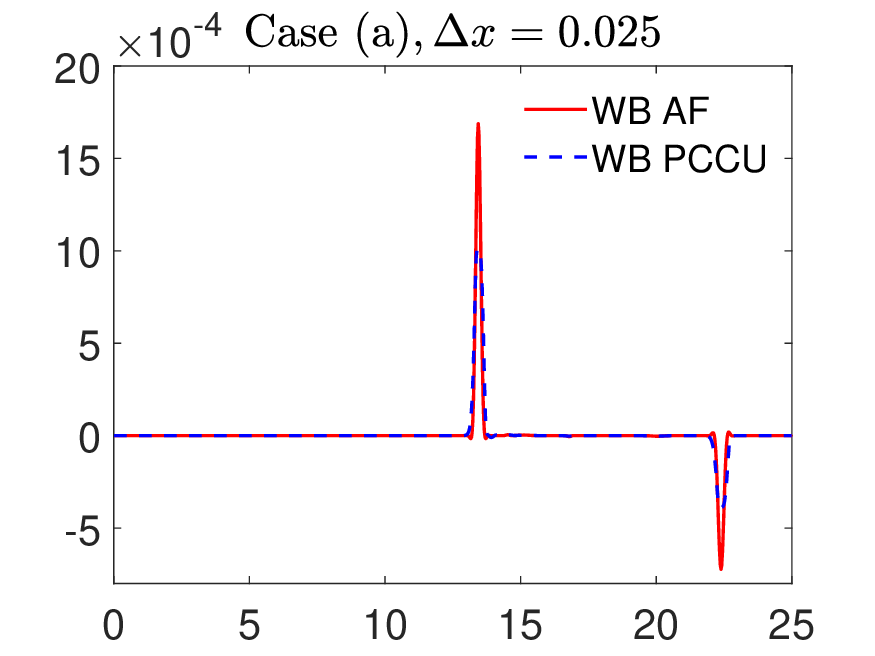}\hspace*{0.15cm}
\includegraphics[trim=0.8cm 0.25cm 0.9cm 0.2cm,clip,width=4.0cm]{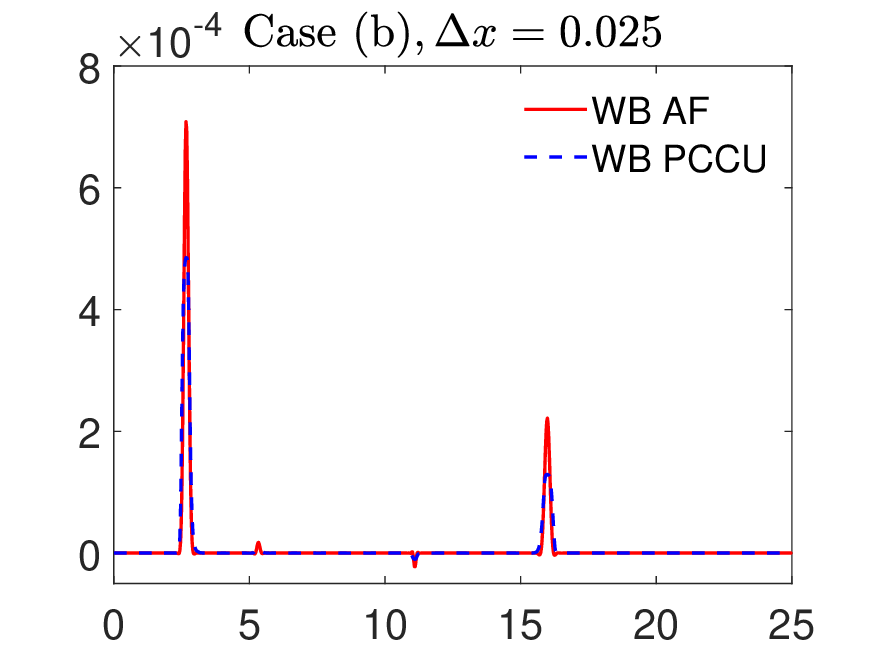}\hspace*{0.15cm}
\includegraphics[trim=0.8cm 0.25cm 0.9cm 0.2cm,clip,width=4.0cm]{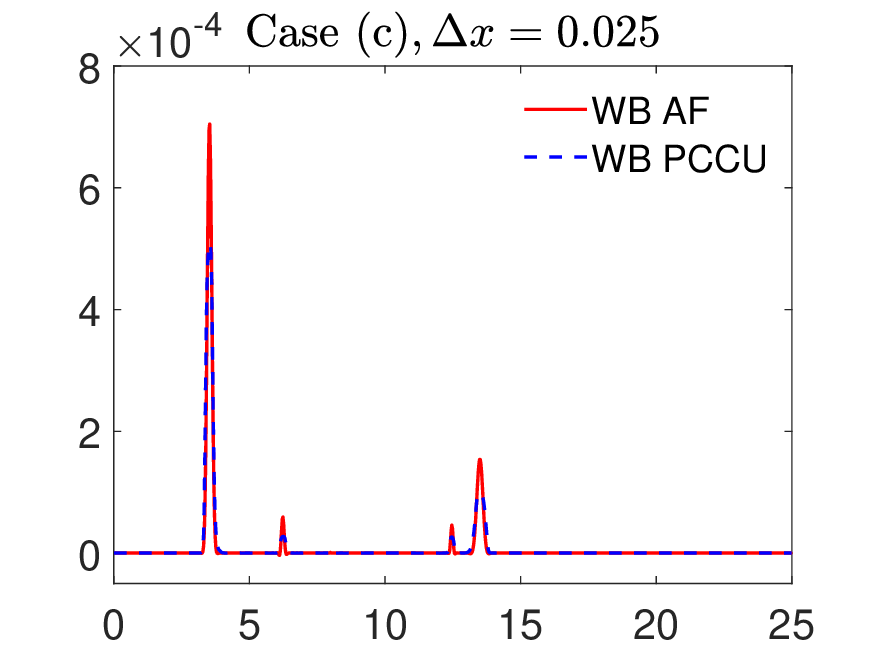}}
\caption{\sf Small perturbation: same as in Figure \ref{fig2} but over discontinuous bathymetry \eref{zz2}.\label{fig2a}}
\end{figure}

\paragraph{Test 1} In this test, we consider the same initial setting as in Example 3 but over a discontinuous bathymetry:
\begin{equation}\label{zz2}
Z(x)=\left\{\begin{aligned}
&0.2&&\mbox{if}~8\le x\le12,\\
&0&&\mbox{otherwise}.
\end{aligned}\right.
\end{equation}
In this case, the relevant steady-state water depth $h_{\rm eq}(x,0)$ for case (c) in \eref{3.3} becomes
\begin{equation*}
h_{\rm eq}(x,0)=\left\{\begin{aligned}
&-\frac{a_0}{3}\Big[2\cos\big(\frac{\theta+4\pi}{3}\big)+1\Big]&&\quad \mbox{if}~x>12,\\
&-\frac{a_0}{3}\Big[2\cos\frac{\theta}{3}+1\Big]&&\quad \mbox{if}~x<8,\\
&\Big(\frac{q_{\rm eq}^2}{g}\Big)^{\frac{1}{3}}&&\quad \mbox{else}.
\end{aligned}\right.
\end{equation*}
The same simulation as it was done in Example 3 is performed and we report the results in Table \ref{tab3a} and Figure \ref{fig2a}. As one can see, the proposed scheme is able to exactly preserve the moving-water equilibrium and accurately capture the small perturbations even over a discontinuous bathymetry.
\begin{table}[ht!]
\caption{\sf $L^1$- and $L^\infty$-errors in average values of $h$ and $q$ over discontinuous bathymetry \eref{zz2}.\label{tab3a}}
\begin{center}
\begin{tabular}{|c|c|c|c|c|}
\hline
&$L^1$-error in $h$&$L^\infty$-error in $h$&$L^1$-error in $q$&$L^\infty$-error in $q$\\
\hline
Case (a)&3.57e-13&3.06e-14&3.79e-12&3.23e-13\\
\hline
Case (b)&$<$1.0e-16&$<$1.0e-16&$<$1.0e-16&$<$1.0e-16\\
\hline
Case (c)&3.16e-15&2.22e-16&5.55e-15&2.22e-16\\
\hline
\end{tabular}
\end{center}
\end{table}

\begin{figure}[ht!]
\centerline{\includegraphics[trim=0cm 0.25cm 0.9cm 0.12cm,clip,width=4.2cm]{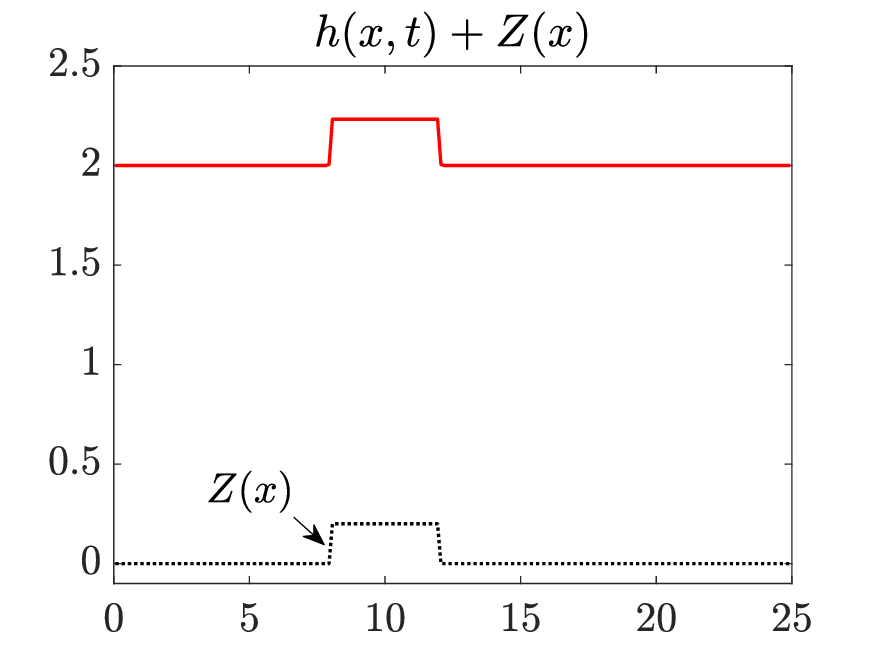}\hspace*{0.15cm}
\includegraphics[trim=0cm 0.25cm 0.9cm 0.12cm,clip,width=4.2cm]{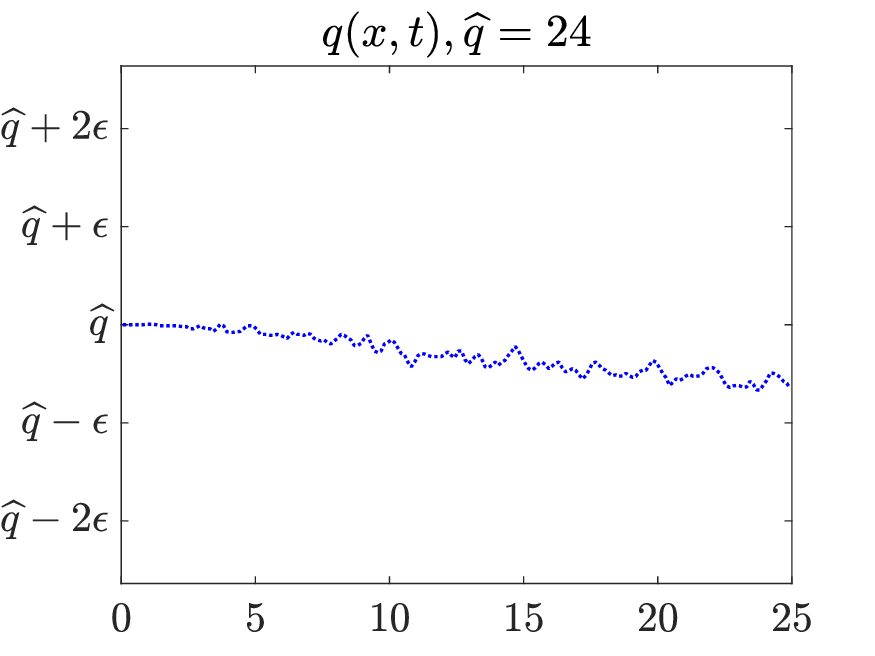}\hspace*{0.15cm}
\includegraphics[trim=0cm 0.25cm 0.9cm 0.12cm,clip,width=4.2cm]{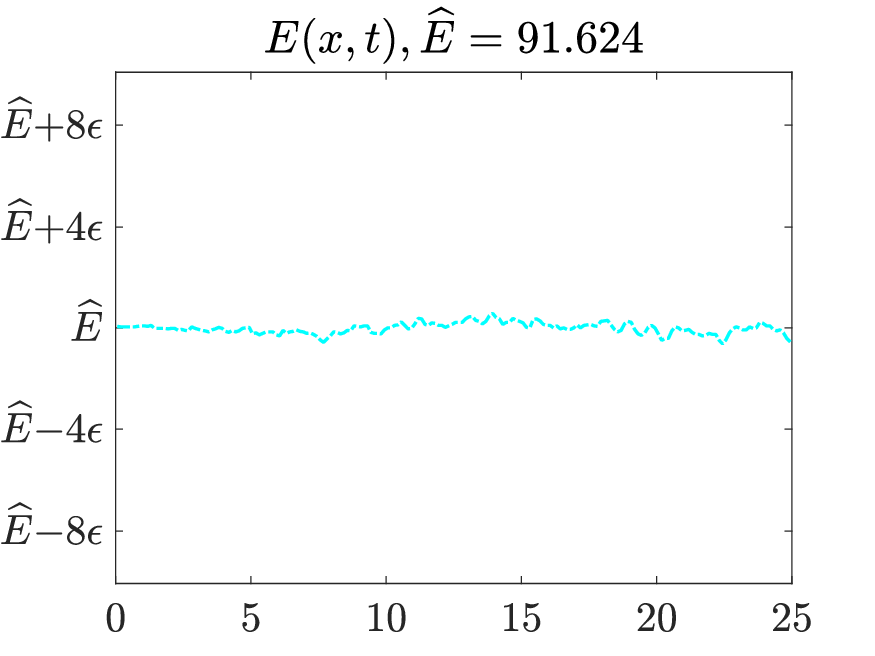}}
\vskip5pt
\centerline{\includegraphics[trim=0cm 0.25cm 0.9cm 0.12cm,clip,width=4.2cm]{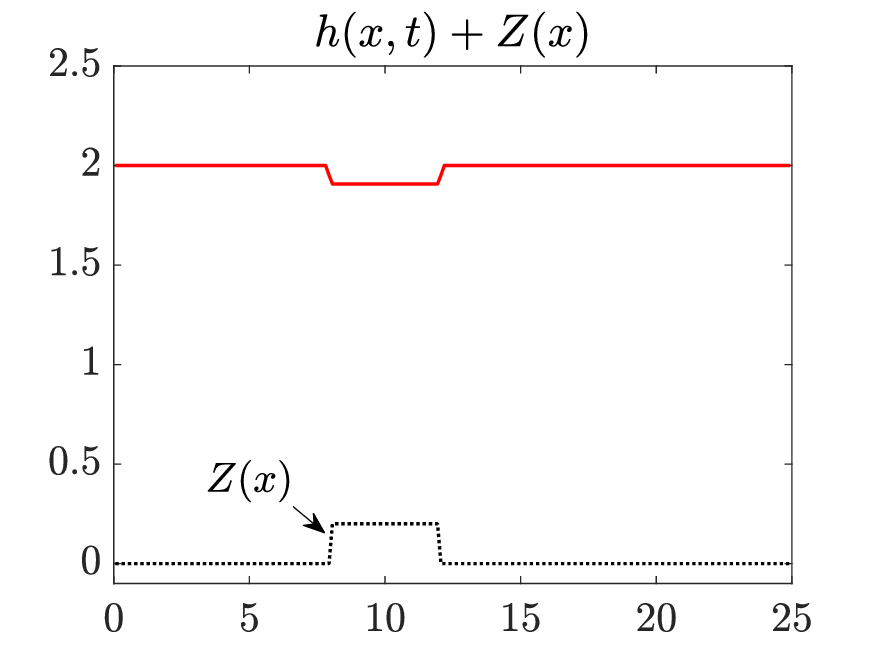}\hspace*{0.15cm}
\includegraphics[trim=0cm 0.25cm 0.9cm 0.12cm,clip,width=4.2cm]{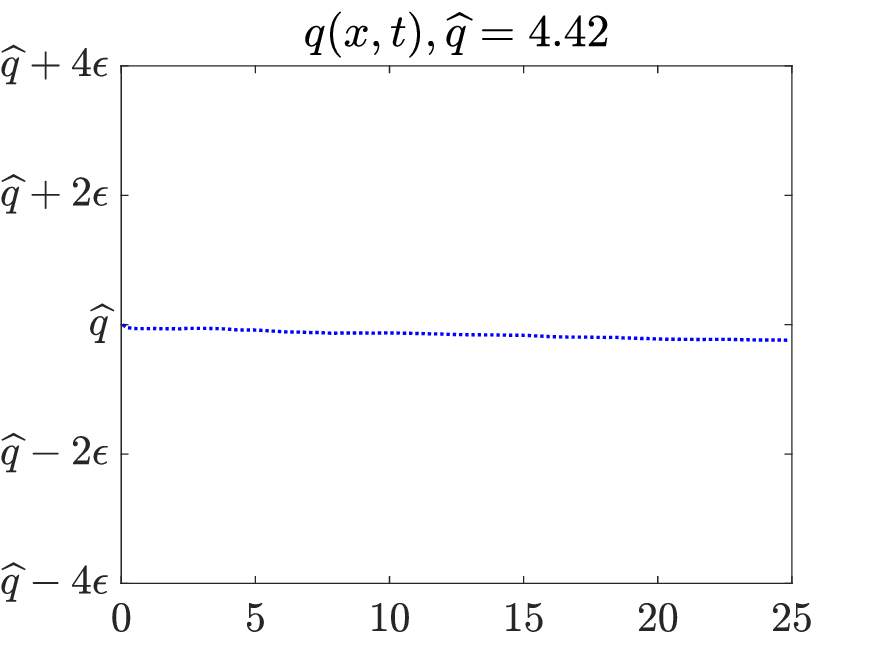}\hspace*{0.15cm}
\includegraphics[trim=0cm 0.25cm 0.9cm 0.12cm,clip,width=4.2cm]{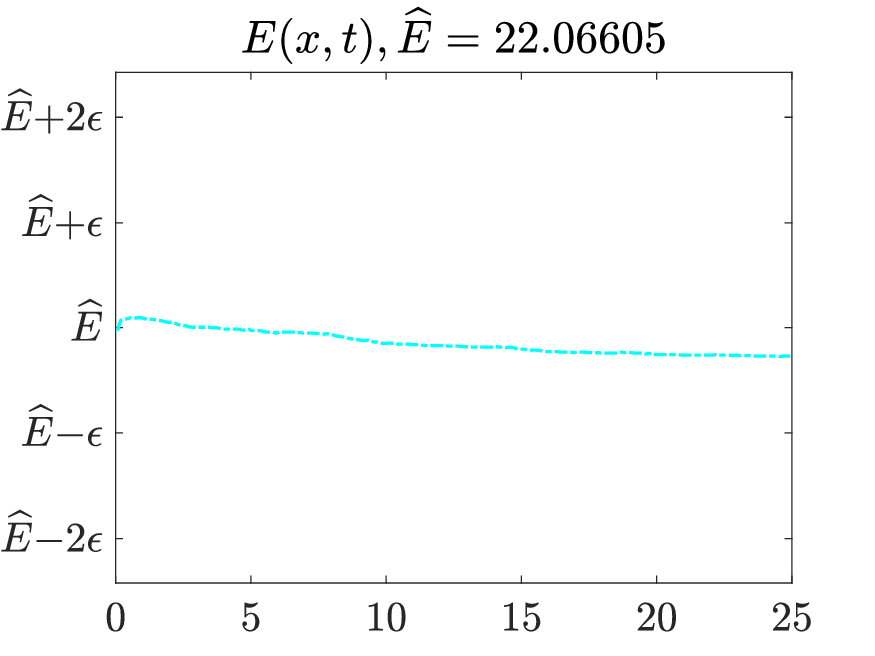}}
\vskip5pt
\centerline{\includegraphics[trim=0cm 0.25cm 0.9cm 0.12cm,clip,width=4.2cm]{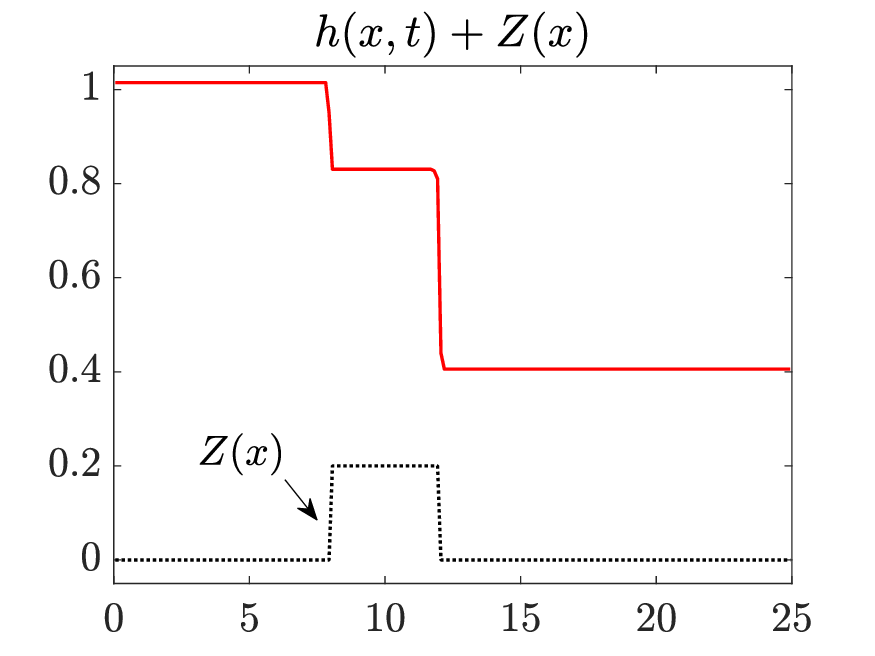}\hspace*{0.15cm}
\includegraphics[trim=0cm 0.25cm 0.9cm 0.12cm,clip,width=4.2cm]{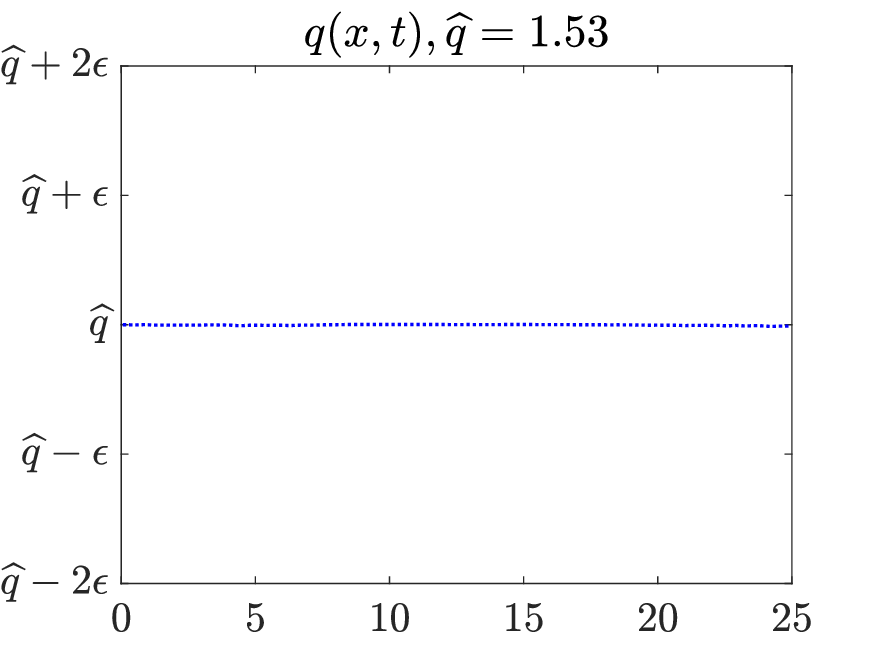}\hspace*{0.15cm}
\includegraphics[trim=0cm 0.25cm 0.9cm 0.12cm,clip,width=4.2cm]{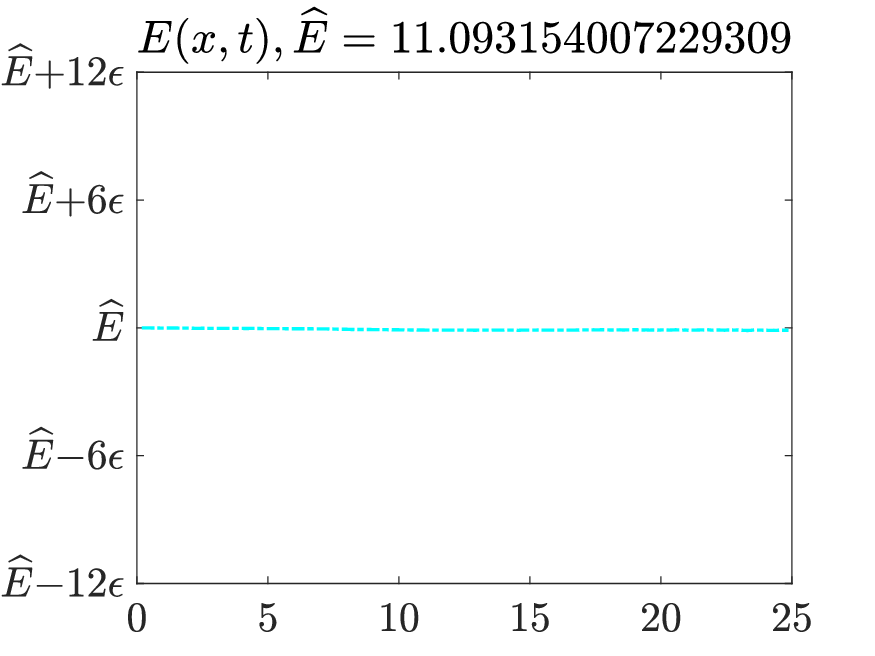}}
\vskip5pt
\centerline{\includegraphics[trim=0cm 0.25cm 0.9cm 0.12cm,clip,width=4.2cm]{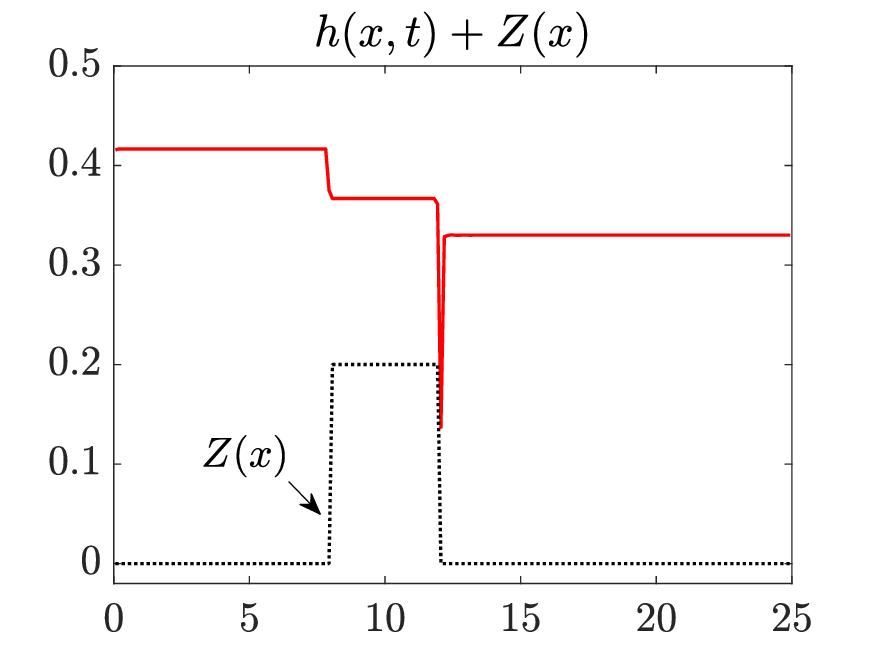}\hspace*{0.15cm}
\includegraphics[trim=0cm 0.25cm 0.9cm 0.12cm,clip,width=4.2cm]{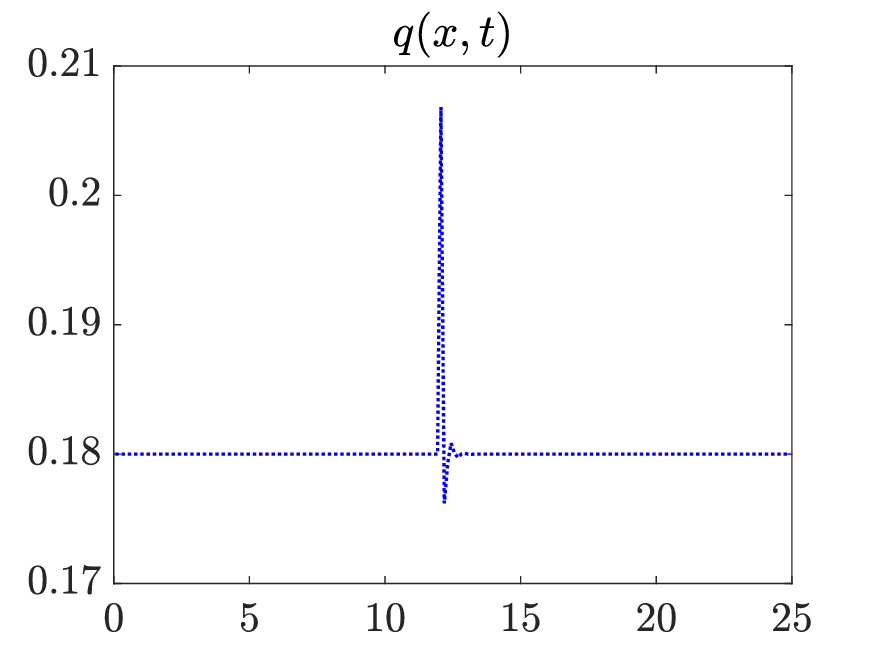}\hspace*{0.15cm}
\includegraphics[trim=0cm 0.25cm 0.9cm 0.12cm,clip,width=4.2cm]{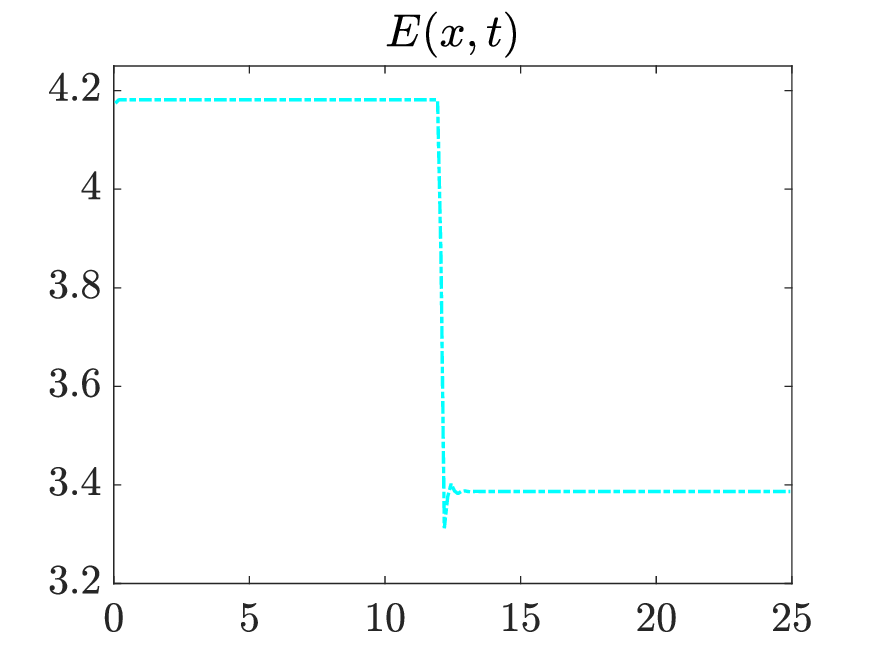}}
\caption{\sf Convergent solutions: Same as in Figure \ref{fig4} but over discontinuous bathemetry \eqref{zz2}. \label{fig4db}}
\end{figure}

\paragraph{Test 2} In this test, we repeat the simulation as it was done in Example 4, but with discontinuous bathymetry given by \eref{zz2}. The obtained results are reported in Figure \ref{fig4db}. As we can see, the expected convergent solutions are obtained and consistent with those produced by the well-balanced path-conservative central-upwind scheme in \cite{CKLX}.

\begin{figure}[ht!]
\centerline{\includegraphics[trim=0.9cm 0.1cm 0.6cm 0.2cm,clip,width=5.2cm]{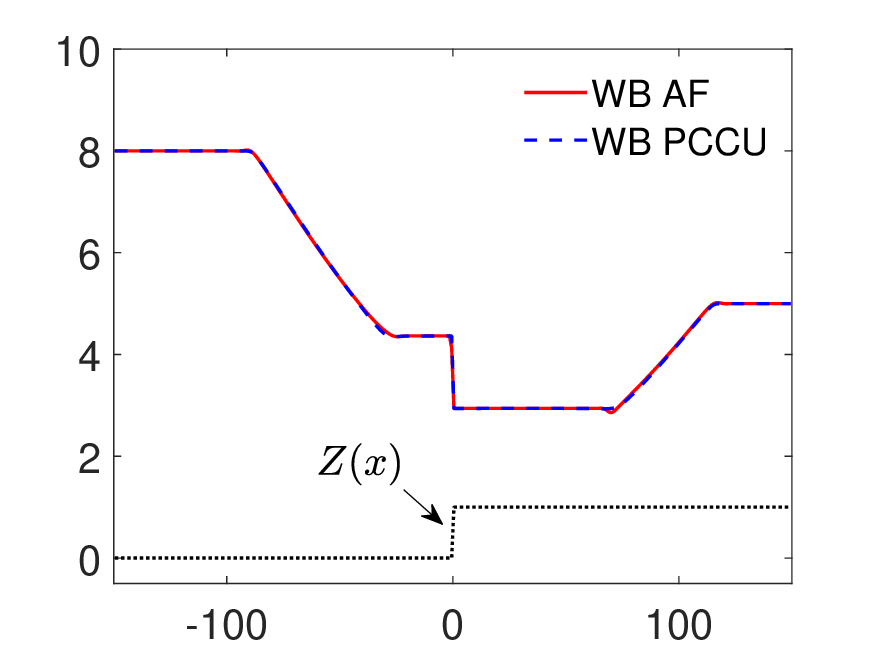}\hspace*{0.3cm}
\includegraphics[trim=0.9cm 0.1cm 0.6cm 0.2cm,clip,width=5.2cm]{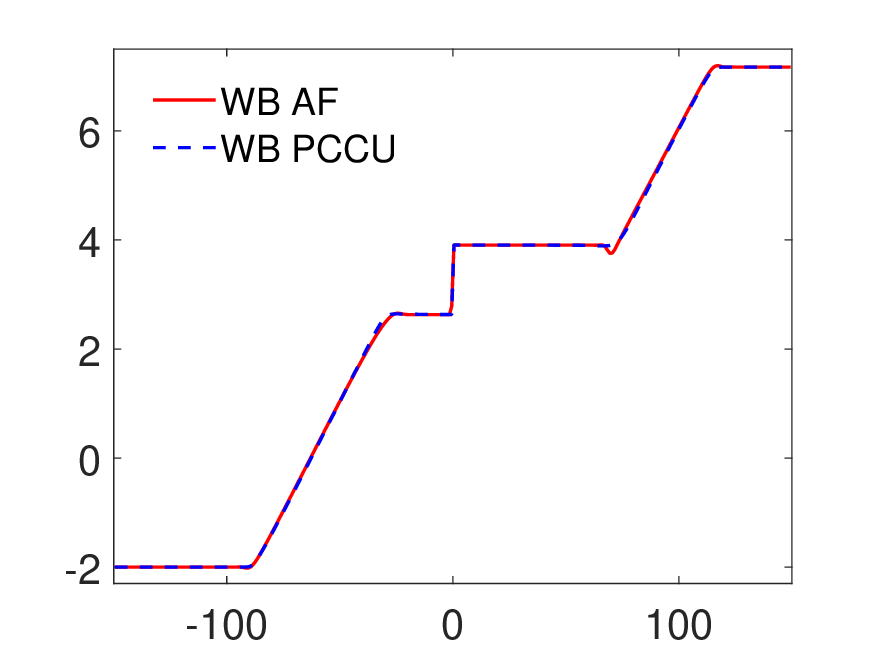}}
\caption{\sf Water surface $h+Z$ together with $Z$ (left) and velocity $u$ (right) computed by the \bla{WB AF} and \bla{WB PCCU} schemes.\label{Ex52fig}}
\end{figure}

\paragraph{Test 3} In this test, we study a problem associated with two rarefaction waves moving toward opposite directions. The initial data and bottom topography are
\begin{equation*}
  (h(x,0), u(x,0), Z(x))=\left\{\begin{aligned}
  &(8, -2, 0),\quad &&x<0,\\
  &(5, 7.1704, 1),\quad &&x\geq0.\\
  \end{aligned}\right.
\end{equation*}
The exact solution of this problem is given by a left-propagating 1-Rarefaction, a bottom step discontinuity, and a right-propagating 2-Rarefaction. We use both the \bla{WB AF} and \bla{WB PCCU} schemes to compute the numerical solutions at the final time $t=8$ in the computational domain $[-150,150]$ covered with $300$ uniform cells. Results are shown in Figure \ref{Ex52fig}. As one can see, the obtained water surface $h+Z$ profiles provide expected results with a left-directed rarefaction wave, a stationary shock, and a right-directed rarefaction wave. The numerical solutions computed by the two studied schemes are in a satisfactory agreement.

\bibliographystyle{siamplain}
\bibliography{reference}

\begin{thebibliography}{10}

\bibitem{Abgrall_camc}
{\sc R.~Abgrall}, {\em A combination of {R}esidual {D}istribution and the
  {A}ctive {F}lux formulations or a new class of schemes that can combine
  several writings of the same hyperbolic problem: application to the 1{D}
  {E}uler equation}, Commun. Appl. Math. Comput., 5 (2023), pp.~370--402.

\bibitem{AT}
{\sc R.~Abgrall and D.~Torlo}, {\em Some preliminary results on a high order
  asymptotic preserving computationally explicit kinetic scheme}, Commun. Math.
  Sci., 20 (2022), pp.~297--326.

\bibitem{ABBKP}
{\sc E.~Audusse, F.~Bouchut, M.-O. Bristeau, R.~Klein, and B.~Perthame}, {\em A
  fast and stable well-balanced scheme with hydrostatic reconstruction for
  shallow water flows}, SIAM J. Sci. Comput., 25 (2004), pp.~2050--2065.

\bibitem{Audusse15}
{\sc E.~Audusse, C.~Chalons, and P.~Ung}, {\em A simple well-balanced and
  positive numerical scheme for the shallow-water system}, Commun. Math. Sci.,
  13 (2015), pp.~1317--1332.

\bibitem{BB}
{\sc W.~Barsukow and J.~P. Berberich}, {\em A well-balanced {A}ctive {F}lux
  method for the shallow watere quations with wetting and drying}, Commun.
  Appl. Math. Comput.,  (2023).
\newblock https://doi.org/10.1007/s42967-022-00241-x.

\bibitem{BC_FWB}
{\sc C.~Berthon and C.~Chalons}, {\em A fully well-balanced, positive and
  entropy-satisfying {G}odunov-type method for the shallow-water equations},
  Mathematics of Computation, 85 (2016), pp.~1281--1307.

\bibitem{BerthonF}
{\sc C.~Berthon and F.~Foucher}, {\em Efficient well-balanced hydrostatic
  upwind schemes for shallow-water equations}, J. Comput. Phys., 231 (2012),
  pp.~4993--5015.

\bibitem{BMLS}
{\sc C.~Berthon, M.~M'Baye, M.~H. Le, and D.~Seck}, {\em A well-defined moving
  steady states capturing {G}odunov-type scheme for shallow-water model}, Int.
  J. Finite Vol., 15 (2021).

\bibitem{BerthonM_23}
{\sc C.~Berthon and V.~Michel-Dansac}, {\em A fully well-balanced hydrodynamic
  reconstruction},  (2023).
\newblock https://hal.science/hal-04083181.

\bibitem{BCKN}
{\sc A.~Bollermann, G.~Chen, A.~Kurganov, and S.~Noelle}, {\em A well-balanced
  reconstruction of wet/dry fronts for the shallow water equations}, J. Sci.
  Comput., 56 (2013), pp.~267--290.

\bibitem{BNL}
{\sc A.~Bollermann, S.~Noelle, and M.~Luk\'a\v{c}ov\'a-Medvi\v{d}ov\'a}, {\em
  Finite volume evolution {G}alerkin methods for the shallow water equations
  with dry beds}, Commun. Comput. Phys., 10 (2011), pp.~371--404.

\bibitem{CSDT}
{\sc A.~Canestrelli, A.~Siviglia, M.~Dumbser, and E.~F. Toro}, {\em
  Well-balanced high-order centred schemes for non-conservative hyperbolic
  systems. applications to shallow water equations with fixed and mobile bed},
  Adv. Water Resour., 32 (2009), pp.~834--844.

\bibitem{CKLX}
{\sc Y.~Cao, A.~Kurganov, Y.~Liu, and R.~Xin}, {\em Flux globalization based
  well-balanced path-conservative central-upwind schemes for shallow water
  models}, J. Sci. Comput., 92 (2022).
\newblock Paper No. 69, 31 pp.

\bibitem{CLPwb13}
{\sc M.~J. Castro~D\'{\i}az, J.~A. L\'{o}pez-Garc\'{\i}a, and C.~Par\'{e}s},
  {\em High order exactly well-balanced numerical methods for shallow water
  systems}, J. Comput. Phys., 246 (2013), pp.~242--264.

\bibitem{CPMP07}
{\sc M.~J. Castro~D\'{\i}az, A.~Pardo~Milan{\'e}s, and C.~Par{\'e}s}, {\em
  Well-balanced numerical schemes based on a generalized hydrostatic
  reconstruction technique}, Math. Models Methods Appl. Sci., 17 (2007),
  pp.~2055--2113.

\bibitem{CCHKW_18}
{\sc Y.~Cheng, A.~Chertock, M.~Herty, A.~Kurganov, and T.~Wu}, {\em A new
  approach for designing moving-water equilibria preserving schemes for the
  shallow water equations}, J. Sci. Comput., 80 (2019), pp.~538--554.

\bibitem{CK16}
{\sc Y.~Cheng and A.~Kurganov}, {\em Moving-water equilibria preserving
  central-upwind schemes for the shallow water equations}, Commun. Math. Sci.,
  14 (2016), pp.~1643--1663.

\bibitem{CCKW}
{\sc A.~Chertock, S.~Cui, A.~Kurganov, and T.~Wu}, {\em Well-balanced
  positivity preserving central-upwind scheme for the shallow water system with
  friction terms}, Int. J. Numer. Methods Fluids, 78 (2015), pp.~355--383.

\bibitem{CKLLW}
{\sc A.~Chertock, A.~Kurganov, X.~Liu, Y.~Liu, and T.~Wu}, {\em Well-balancing
  via flux globalization: {A}pplications to shallow water equations with
  wet/dry fronts}, J. Sci. Comput., 90 (2022).
\newblock Paper No. 9, 21 pp.

\bibitem{CMMOT}
{\sc M.~Ciallella, L.~Micalizzi, V.~Michel-Dansac, P.~{\"O}ffner, and
  D.~Torlo}, {\em A high-order, fully well-balanced, unconditionally
  positivity-preserving finite volume framework for flood simulations},
  (2024).
\newblock arXiv preprint arXiv:2402.12248.

\bibitem{CTR}
{\sc M.~Ciallella, D.~Torlo, and M.~Ricchiuto}, {\em Arbitrary high order
  {WENO} finite volume scheme with flux globalization for moving equilibria
  preservation}, J. Sci. Comput., 96 (2023).
\newblock Paper No. 53, 28 pp.

\bibitem{CDL}
{\sc S.~Clain, S.~Diot, and R.~Loub\`{e}re}, {\em A high-order finite volume
  method for systems of conservation laws---multi-dimensional optimal order
  detection ({MOOD})}, J. Comput. Phys., 230 (2011), pp.~4028--4050.

\bibitem{DLM}
{\sc G.~Dal~Maso, P.~G. Lefloch, and F.~Murat}, {\em Definition and weak
  stability of nonconservative products}, J. Math. Pures Appl., 74 (1995),
  pp.~483--548.

\bibitem{DM_FWBC}
{\sc V.~Desveaux and A.~Masset}, {\em A fully well-balanced scheme for shallow
  water equations with {C}oriolis force}, Commun. Math. Sci., 20 (2022),
  pp.~1875--1900.

\bibitem{EPD_08}
{\sc A.~Ern, S.~Piperno, and K.~Djadel}, {\em A well-balanced {R}unge-{K}utta
  discontinuous {G}alerkin method for the shallow-water equations with flooding
  and drying}, International journal for numerical methods in fluids, 58
  (2008), pp.~1--25.

\bibitem{FMT11}
{\sc U.~S. Fjordholm, S.~Mishra, and E.~Tadmor}, {\em Well-balanced and energy
  stable schemes for the shallow water equations with discontinuous
  topography}, J. Comput. Phys., 230 (2011), pp.~5587--5609.

\bibitem{GPC}
{\sc J.~Gallardo, C.~Par\'{e}s, and M.~Castro}, {\em On a well-balanced
  high-order finite volume scheme for shallow water equations with topography
  and dry areas}, J. Comput. Phys., 227 (2007), pp.~574--601.

\bibitem{Bueno}
{\sc I.~G{\'o}mez-Bueno, M.~J. D{\'\i}az~Castro, C.~Par{\'e}s, and G.~Russo},
  {\em Collocation methods for high-order well-balanced methods for systems of
  balance laws}, Mathematics, 9 (2021), p.~1799.

\bibitem{GKS}
{\sc S.~Gottlieb, D.~Ketcheson, and C.-W. Shu}, {\em Strong stability
  preserving {R}unge-{K}utta and multistep time discretizations}, World
  Scientific Publishing Co. Pte. Ltd., Hackensack, NJ, 2011.

\bibitem{GST}
{\sc S.~Gottlieb, C.-W. Shu, and E.~Tadmor}, {\em Strong stability-preserving
  high-order time discretization methods}, SIAM Rev., 43 (2001), pp.~89--112.

\bibitem{Guo2023}
{\sc W.~Guo, Z.~Chen, S.~Qian, G.~Li, and Q.~Niu}, {\em A new well-balanced
  finite volume {CWENO} scheme for shallow water equations over bottom
  topography}, Advances in Applied Mathematics and Mechanics, 15 (2023),
  pp.~1515--1539.

\bibitem{Dumbser19}
{\sc M.~Ioriatti and M.~Dumbser}, {\em A posteriori sub-cell finite volume
  limiting of staggered semi-implicit discontinuous {G}alerkin schemes for the
  shallow water equations}, Applied Numerical Mathematics, 135 (2019),
  pp.~443--480.

\bibitem{JW}
{\sc S.~Jin and X.~Wen}, {\em Two interface-type numerical methods for
  computing hyperbolic systems with geometrical source terms having
  concentrations}, SIAM J. Sci. Comput., 26 (2005), pp.~2079--2101.
\newblock electronic.

\bibitem{KL}
{\sc G.~Kesserwani and Q.~Liang}, {\em Well-balanced {RKDG}2 solutions to the
  shallow water equations over irregular domains with wetting and drying},
  Computers \& Fluids, 39 (2010), pp.~2040--2050.

\bibitem{KKLZ}
{\sc C.~Klingenberg, A.~Kurganov, Y.~Liu, and M.~Zenk}, {\em Moving-water
  equilibria preserving {HLL}-type schemes forthe shallow water equations},
  Commun. Math. Res., 36 (2020), pp.~247--271.

\bibitem{KPshw}
{\sc A.~Kurganov and G.~Petrova}, {\em A second-order well-balanced positivity
  preserving central-upwind scheme for the {S}aint-{V}enant system}, Commun.
  Math. Sci., 5 (2007), pp.~133--160.

\bibitem{LPG2002}
{\sc P.~G. LeFloch}, {\em Hyperbolic systems of conservation laws}, in The
  theory of classical and nonclassical shock waves, Lectures in Mathematics ETH
  Z\"urich, Birkh\"{a}user Verlag, Basel, 2002.

\bibitem{LPG2004}
{\sc P.~G. LeFloch}, {\em Graph solutions of nonlinear hyperbolic systems}, J.
  Hyperbolic Differ. Equ., 1 (2004), pp.~643--689.

\bibitem{LeV98}
{\sc R.~J. LeVeque}, {\em Balancing source terms and flux gradients in
  high-resolution {G}odunov methods: the quasi-steady wave-propagation
  algorithm}, J. Comput. Phys., 146 (1998), pp.~346--365.

\bibitem{LCJKY}
{\sc X.~Liu, X.~Chen, S.~Jin, A.~Kurganov, and H.~Yu}, {\em Moving-water
  equilibria preserving partial relaxation scheme for the {S}aint-{V}enant
  system}, SIAM J. Sci. Comput., 42 (2020), pp.~A2206--A2229.

\bibitem{LLTX}
{\sc Y.~Liu, J.~Lu, Q.~Tao, and Y.~Xia}, {\em An oscillation-free
  {D}iscontinuous {G}alerkin method for shallow water equations}, J. Sci.
  Comput., 92 (2022).
\newblock Paper No. 109, 24 pp.

\bibitem{LNK_2007}
{\sc M.~Luk{\'a}{\v{c}}ov{\'a}-Medvid’ov{\'a}, S.~Noelle, and M.~Kraft}, {\em
  Well-balanced finite volume evolution {G}alerkin methods for the shallow
  water equations}, J. Comput. Phys., 221 (2007), pp.~122--147.

\bibitem{MBerthon24}
{\sc L.~Martaud and C.~Berthon}, {\em ({F}ully) well-balanced entropy stable
  {G}odunov numerical schemes for the shallow water equations},  (2024).
\newblock https://hal.science/hal-04394378v2.

\bibitem{MRA_WB}
{\sc L.~Micalizzi, M.~Ricchiuto, and R.~Abgrall}, {\em Novel well-balanced
  continuous interior penalty stabilizations},  (2023).
\newblock arXiv preprint arXiv:2307.09697.

\bibitem{Victor2016}
{\sc V.~Michel-Dansac, C.~Berthon, S.~Clain, and F.~Foucher}, {\em A
  well-balanced scheme for the shallow-water equations with topography},
  Comput. Math. Appl., 72 (2016), pp.~568--593.

\bibitem{Victor2017}
{\sc V.~Michel-Dansac, C.~Berthon, S.~Clain, and F.~Foucher}, {\em A
  well-balanced scheme for the shallow-water equations with topography or
  {M}anning friction}, J. Comput. Phys., 335 (2017), pp.~115--154.

\bibitem{NXS}
{\sc S.~Noelle, Y.~Xing, and C.-W. Shu}, {\em High-order well-balanced finite
  volume {WENO} schemes for shallow water equation with moving water}, J.
  Comput. Phys., 226 (2007), pp.~29--58.

\bibitem{Puppo16}
{\sc G.~Puppo and M.~Semplice}, {\em Well-balanced high order 1{D} schemes on
  non-uniform grids and entropy residuals}, J. Sci. Comput., 66 (2016),
  pp.~1052--1076.

\bibitem{Ric15}
{\sc M.~Ricchiuto}, {\em An explicit residual based approach for shallow water
  flows}, J. Comput. Phys., 280 (2015), pp.~306--344.

\bibitem{RB09}
{\sc M.~Ricchiuto and A.~Bollermann}, {\em Stabilized residual distribution for
  shallow water simulations}, J. Comput. Phys., 228 (2009), pp.~1071--1115.

\bibitem{ToumaK}
{\sc R.~Touma and S.~Khankan}, {\em Well-balanced unstaggered central schemes
  for one and two-dimensional shallow water equation systems}, Appl. Math.
  Comput., 218 (2012), pp.~5948--5960.

\bibitem{TMK}
{\sc R.~Touma, E.~Malaeb, and C.~Klingenberg}, {\em Combining a central scheme
  with the subtraction method for shallow water equations}, Int. J. Comput.
  Methods,  (2023), p.~2350035.

\bibitem{Vilar}
{\sc F.~Vilar}, {\em A posterior correction of high-order discontinuous
  {G}alerkin scheme through subcell finite volume formulation and flux
  reconstruction}, J. Comput. Phys., 387 (2019), pp.~245--279.

\bibitem{Xing14}
{\sc Y.~Xing}, {\em Exactly well-balanced discontinuous {G}alerkin methods for
  the shallow water equations with moving water equilibrium}, J. Comput. Phys.,
  257 (2014), pp.~536--553.

\bibitem{XS14}
{\sc Y.~Xing and C.-W. Shu}, {\em A survey of high order schemes for the
  shallow water equations}, J. Math. Study, 47 (2014), pp.~221--249.

\bibitem{XZS}
{\sc Y.~Xing, X.~Zhang, and C.-W. Shu}, {\em Positivity-preserving high order
  well-balanced discontinuous {G}alerkin methods for the shallow water
  equations}, Adv. Water Resour., 33 (2010), pp.~476--1493.

\bibitem{XS_AWENO}
{\sc Z.~Xu and C.-W. Shu}, {\em A high-order well-balanced alternative finite
  difference {WENO} ({A}-{WENO}) method with the exact conservation property
  for systems of hyperbolic balance laws},  (2024).
\newblock
  https://bpb-us-w2.wpmucdn.com/sites.brown.edu/dist/0/348/files/2023/07/A-high-order-well-balanced-alternative-finite-difference-WENO-A-WENO-method-with-the-exact-conservation-property.pdf.

\bibitem{XS_DGGL}
{\sc Z.~Xu and C.-W. Shu}, {\em A high-order well-balanced discontinuous
  {G}alerkin method for hyperbolic balance laws based on the {G}auss-{L}obatto
  quadrature rules},  (2024).
\newblock
  https://bpb-us-w2.wpmucdn.com/sites.brown.edu/dist/0/348/files/2023/09/A-high-rder-well-balanced-discontinuous-Galerkin-method-for-hyperbolic-balance-laws.pdf.

\bibitem{Yan2023}
{\sc R.~Yan, W.~Tong, and G.~Chen}, {\em A mass conservative, well balanced and
  positivity-preserving central scheme for shallow water equations}, Appl.
  Math. Comput., 443 (2023), p.~127778.

\bibitem{Zhang2023}
{\sc J.~Zhang, Y.~Xia, and Y.~Xu}, {\em Moving water equilibria preserving
  discontinuous {G}alerkin method for the shallow water equations}, J. Sci.
  Comput., 95 (2023), p.~48.

\bibitem{Zhao2023}
{\sc Z.~Zhao and M.~Zhang}, {\em Well-balanced fifth-order finite difference
  {H}ermite {WENO} scheme for the shallow water equations}, J. Comput. Phys.,
  475 (2023), p.~111860.

\end{thebibliography}
\end{document}